\documentclass[onefignum,onetabnum]{siamart190516}

\usepackage{lipsum}
\usepackage{amsfonts}
\usepackage{graphicx}
\usepackage{epstopdf}
\usepackage{algorithm}
\usepackage[noend]{algpseudocode}
\usepackage{multicol}
\usepackage{siunitx}
\usepackage{subcaption}
\ifpdf
  \DeclareGraphicsExtensions{.eps,.pdf,.png,.jpg}
\else
  \DeclareGraphicsExtensions{.eps}
\fi
\usepackage{tikz}
\usetikzlibrary{calc}
\usetikzlibrary{quotes}
\usetikzlibrary{shapes}
\usetikzlibrary{external}
\tikzstyle{vertexstyle}=[thick, draw, rectangle, fill opacity=0.4, text opacity=1]
\tikzstyle{edgestyle}=[thick, color=lightgray]

\usepackage{bm}
\newcommand{\vct}[1]{\bm{\mathsf{#1}}}
\newcommand{\mtx}[1]{\bm{\mathsf{#1}}}
\newcommand{\mA}{\mtx{A}}
\newcommand{\mB}{\mtx{B}}
\newcommand{\mG}{\mtx{G}}
\newcommand{\mI}{\mtx{I}}
\newcommand{\mP}{\mtx{P}}
\newcommand{\mQ}{\mtx{Q}}
\newcommand{\mR}{\mtx{R}}
\newcommand{\mS}{\mtx{\Sigma}}
\newcommand{\mU}{\mtx{U}}
\newcommand{\mV}{\mtx{V}}
\newcommand{\mW}{\mtx{W}}
\newcommand{\mY}{\mtx{Y}}
\newcommand{\mZ}{\mtx{Z}}
\newcommand{\mOmega}{\mtx{\Omega}}
\newcommand{\mPsi}{\mtx{\Psi}}
\newcommand{\mzero}{\mtx{0}}
\newcommand{\ushort}{\mU} %
\newcommand{\vshort}{\mV} %
\newcommand{\ulong}{\mtx{\mathcal{U}}} %
\newcommand{\vlong}{\mtx{\mathcal{V}}} %
\newcommand{\qq}{\vct{q}}
\newcommand{\uu}{\vct{u}}
\newcommand{\vv}{\vct{v}}
\newcommand{\xx}{\vct{x}}

\newcommand{\domain}{X}
\newcommand{\R}{\mathbb{R}}                                   %
\newcommand{\Z}{\mathbb{Z}}                                   %
\newcommand{\h}{\mathcal{H}}                                  %
\newcommand{\kfunc}{\mathcal{K}}                              %
\newcommand{\alphanear}{\mathcal{L}_\alpha^\text{nei}}        %
\newcommand{\alphainteraction}{\mathcal{L}_\alpha^\text{int}} %
\newcommand{\betainteraction}{\mathcal{L}_\beta^\text{int}}   %
\newcommand{\taunear}{\mathcal{L}_\tau^\text{nei}}            %
\newcommand{\tauinteraction}{\mathcal{L}_\tau^\text{int}}     %
\newcommand{\svd}{\texttt{svd}}                               %
\newcommand{\qr}{\texttt{qr}}                                 %
\newcommand{\fref}[1]{Figure~\ref{#1}}                        %

\newcommand{\bigO}{\mathcal{O}}
\newcommand{\graphdegree}{\deg(G)}
\newcommand{\tl}{T_l} %
\newcommand{\tAl}{T_{A^{(l)}}} %
\newcommand{\tmult}{T_\text{mult}} %
\newcommand{\tflop}{T_\text{flop}}
\newcommand{\tcompress}{T_\text{compress}} %
\newcommand{\tcolor}{T_\text{color}} %
\newcommand{\colorsnonunif}{\chi_\text{nonunif}} %
\newcommand{\colorsnear}{\chi_\text{leaf}} %
\newcommand{\colorsunif}{\chi_\text{unif}} %

\def\lrcolor{white!40!blue} %
\def\frcolor{white!30!red} %

\newsiamremark{remark}{Remark}
\newsiamremark{hypothesis}{Hypothesis}
\crefname{hypothesis}{Hypothesis}{Hypotheses}
\newsiamthm{claim}{Claim}

\headers{Rank-structured Matrix Compression with Graph Coloring}{J. Levitt and P. G. Martinsson}

\title{
    Randomized Compression of Rank-Structured Matrices Accelerated with Graph Coloring
    \thanks{Submitted to the editors DATE.
}}

\author{
    James Levitt\thanks{
        Oden Institute for Computational Engineering and Sciences,
        University of Texas at Austin, Austin, TX
        (\email{jlevitt@utexas.edu}, \email{pgm@oden.utexas.edu}).
    }
    \and
    Per-Gunnar Martinsson\footnotemark[2]
}

\usepackage{amsopn}

 \newcommand{\pvct}[1]{#1}

\ifpdf
\hypersetup{
    pdftitle={Randomized Compression of Rank-Structured Matrices
    Accelerated with Graph Coloring},
    pdfauthor={J. Levitt and P. G. Martinsson}
}
\fi

\begin{document}

\maketitle

\begin{abstract}
A randomized algorithm for computing a data sparse representation of a given
rank-structured matrix $\mtx{A}$ (a.k.a.~an $\mathcal{H}$-matrix) is presented. 
The algorithm draws on the randomized singular value decomposition (RSVD), and
operates under the assumption that methods for rapidly applying $\mtx{A}$
and $\mtx{A}^{*}$ to vectors are available. 
The algorithm uses graph coloring algorithms to analyze the hierarchical tree that defines the rank structure to generate a tailored probability distribution from which to draw the random test matrices. 
The matrix is then applied to the test matrices, and in a final step the matrix
itself is reconstructed by the observed input-output pairs.
The method presented is an evolution of the ``peeling algorithm'' of 
\textit{L. Lin, J. Lu, and L. Ying, ``Fast construction of hierarchical matrix representation from matrix–vector multiplication,'' JCP, \textbf{230}(10), 2011}.
For the case of uniform trees, the new method substantially reduces the pre-factor
of the original peeling algorithm.
More significantly, the new technique leads to dramatic acceleration for many
non-uniform trees since it constructs test matrices that are optimized for a
given tree. 
The algorithm is particularly effective for kernel matrices involving a set of
points restricted to a lower dimensional object than the ambient space, such as
a boundary integral equation defined on a surface in three dimensions.
\end{abstract}



\section{Introduction} \label{sec:intro}
\begin{figure}[tb]
  \centering
  \includegraphics[scale=.2]{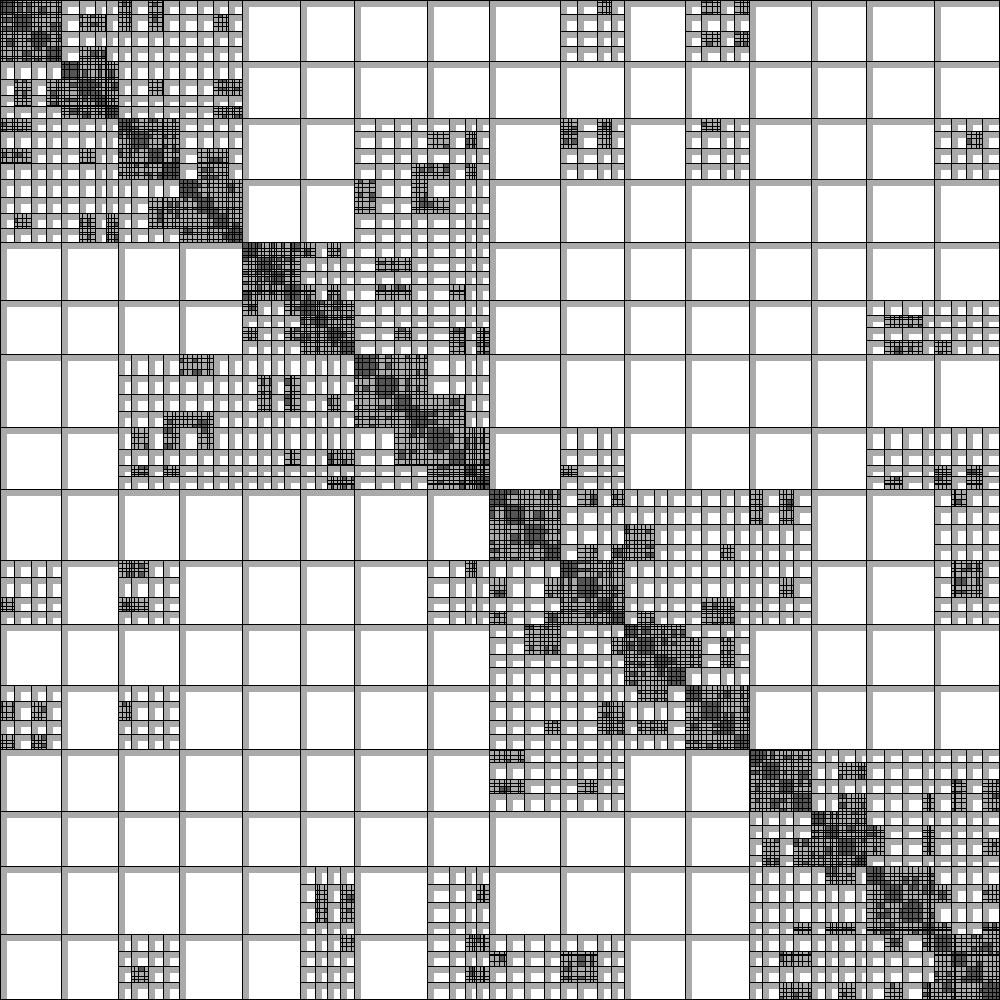}
  \caption{
      An $\h^1$ matrix for a quadtree over a uniform grid in the plane.
      Dense blocks are shown in dark gray,
      and low-rank blocks are represented with a white background
      and light gray rectangles representing the shapes of the low-rank factors.
    }
  \label{fig:h1matrix}
\end{figure}

We present a set of efficient algorithms for computing
data-sparse representations of large dense matrices that have 
\textit{rank structure}. This property is often defined to mean that
an $N\times N$ matrix can be tessellated  into $O(N)$ blocks in such 
a way that each block is either small or of low numerical rank,
cf.~Figures \ref{fig:h1matrix} and \ref{fig:matrix_strong_1d_lvl3}.
Matrices with such structure can be stored and applied to vectors efficiently, 
often with cost that scales linearly  or close to linearly with $N$. 
It is in many environments also possible to compute an approximate 
inverse or LU factorization in linear or close to linear time. 
The notion of rank-structured matrices appears implicitly in 
techniques such as Barnes-Hut \cite{1986_barnes_hut} and the 
Fast Multipole Method \cite{rokhlin1987}.
The idea was then
brought to the forefront as a linear algebraic structure by 
Hackbusch and co-workers \cite{hackbusch1999sparse,2002_hackbusch_H2},
and is now a subject of research under names such as 
$\mathcal{H}$-matrices~\cite{bebendorf2008hierarchical,2010_borm_book,hackbusch1999sparse};
HODLR matrices~\cite{ambikasaran2013mathcal,martinsson2009fast},
Hierarchically Semi-Separable (HSS) a.k.a.~Hierarchically Block-Separable (HBS) matrices~\cite{chandrasekaran2007fast,chandrasekaran2006fast,gillman2012direct,xia2010superfast},
Recursive Skeletonization~\cite{ho2012fast,2005_martinsson_fastdirect,minden2017recursive},
and many more.

The question we address is the following: Suppose that you are given
an $N\times N$ matrix $\mtx{A}$ that you know is rank structured under some given partitioning, 
but you do not have direct access to the low-rank factors that define 
the compressible off-diagonal blocks. Instead, you have access to some
efficient procedure that, given tall thin matrices 
$\mtx{\Omega},\mtx{\Psi} \in \mathbb{R}^{N\times \ell}$, 
can evaluate the matrix-matrix products
$$
\mtx{Y} = \mtx{A\Omega},
\qquad\mbox{and}\qquad
\mtx{Z} = \mtx{A}^{*}\mtx{\Psi}.
$$
Your task is to build two random (but structured) matrices 
$\mtx{\Omega}$ and $\mtx{\Psi}$ that have the property that
$\mtx{A}$ can be recovered from the information in the set
$\{\mtx{Y},\,\mtx{\Omega},\,\mtx{Z},\,\mtx{\Psi}\}$ in a 
computationally efficient manner.
The algorithms described here solve the reconstruction problem using
$\ell \sim k\,\log(N)$ sample vectors, where $k$ is an upper bound on 
the ranks of the off-diagonal blocks. 
The key novelty is the formulation of a graph coloring problem to analyze the
cluster tree that defines the rank structure to build test matrices that
are optimized for the specific problem under consideration. (We assume
that the correct partition is already known, typically from the geometry
of an underlying physical problem.)

The method that we describe relies on the user having access to an efficient
algorithm for applying the rank-structured matrix. An important application
of such a method concerns the compression of the dense Schur complements
that arise in the LU factorization of sparse matrices associated with finite
difference and finite element discretizations of elliptic PDEs.
\cite{xia2012superfast,2019_amestoy_BLR_multifrontal,2016_ghysels_randomized_multifrontal,gillman2012direct}. During the LU factorization, the user needs to form
a data-sparse representation of a Schur complement such as 
$\mtx{S}_{22} = \mtx{A}_{21}\mtx{A}_{11}^{-1}\mtx{A}_{12}$,
representing the elimination of the diagonal block $\mtx{A}_{11}$.
Commonly, the matrices $\mtx{A}_{21}$,\,$\mtx{A}_{11}$,\,$\mtx{A}_{12}$ 
are already held in a sparse or data-sparse form, which means that
the map $\vct{x} \mapsto \mtx{A}_{21}\mtx{A}_{11}^{-1}\mtx{A}_{12}\vct{x}$ can be evaluated rapidly.
In contrast, explicitly forming $\mtx{A}_{11}^{-1}$ and the matrix products
is much more challenging. Another application of our technique concerns
the multiplication of dense operators in potential theory, where each
factor operator can be applied using a legacy technique such as the FMM
\cite{rokhlin1987}, and where inverses can often be applied via ``fast
direct solvers'' \cite{2019_martinsson_fast_direct_solvers}.

The method we describe is inspired by the ``peeling algorithm'' of \cite{lin2011fast}, 
which to the best of our knowledge was the first true black-box algorithm described
in the literature. The method of \cite{lin2011fast} has the same asymptotic sample
complexity $\ell \sim k\,\log(N)$ as our method, but involves substantially larger
pre-factors. To be precise, \cite{lin2011fast} is targeted specifically for 
$\mathcal{H}^1$- and $\mathcal{H}^2$-matrices arising from the discretization of
integral equations. Strong admissibility, and regular tree structures are used.
In this environment, the method requires $\ell \sim k\,8^d\,\log(N)$
matrix-vector products involving $\mA$ and $\mA^*$,
where $d$ is the dimension of space in which the underlying integral equation
is defined. In contrast, the
method presented here has complexity $\ell \sim k\,6^d\,\log(N)$ for fully
populated uniform trees. For more general trees, the acceleration over the method of 
\cite{lin2011fast} can be even more dramatic, since the adaptivity of our method
enables it to ``discover'' intrinsic low dimensional structure in a given problem.
As an illustration, Section
\ref{sec:experiments} reports on experiments involving a boundary integral equation
defined on a 2D surface in three dimensional space, as well as examples in higher
dimensions.

Another advantage of the techniques presented here
is that they are not limited to the $\h^1$ structure.
To compress uniform $\h^1$ and $\h^2$ matrices,
the presented algorithms obtain uniform basis matrices
(cf.~\cite[Sec.~4.6]{hackbusch1999sparse})
by sampling the interactions of a box
with its entire interaction list collectively.
This process results in higher quality samples
while requiring fewer matrix-vector products
compared to existing methods
that approximate interactions between boxes separately
and then apply a recompression step to obtain uniform basis matrices.
More generally, the formulation of a graph coloring problem
can be used to design test matrices for any tessellation
(e.g., arising from some other geometric or algebraic admissibility condition).

\begin{remark}[Related work]
A kindred class of algorithms for computing a rank-structured matrix
by observing its action on random vectors was 
described in \cite{2008_martinsson_randomhudson,martinsson2011fast}.
These techniques have actual
linear complexity (no logarithmic terms), and tend to be very fast
in practice. However, they are not true black-box algorithms, as they
require the direct evaluation of a small number of entries of the
matrix. In contrast, the method presented here is truly black box,
like the methods of \cite{lin2011fast,martinsson2016compressing}.

The recent work~\cite{levitt2022linear}
describes an algorithm for compressing
Hierarchically Block Separable (HBS) matrices
(also known as Hierarchically Semiseparable or HSS matrices),
which is truly black box and has linear complexity.
However, that work relies on the assumption
that the basis matrices are nested,
and it only considers the HBS format,
which is subject to a weak admissibility condition.
The present work focuses on rank-structured matrices
subject to strong admissibility,
including the $\h^1$ and uniform $\h^1$ formats,
which do not assume nested basis matrices.

The paper \cite{schafer2021sparse} explores the
problem of sparse direct solvers for discretized elliptic operators
using techniques related to those presented here. Specifically, 
\cite{schafer2021sparse} relies on graph coloring to analyze
the sparsity graph associated with a discretized PDE, and uses a
``peeling algorithm'' to build an approximate, sparsity preserving
factorization.
\end{remark}

The manuscript is structured as follows:
Section \ref{sec:preliminaries} surveys some basic linear algebraic techniques that
we rely on.
Section  \ref{sec:rank_structured} introduces our formalism for 
rank-structured matrices.
Section \ref{sec:compress} describes the new algorithm, and analyzes its
asymptotic complexity.
Section \ref{sec:experiments} describes numerical results.

\section{Preliminaries} \label{sec:preliminaries}

\subsection{Notation}

Throughout the paper, we measure a vector $\xx \in \R^n$
by its Euclidean norm
$\Vert \xx \Vert = \left( \sum_i \vert x_i \vert^2 \right)^\frac{1}{2}$.
We measure a matrix $\mA \in \R^{m \times n}$
with the corresponding operator norm
$\Vert \mA \Vert = \sup_{\Vert \xx \Vert = 1} \Vert \mA \xx \Vert$,
and in some cases with the Frobenius norm
$\Vert \mA \Vert_\text{Fro}
= (\sum_{i, j} \vert \mA(i, j) \vert^2)^{1/2}$.
To denote submatrices,
we use the notation of Golub and Van Loan~\cite{golub2013matrix}:
If $\mA$ is an $m \times n$ matrix,
and $I = \lbrack i_1, i_2, \dots, i_k \rbrack$
and $J = \lbrack j_1, j_2, \dots, j_l \rbrack$,
then $\mA(I, J)$ denotes the $k \times l$ 
submatrix of $\mA$ at the intersection of the
rows in $I$ and the columns in $J$.
We let $\mA(I,:)$ denote the submatrix
$\mA(I, \lbrack 1, 2, \dots, n \rbrack)$
and define $\mA(:, J)$ analogously.
We define the complement $I^c$ of row index set $I$
to be the set of row indices of $\mA$ not included in $I$,
$I^c = \{ 1, 2, \dots, m \} \setminus I$,
and we define the complement of a column index set analogously.
We let $\mA^*$ denote the transpose of $\mA$,
and we say that matrix $\mU$ is \emph{orthonormal}
if its columns are orthonormal, $\mU^* \mU = \mI$.

\subsection{The QR factorization}
The column-pivoted \emph{QR factorization}
of a matrix $\mA$ of size $m \times n$
takes the form
\begin{equation}
\label{eq:QR}
\begin{array}{ccccccc}
    \mA        & \mP        &=& \mQ        & \mR, \\
    m \times n & n \times n &       & m \times r & r \times n
\end{array}
\end{equation}
where
$r = \min(m, n)$,
$\mQ$ is orthonormal,
$\mR$ is upper-triangular,
and
$\mP$ is a permutation matrix.
Representing the permutation matrix $\mP$
as the vector $J \subset \Z_+^n$
of column indices such that
$\mP = \mI(:, J)$,
the factorization~\cref{eq:QR}
can be expressed as
\begin{equation*}
\begin{array}{ccccccc}
    \mA(:, J)  & &=& \mQ        & \mR. \\
    m \times n & &       & m \times r & r \times n
\end{array}
\end{equation*}
For a matrix that is numerically low-rank,
a rank-$k$ approximation of $\mA$
is given by a ``partial QR factorization of $\mA$,''
\begin{equation*}
\begin{array}{ccccccc}
    \mA(:, J)  & &\approx& \mQ_k      & \mR_k. \\
    m \times n & &       & m \times k & k \times n
\end{array}
\end{equation*}

\subsection{The singular value decomposition (SVD)}

The singular value decomposition
of a matrix $\mA$ of size $m \times n$
takes the form
\begin{equation}
\label{eq:SVD}
\begin{array}{ccccccc}
    \mA        &=& \mU        & \mS        & \mV^*,\\
    m \times n &         & m \times r & r \times r & r \times n
\end{array}
\end{equation}
where
$r = \min(m, n)$,
$\mU$ and $\mV$ are orthonormal matrices,
and $\mS$ is a diagonal matrix
with diagonal elements $\{\sigma_j\}_{j=1}^r$
ordered such that
$\sigma_1 \geq \sigma_2 \geq \dots \geq \sigma_r \geq 0$.
The columns of $\mU$ and $\mV$,
denoted by
$\{\uu_i\}_{i=1}^r$
and
$\{\vv_i\}_{i=1}^r$,
are the left and right singular vectors of $\mA$.
We let $\mA_k$ denote the rank-$k$ approximation
obtained by truncating the SVD to its first $k$ terms,
so that $\mA_k = \sum_{j=1}^k \sigma_j \uu_j \vv_j^*$.
It then follows that
\begin{equation*}
    \Vert \mA - \mA_k \Vert = \sigma_{k + 1}
    \qquad
    \text{and that}
    \qquad
    \Vert \mA - \mA_k \Vert_\text{Fro}
    =
    \left( \sum_{j = k + 1}^{\min(m, n)} \sigma_j^2 \right)^{1/2}.
\end{equation*}
The Eckart--Young theorem asserts that $\mA_k$
achieves the smallest possible approximation error of $\mA$,
in both the operator and Frobenius norms,
of any rank-$k$ matrix.

\subsection{Functions for low-rank factorizations}

We introduce the following notation to denote calls to functions
that return the results of
QR factorizations and singular value decompositions.
Function calls to evaluate the full factorizations are written as
$$
\lbrack \mQ, \mR, J \rbrack = \qr(\mA),
\qquad
\lbrack \mU, \mS, \mV \rbrack = \svd(\mA).
$$
Calls to evaluate the rank-$k$ truncated factorizations are written as
$$
\lbrack \mQ, \mR, J \rbrack = \qr(\mA, k),
\qquad
\lbrack \mU, \mS, \mV \rbrack = \svd(\mA, k).
$$
We write
$$
\lbrack \mQ, \mR \rbrack = \qr(\mA),
\qquad
\mQ = \qr(\mA)
$$
to compute an unpivoted QR factorization,
and just the factor $\mQ$, respectively.

\subsection{Randomized compression}

In this section, we give a brief review of randomized low-rank approximation,
following the presentation of~\cite{martinsson2016compressing}.
Let $\mA$ be an $m \times n$ matrix
that can be accurately approximated by a matrix of rank $k$,
and suppose we seek to determine a matrix $\mQ$
with orthonormal columns (as few as possible) such that
$\Vert \mA - \mQ \mQ^* \mA \Vert$
is small.
In other words, we seek a matrix $\mQ$
whose columns form an approximate orthornomal basis (ON-basis)
for the column space of $\mA$.
This task can efficiently be solved via the following randomized procedure:
\begin{enumerate}
    \item Pick a small integer $p$ representing how much ``oversampling'' is done.
        ($p = 10$ is often good.)
    \item Form an $n \times (k + p)$ matrix $\mG$
        whose entries are drawn independently from a
        standardized normal (``Gaussian'') distribution.
    \item Form the ``sample matrix'' $\mY = \mA \mG$ of size $m \times (k + p)$.
    \item Construct an $m \times (k + p)$ matrix $\mQ$ whose columns form an ON basis
        for the columns of $\mY$.
\end{enumerate}

Note that each column of the sample matrix $\mY$
is a random linear combination of the columns of $\mA$.
The probability of the algorithm
producing an accurate result
approaches 1 extremely rapidly as $p$ increases,
and, remarkably, this probability depends only on $p$
(not on $m$ or $n$, or any other properties of $\mA$);
cf.~\cite{halko2011finding}.

Now, suppose we would like to use the above sampling procedure
to compute a low-rank factorization of the form
\begin{equation}
\label{eq:UBV}
\begin{array}{ccccccc}
    \mA        &= & \mU        & \mB        & \mV,\\
    m \times n &  & m \times r & r \times r & r \times n
\end{array}
\end{equation}
where $\mU$ and $\mV$ have orthonormal columns
and $\mB$ is a square matrix, not necessarily diagonal.
We review two methods of completing this task
using randomized sampling.
The algorithm summarized in \cref{alg:compress_double_random}
involves two randomized samples,
one of $\mA$ and one of $\mA^*$.
It is based on so-called ``single-view'' algorithms
for matrix compression~\cite{tropp2017practical,martinsson2020randomized}.
The second method~\cite{halko2011finding},
summarized in \cref{alg:compress_two_stage},
uses a randomized sample of $\mA$
to compute the column basis matrix $\mU$,
and then operates on the product $\mA^* \mU$
to obtain $\mB$ and $\mV$.
We will use both methods
in the algorithms described in \cref{sec:compress}.
(Though it is sometimes prudent to do slightly
more over-sampling in the ``single view'' environment
\cite{tropp2017practical,martinsson2020randomized},
we use the same $p$ in both cases.)

\begin{algorithm}
    \caption{Compress $\mA \approx \mU \mB \mV^*$ via two randomized samples}
\label{alg:compress_double_random}
\begin{algorithmic}
    \State{Form $n \times (k + p)$ Gaussian random matrix $\mG_1$.}
    \State{Multiply $\mY = \mA \mG_1$.}
    \State{Orthonormalize
        $\mU = \qr(\mY, k)$.}
    \State{Form $m \times (k + p)$ Gaussian random matrix $\mG_2$.}
    \State{Multiply $\mZ = \mA^* \mG_2$.}
    \State{Orthonormalize
        $\mV = \qr(\mZ, k)$.}
    \State{Solve $\mB = (\mG_2^* \mU)^\dagger \mG_2^* \mA \mG_1 (\mV^* \mG_1)^\dagger$.}
    \State \Return{$\mU, \mB, \mV$.}
\end{algorithmic}
\end{algorithm}

\begin{algorithm}
    \caption{Compress $\mA \approx \mU \mB \mV^*$ with one randomized and one deterministic sample}
\label{alg:compress_two_stage}
\begin{algorithmic}
    \State{Form an $n \times (k + p)$ Gaussian random matrix $\mG$.}
    \State{Multiply $\mY = \mA \mG$.}
    \State{Orthonormalize $\mU = \qr(\mY)$.}
    \State{Multiply $\mW = \mA^* \mU$.}
    \State{Orthonormalize $\lbrack \mV, \mB^* \rbrack = \qr(\mW)$.}
    \State \Return{$\mU, \mB, \mV$.}
\end{algorithmic}
\end{algorithm}

\subsection{The degree of saturation graph coloring algorithm}

A vertex coloring of a graph
is an assigment of a color to each vertex
in such a way that no pair of adjacent vertices shares the same color.
The problem of finding a vertex coloring
with the minimum number of colors is NP-hard,
but there exist a number of algorithms
for efficiently coloring graphs,
though they may produce colorings with more than the minimum number of colors.

Greedy graph coloring algorithms process the vertices in sequence,
at each iteration assigning to one vertex the first available color
that hasn't already been assigned to one of its neighbors.
In the \emph{degree of saturation} (DSatur) algorithm~\cite{brelaz1979new},
the choice of which vertex to color at each step
is made by selecting from the remaining uncolored vertices
the one whose neighbors have the greatest number of distinct colors
(the so-called degree of saturation).
The procedure is summarized in \cref{alg:dsatur}.

\begin{algorithm}
    \caption{Greedy graph coloring with DSatur}
\label{alg:dsatur}
\begin{algorithmic}
    \State{Initialize priority queue $q$ of vertices keyed by degree of saturation (initially all zero)}
    \State{Initialize for each vertex a set of invalid colors (initially all empty)}
    \While{$q$ is not empty}
        \State{Pop from $q$ the vertex $v$ with highest degree of saturation}
        \Comment{$\bigO(\log{\vert V \vert})$}

        \State{Assign a color to $v$, creating a new one if necessary}
        \Comment{$\bigO(\graphdegree)$}
        \For{each vertex $w$ adjacent to $v$}
            \State{Add the color of $v$ to the set of invalid colors for $w$}
            \Comment{$\bigO(1)$}

            \State{Update the priority of $w$ within $q$}
            \Comment{$\bigO(\log{\vert V \vert)}$}
        \EndFor
    \EndWhile
\end{algorithmic}
\end{algorithm}

While the asymptotic complexity of DSatur
is often described as quadratic in the number of vertices,
the algorithm can be implemented
with quasilinear complexity,
assuming the degree of the graph is bounded.
The \emph{degree of a vertex} is the number of edges incident to that vertex,
and the \emph{degree of a graph} is the maximum degree over all of its vertices.
Summing the costs listed in \cref{alg:dsatur},
we find that the asymptotic complexity is
\begin{equation}
    \label{eq:dsaturcomplexity}
    \tcolor \sim \graphdegree \vert V \vert \log \vert V \vert,
\end{equation}
where $\graphdegree$ denotes the degree of the graph.
For the step of assigning a color to vertex $v$,
we note that a greedy coloring algorithm
uses no more than $\graphdegree + 1$ colors in the worst case,
so checking the existing colors
for whether they belong to the set of invalid colors for $v$
requires $\bigO(\graphdegree)$ operations.

\section{Rank-structured matrices} \label{sec:rank_structured}
This section provides definitions
of the rank-structured matrix formats that we use,
following the notational framework of~\cite{2019_martinsson_fast_direct_solvers}.

Let $\{x_i\}_{i=1}^N \subset \domain$ be a set of points
within a $d$-dimensional hypercube $\domain = [0, 1]^d$,
and let $\mA$ be an $N$-by-$N$ matrix where
\begin{equation*}
    \mA(i, j) = \kfunc(x_i, x_j) \qquad 1 \leq i \leq N, 1 \leq j \leq N,
\end{equation*}
for some kernel function $\kfunc$.
That is, the $(i,j)$ entry of $\mA$ represents
an interaction between points $x_i, x_j$
as defined by $\kfunc$.

We divide the domain into a number of boxes
and arrange the boxes into a tree structure,
where the levels of the tree represent
successively finer partitions of $\domain$.
Level 0 only contains a single box,
which encloses all of $\domain$.
The boxes belonging to level $l+1$
are obtained by bisecting each of the boxes in level $l$
along each dimension to form $2^d$ smaller boxes.
The $2^d$ boxes in level $l+1$ obtained by splitting a box in level $l$
are designated as the children of that box,
giving rise to the tree structure.
Boxes that do not contain any points are omitted from the tree,
so the branching factor of the tree may be less than $2^d$,
depending on the distribution of points.
The splitting procedure is applied recursively
to boxes that contain more than $m$ points,
where $m$ is a prespecified maximum number of points
to be allowed in a leaf box.
We let $L$ denote the depth of the tree,
and assume that the points are distributed evenly enough
across the domain so that $L \sim \log N$.
For each box $\tau$,
we define $I_\tau = \{ i : x_i \in \tau\}$
to be the set of indices
of the points contained by $\tau$.

\begin{figure}[tb]
    \centering
    \minipage{0.1\textwidth}
        Level 0
        \vspace{0.5cm}

        Level 1
        \vspace{0.5cm}

        Level 2
    \endminipage
    \hfill
    \minipage{0.35\textwidth}
        \begin{tikzpicture}[level distance=1.0cm,
            level 1/.style={sibling distance=2.5cm},
            level 2/.style={sibling distance=1.25cm}]
            \node{$I_1$}
            child{
                node{$I_2$}
                child{node{$I_4$}}
                child{node{$I_5$}}
            }
            child{
                node{$I_3$}
                child{node{$I_6$}}
                child{node{$I_7$}}
            };
        \end{tikzpicture}
    \endminipage
    \hfill
    \minipage{0.45\textwidth}
        $I_1 = [1, 2, ..., 400]$
        \vspace{0.5cm}

        $I_2 = [1, ..., 200]$,
        $I_3 = [201, ..., 400]$

        \vspace{0.5cm}
        $I_4 = [1, ..., 100]$,
        $I_5 = [101, ..., 200]$, ...
    \endminipage
    \caption{A binary tree structure, where the levels of the tree represent
    successively refined partitions of the index vector $[1, ..., 400]$.
    }
    \label{fig:tree}
\end{figure}
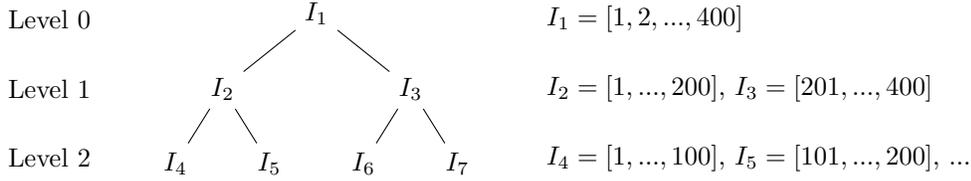

Two boxes that belong to the same level and have overlapping boundaries
are said to be \emph{neighbors}.
The neighbors of box $\tau$,
of which there are up to $3^d$ (including $\tau$ itself),
are stored in the \emph{neighbor list} of $\tau$,
denoted by $\taunear$.
The \emph{interaction list} of $\tau$,
denoted by $\tauinteraction$
contains the children of neighbors of the parent of $\tau$,
excluding those that are neighbors of $\tau$.
The interaction list of a box contains up to $6^d - 3^d$ boxes.
\cref{fig:tree} shows an example tree structure
and index lists.

Matrix $\mA$ is said to have $\h^1$ structure
if block $\mA(I_\alpha, I_\beta)$ is numerically low-rank,
for all pairs of boxes $\alpha, \beta$
belonging to the interaction lists of one another.
Such blocks, and the corresponding pairs of boxes,
are said to be \emph{admissible}.
\cref{fig:matrix_strong_1d_lvl3} shows a tessellation of a matrix
consisting of the admissible blocks
as well as a number of \emph{inadmissible} blocks
of interactions between neighboring boxes in level $L$,
which are not necessarily low-rank.
Constructing a compressed representation of an $\h^1$ matrix
consists of finding low rank approximations
\begin{equation*}
    \mA(I_\alpha, I_\beta) \approx \ulong_{\alpha, \beta} \mB_{\alpha, \beta} \vlong_{\alpha, \beta}
\end{equation*}
for each admissible pair $(\alpha, \beta)$
and storing the inadmissible blocks
corresponding to neighbor interactions of boxes in level $L$.

\begin{figure}[tb]
    \centering
    \def\n{8cm}
    \begin{tikzpicture}[level distance=\n/10,
        level 1/.style={sibling distance=\n/2},
        level 2/.style={sibling distance=\n/4},
        level 3/.style={sibling distance=\n/8}]
        \begin{scope}[local bounding box=scope1, line width=1pt]

            \draw (0, 0) rectangle (\n, \n);

            \draw[fill=\lrcolor] (\n*2/4, \n*3/4) rectangle (\n*3/4, \n*4/4) node[midway] {$\mA_{4, 6}$};
            \draw[fill=\lrcolor] (\n*3/4, \n*3/4) rectangle (\n*4/4, \n*4/4) node[midway] {$\mA_{4, 7}$};
            \draw[fill=\lrcolor] (\n*3/4, \n*2/4) rectangle (\n*4/4, \n*3/4) node[midway] {$\mA_{5, 7}$};
            \draw[fill=\lrcolor] (\n*0/4, \n*1/4) rectangle (\n*1/4, \n*2/4) node[midway] {$\mA_{6, 4}$};
            \draw[fill=\lrcolor] (\n*0/4, \n*0/4) rectangle (\n*1/4, \n*1/4) node[midway] {$\mA_{7, 4}$};
            \draw[fill=\lrcolor] (\n*1/4, \n*0/4) rectangle (\n*2/4, \n*1/4) node[midway] {$\mA_{7, 5}$};

            \draw[fill=\lrcolor] (\n*2/8, \n*7/8) rectangle (\n*3/8, \n*8/8) node[midway] {$\mA_{ 8, 10}$};
            \draw[fill=\lrcolor] (\n*3/8, \n*7/8) rectangle (\n*4/8, \n*8/8) node[midway] {$\mA_{ 8, 11}$};
            \draw[fill=\lrcolor] (\n*3/8, \n*6/8) rectangle (\n*4/8, \n*7/8) node[midway] {$\mA_{ 9, 11}$};
            \draw[fill=\lrcolor] (\n*0/8, \n*5/8) rectangle (\n*1/8, \n*6/8) node[midway] {$\mA_{10,  8}$};
            \draw[fill=\lrcolor] (\n*4/8, \n*5/8) rectangle (\n*5/8, \n*6/8) node[midway] {$\mA_{10, 12}$};
            \draw[fill=\lrcolor] (\n*5/8, \n*5/8) rectangle (\n*6/8, \n*6/8) node[midway] {$\mA_{10, 13}$};
            \draw[fill=\lrcolor] (\n*0/8, \n*4/8) rectangle (\n*1/8, \n*5/8) node[midway] {$\mA_{11,  8}$};
            \draw[fill=\lrcolor] (\n*1/8, \n*4/8) rectangle (\n*2/8, \n*5/8) node[midway] {$\mA_{11,  9}$};
            \draw[fill=\lrcolor] (\n*5/8, \n*4/8) rectangle (\n*6/8, \n*5/8) node[midway] {$\mA_{11, 13}$};
            \draw[fill=\lrcolor] (\n*2/8, \n*3/8) rectangle (\n*3/8, \n*4/8) node[midway] {$\mA_{12, 10}$};
            \draw[fill=\lrcolor] (\n*6/8, \n*3/8) rectangle (\n*7/8, \n*4/8) node[midway] {$\mA_{12, 14}$};
            \draw[fill=\lrcolor] (\n*7/8, \n*3/8) rectangle (\n*8/8, \n*4/8) node[midway] {$\mA_{12, 15}$};
            \draw[fill=\lrcolor] (\n*2/8, \n*2/8) rectangle (\n*3/8, \n*3/8) node[midway] {$\mA_{13, 10}$};
            \draw[fill=\lrcolor] (\n*3/8, \n*2/8) rectangle (\n*4/8, \n*3/8) node[midway] {$\mA_{13, 11}$};
            \draw[fill=\lrcolor] (\n*7/8, \n*2/8) rectangle (\n*8/8, \n*3/8) node[midway] {$\mA_{13, 15}$};
            \draw[fill=\lrcolor] (\n*4/8, \n*1/8) rectangle (\n*5/8, \n*2/8) node[midway] {$\mA_{14, 12}$};
            \draw[fill=\lrcolor] (\n*4/8, \n*0/8) rectangle (\n*5/8, \n*1/8) node[midway] {$\mA_{15, 12}$};
            \draw[fill=\lrcolor] (\n*5/8, \n*0/8) rectangle (\n*6/8, \n*1/8) node[midway] {$\mA_{15, 13}$};

            \draw[fill=\frcolor] (\n*1/8, \n*7/8) rectangle (\n*2/8, \n*8/8) node[midway] {$\mA_{ 8,  9}$};
            \draw[fill=\frcolor] (\n*2/8, \n*6/8) rectangle (\n*3/8, \n*7/8) node[midway] {$\mA_{ 9, 10}$};
            \draw[fill=\frcolor] (\n*3/8, \n*5/8) rectangle (\n*4/8, \n*6/8) node[midway] {$\mA_{10, 11}$};
            \draw[fill=\frcolor] (\n*4/8, \n*4/8) rectangle (\n*5/8, \n*5/8) node[midway] {$\mA_{11, 12}$};
            \draw[fill=\frcolor] (\n*5/8, \n*3/8) rectangle (\n*6/8, \n*4/8) node[midway] {$\mA_{12, 13}$};
            \draw[fill=\frcolor] (\n*6/8, \n*2/8) rectangle (\n*7/8, \n*3/8) node[midway] {$\mA_{13, 14}$};
            \draw[fill=\frcolor] (\n*7/8, \n*1/8) rectangle (\n*8/8, \n*2/8) node[midway] {$\mA_{14, 15}$};

            \draw[fill=\frcolor] (\n*0/8, \n*6/8) rectangle (\n*1/8, \n*7/8) node[midway] {$\mA_{ 9,  8}$};
            \draw[fill=\frcolor] (\n*1/8, \n*5/8) rectangle (\n*2/8, \n*6/8) node[midway] {$\mA_{10,  9}$};
            \draw[fill=\frcolor] (\n*2/8, \n*4/8) rectangle (\n*3/8, \n*5/8) node[midway] {$\mA_{11, 10}$};
            \draw[fill=\frcolor] (\n*3/8, \n*3/8) rectangle (\n*4/8, \n*4/8) node[midway] {$\mA_{12, 11}$};
            \draw[fill=\frcolor] (\n*4/8, \n*2/8) rectangle (\n*5/8, \n*3/8) node[midway] {$\mA_{13, 12}$};
            \draw[fill=\frcolor] (\n*5/8, \n*1/8) rectangle (\n*6/8, \n*2/8) node[midway] {$\mA_{14, 13}$};
            \draw[fill=\frcolor] (\n*6/8, \n*0/8) rectangle (\n*7/8, \n*1/8) node[midway] {$\mA_{15, 14}$};

            \draw[fill=\frcolor] (\n*0/8, \n*7/8) rectangle (\n*1/8, \n*8/8) node[midway] {$\mA_{ 8,  8}$};
            \draw[fill=\frcolor] (\n*1/8, \n*6/8) rectangle (\n*2/8, \n*7/8) node[midway] {$\mA_{ 9,  9}$};
            \draw[fill=\frcolor] (\n*2/8, \n*5/8) rectangle (\n*3/8, \n*6/8) node[midway] {$\mA_{10, 10}$};
            \draw[fill=\frcolor] (\n*3/8, \n*4/8) rectangle (\n*4/8, \n*5/8) node[midway] {$\mA_{11, 11}$};
            \draw[fill=\frcolor] (\n*4/8, \n*3/8) rectangle (\n*5/8, \n*4/8) node[midway] {$\mA_{12, 12}$};
            \draw[fill=\frcolor] (\n*5/8, \n*2/8) rectangle (\n*6/8, \n*3/8) node[midway] {$\mA_{13, 13}$};
            \draw[fill=\frcolor] (\n*6/8, \n*1/8) rectangle (\n*7/8, \n*2/8) node[midway] {$\mA_{14, 14}$};
            \draw[fill=\frcolor] (\n*7/8, \n*0/8) rectangle (\n*8/8, \n*1/8) node[midway] {$\mA_{15, 15}$};
        \end{scope}
        \begin{scope}[shift={($(scope1.west)+(-\n*7/20, 0)$)}, grow=right]
            \node{$I_1$}
            child{
                node{$I_3$}
                child{
                    node{$I_7$}
                    child{node{$I_{15}$}}
                    child{node{$I_{14}$}}
                }
                child{
                    node{$I_6$}
                    child{node{$I_{13}$}}
                    child{node{$I_{12}$}}
                }
            }
            child{
                node{$I_2$}
                child{
                    node{$I_5$}
                    child{node{$I_{11}$}}
                    child{node{$I_{10}$}}
                }
                child{
                    node{$I_4$}
                    child{node{$I_{9}$}}
                    child{node{$I_{8}$}}
                }
            };
        \end{scope}
        \begin{scope}[shift={($(scope1.north)+(0, \n*7/20)$)}, grow=down]
            \node{$I_1$}
            child{
                node{$I_2$}
                child{
                    node{$I_4$}
                    child{node{$I_{8}$}}
                    child{node{$I_{9}$}}
                }
                child{
                    node{$I_5$}
                    child{node{$I_{10}$}}
                    child{node{$I_{11}$}}
                }
            }
            child{
                node{$I_3$}
                child{
                    node{$I_6$}
                    child{node{$I_{12}$}}
                    child{node{$I_{13}$}}
                }
                child{
                    node{$I_7$}
                    child{node{$I_{14}$}}
                    child{node{$I_{15}$}}
                }
            };
        \end{scope}
    \end{tikzpicture}
    \caption{
        An $\h^1$ matrix with depth 3
        based on a grid over $\lbrack 0, 1 \rbrack$.
        Admissible blocks are shown in blue,
        and inadmissible blocks are shown in red.
    }
    \label{fig:matrix_strong_1d_lvl3}
\end{figure}

An $\h^1$ matrix $\mA$ is said to have uniform $\h^1$ structure
if it further satisfies the condition that there exist low-rank basis matrices
$\ulong_\alpha$ spanning
$\mA \left( I_\alpha, \cup_{\beta \in \alphainteraction} I_\beta \right)$
and
$\vlong_\alpha$ spanning
$\mA \left( \cup_{\beta \in \alphainteraction} I_\beta, I_\alpha \right)$
for each box $\alpha$.
Constructing a compressed representation of a uniform $\h^1$ matrix
consists of finding the matrices $\ulong_\alpha, \vlong_\alpha$ for each box $\alpha$
and the matrices $\mB_{\alpha, \beta}$ such that
\begin{equation*}
    \mA(I_\alpha, I_\beta) \approx \ulong_\alpha \mB_{\alpha, \beta} \vlong_\beta
\end{equation*}
for each admissible pair $(\alpha, \beta)$
and storing the inadmissible blocks
corresponding to neighbor interactions of boxes in level $L$.

A uniform $\h^1$ matrix $\mA$ is said to have $\h^2$ structure
if there exist low-rank basis matrices $\ulong_\alpha$ spanning
$\mA \left( I_\alpha, \left(\cup_{\beta \in \alphanear} I_\beta \right)^c \right)$
and $\vlong_\alpha$ spanning
$\mA \left( \left( \cup_{\beta \in \alphanear} I_\beta \right)^c, I_\alpha \right)$
for each box $\alpha$.
Then the basis matrices of a non-leaf box
can be expressed in terms of the basis matrices
of its children.
For example, if $\tau$ is a box
with children $\alpha, \beta$,
then
\begin{equation}\label{eq:nested_basis}
    \ulong_\tau
    =
    \begin{bmatrix}
        \ulong_\alpha & \mzero \\
        \mzero & \ulong_\beta \\
    \end{bmatrix}
    \ushort_\tau,
\end{equation}
where $\ulong_\gamma$
is the ``long'' column basis matrix of size
$\vert \gamma \vert \times k$
associated with $\gamma \in \{\tau, \alpha, \beta\}$,
and
$\ushort_\tau$ is a ``short'' basis matrix
of size $2k \times k$.
Analogous relationships must hold
for row basis matrices
$\vlong_\tau, \vlong_\alpha, \vlong_\beta$.
Using such nested basis matrices
eliminates the need to explicitly store
$\ulong_\tau$
since it is fully specified by
$\ulong_\alpha, \ulong_\beta$ and $\ushort_\tau$.
Likewise, if $\alpha$ or $\beta$ have children of their own,
then their basis matrices
will be expressed in terms of their own small basis matrices
and the basis matrices of their children.

\section{Compressing rank-structured matrices with graph coloring} \label{sec:compress}

\subsection{\texorpdfstring{$\h^1$}{H1} matrix compression}
\label{sec:h1}

In this section, we present the process of constructing an $\h^1$ representation of a matrix
by applying \cref{alg:compress_double_random}
to compute low-rank approximations of the admissible blocks
of levels $2, \dots, L$
and then extracting the inadmissible blocks of level $L$.
We apply $\mtx{A}$ and $\mtx{A}^*$ to a set of carefully constructed test matrices,
and from those products we extract randomized samples of each admissible block.
We demonstrate the techniques on a matrix
shown in \cref{fig:matrix_strong_1d_lvl3},
which is based on points in one dimension,
but the algorithm generalizes in a straightforward way
to points in higher dimensions.

\subsubsection{Compressing level 2}
\label{sec:h1_level2}

We construct the compressed representation by processing levels
of the tree in sequence
from the coarsest level to the finest.
There are no admissible blocks associated with levels 0 and 1,
so we begin by computing low-rank approximations of the 6 admissible blocks
associated with level 2 of the tree.
To that end, we define four test matrices of size $N \times r$
\begin{equation*}
    \mOmega_1
    =
    \begin{bmatrix}
        \mG_4 \\
        \mzero \\
        \mzero \\
        \mzero \\
    \end{bmatrix},
    \qquad
    \mOmega_2
    =
    \begin{bmatrix}
        \mzero \\
        \mG_5 \\
        \mzero \\
        \mzero \\
    \end{bmatrix},
    \qquad
    \mOmega_3
    =
    \begin{bmatrix}
        \mzero \\
        \mzero \\
        \mG_6 \\
        \mzero \\
    \end{bmatrix},
    \qquad
    \mOmega_4
    =
    \begin{bmatrix}
        \mzero \\
        \mzero \\
        \mzero \\
        \mG_7 \\
    \end{bmatrix},
\end{equation*}
where
$\mG_4, \mG_5, \mG_6, \mG_7$
are random matrices of size $N/4 \times r$
whose entries are drawn from the standard normal distribution.
We invoke the black-box matrix-vector product routine to evaluate the products
$\mY_i = \mA \mOmega_i, i \in \{ 1, 2, 3, 4 \}$, with the following structures.
Contained within these products are randomized samples
of each admissible block of level 2.
\begin{equation}
\label{eq:h1testmatrices}
\begin{aligned}
    \mY_1
    =
    \mA \mOmega_1
    =
    \begin{bmatrix}
        * \\
        * \\
        \mA_{6, 4} \mG_4 \\
        \mA_{7, 4} \mG_4 \\
    \end{bmatrix},
    \qquad
    \mY_2
    =
    \mA \mOmega_2
    =
    \begin{bmatrix}
        * \\
        * \\
        * \\
        \mA_{7, 5} \mG_5 \\
    \end{bmatrix} \\
    \qquad
    \mY_3
    =
    \mA \mOmega_3
    =
    \begin{bmatrix}
        \mA_{4, 6} \mG_6 \\
        * \\
        * \\
        * \\
    \end{bmatrix},
    \qquad
    \mY_4
    =
    \mA \mOmega_4
    =
    \begin{bmatrix}
        \mA_{4, 7} \mG_7 \\
        \mA_{5, 7} \mG_7 \\
        * \\
        * \\
    \end{bmatrix}
\end{aligned}
\end{equation}

We then obtain a basis matrix for the column space of each admissible block
by orthonormalizing the relevant block of one of the sample matrices.

\begin{equation}
\label{eq:orthonormalize}
\begin{aligned}
    \ulong_{4, 6} &= \qr(\mY_3(I_4, :)) \\
    \ulong_{4, 7} &= \qr(\mY_4(I_4, :)) \\
    \ulong_{5, 7} &= \qr(\mY_4(I_5, :)) \\
    \ulong_{6, 4} &= \qr(\mY_1(I_6, :)) \\
    \ulong_{7, 4} &= \qr(\mY_1(I_7, :)) \\
    \ulong_{7, 5} &= \qr(\mY_2(I_7, :)) \\
\end{aligned}
\end{equation}

To find a basis matrix for the row space of each admissible block,
we follow a similar process using $\mA^*$ in place of $\mA$.
That is, we compute another set of sample matrices
$\mZ_i = \mA^* \mPsi_i, i \in \{ 1, 2, 3, 4 \}$,
from which we orthonormalize the relevant blocks
to obtain row bases $\vlong_{\alpha, \beta}$
for each admissible pair $(\alpha, \beta)$.

Finally, we solve for the matrices $\mB_{\alpha, \beta}$ as follows.
Note that the products $\mA_{\alpha, \beta} \mG_\beta$
have already been obtained
from the samples in \cref{eq:h1testmatrices}.
\begin{equation}\label{eq:solveforb}
    \mB_{\alpha, \beta}
    =
    \left( \mG_\alpha \ulong_\alpha \right)^\dagger
    \mG_\alpha \mA_{\alpha, \beta} \mG_\beta
    \left( \vlong_\beta \mG_\beta \right)^\dagger
    \quad
    \text{for each admissible pair $(\alpha, \beta)$}
\end{equation}

\subsubsection{Compressing levels \texorpdfstring{$3, ..., L$}{3, ..., L}}
\label{sec:h1_level3}

After we have obtained low-rank approximations of the blocks
associated with level 2 of the tree, we proceed to level 3.
One approach would be to extend the procedure for level 2
by using one test matrix corresponding to each box,
for a total of eight test matrices.
That approach would grow to be prohibitively expensive
for finer levels of the tree
as the number of matrix-vector products required
would grow proportionately with the number of boxes.
Instead, we present a more efficient procedure,
which requires a number of matrix-vector products
that is bounded across all levels.

Following~\cite{martinsson2016compressing},
we define the level-$l$ truncated matrix $\mA^{(l)}$
to be the matrix obtained by replacing with zeros every block of $\mA$
corresponding to levels finer than level $l$.
Note that a matrix-vector product involving $\mA^{(l)}$
can be computed inexpensively
using the low-rank approximations already computed
when processing levels $2, \dots, l-1$.
The structures of the level-2 truncated matrix
and the difference $\mA - \mA^{(2)}$ are shown below.

\begin{equation*}
    \mA^{(2)}
    =
    \begin{tikzpicture}[scale=4, line width=1pt, baseline=2.0cm]
        \def\n{1cm}
        \draw (0, 0) rectangle (\n, \n) node[midway] {\mzero};

        \draw[fill=\lrcolor] (\n*2/4, \n*3/4) rectangle (\n*3/4, \n*4/4) node[midway] {$\mA_{4, 6}$};
        \draw[fill=\lrcolor] (\n*3/4, \n*3/4) rectangle (\n*4/4, \n*4/4) node[midway] {$\mA_{4, 7}$};
        \draw[fill=\lrcolor] (\n*3/4, \n*2/4) rectangle (\n*4/4, \n*3/4) node[midway] {$\mA_{5, 7}$};
        \draw[fill=\lrcolor] (\n*0/4, \n*1/4) rectangle (\n*1/4, \n*2/4) node[midway] {$\mA_{6, 4}$};
        \draw[fill=\lrcolor] (\n*0/4, \n*0/4) rectangle (\n*1/4, \n*1/4) node[midway] {$\mA_{7, 4}$};
        \draw[fill=\lrcolor] (\n*1/4, \n*0/4) rectangle (\n*2/4, \n*1/4) node[midway] {$\mA_{7, 5}$};
    \end{tikzpicture}
    \quad,\qquad
    \mA - \mA^{(2)}
    =
    \begin{tikzpicture}[scale=4, line width=1pt, baseline=2.0cm]
        \def\n{1cm}

        \draw (0, 0) rectangle (\n, \n);

        \draw (\n*2/4, \n*3/4) rectangle (\n*3/4, \n*4/4) node[midway] {\mzero};
        \draw (\n*3/4, \n*3/4) rectangle (\n*4/4, \n*4/4) node[midway] {\mzero};
        \draw (\n*3/4, \n*2/4) rectangle (\n*4/4, \n*3/4) node[midway] {\mzero};
        \draw (\n*0/4, \n*1/4) rectangle (\n*1/4, \n*2/4) node[midway] {\mzero};
        \draw (\n*0/4, \n*0/4) rectangle (\n*1/4, \n*1/4) node[midway] {\mzero};
        \draw (\n*1/4, \n*0/4) rectangle (\n*2/4, \n*1/4) node[midway] {\mzero};

        \draw[fill=\lrcolor] (\n*2/8, \n*7/8) rectangle (\n*3/8, \n*8/8); %
        \draw[fill=\lrcolor] (\n*3/8, \n*7/8) rectangle (\n*4/8, \n*8/8); %
        \draw[fill=\lrcolor] (\n*3/8, \n*6/8) rectangle (\n*4/8, \n*7/8); %
        \draw[fill=\lrcolor] (\n*0/8, \n*5/8) rectangle (\n*1/8, \n*6/8); %
        \draw[fill=\lrcolor] (\n*4/8, \n*5/8) rectangle (\n*5/8, \n*6/8); %
        \draw[fill=\lrcolor] (\n*5/8, \n*5/8) rectangle (\n*6/8, \n*6/8); %
        \draw[fill=\lrcolor] (\n*0/8, \n*4/8) rectangle (\n*1/8, \n*5/8); %
        \draw[fill=\lrcolor] (\n*1/8, \n*4/8) rectangle (\n*2/8, \n*5/8); %
        \draw[fill=\lrcolor] (\n*5/8, \n*4/8) rectangle (\n*6/8, \n*5/8); %
        \draw[fill=\lrcolor] (\n*2/8, \n*3/8) rectangle (\n*3/8, \n*4/8); %
        \draw[fill=\lrcolor] (\n*6/8, \n*3/8) rectangle (\n*7/8, \n*4/8); %
        \draw[fill=\lrcolor] (\n*7/8, \n*3/8) rectangle (\n*8/8, \n*4/8); %
        \draw[fill=\lrcolor] (\n*2/8, \n*2/8) rectangle (\n*3/8, \n*3/8); %
        \draw[fill=\lrcolor] (\n*3/8, \n*2/8) rectangle (\n*4/8, \n*3/8); %
        \draw[fill=\lrcolor] (\n*7/8, \n*2/8) rectangle (\n*8/8, \n*3/8); %
        \draw[fill=\lrcolor] (\n*4/8, \n*1/8) rectangle (\n*5/8, \n*2/8); %
        \draw[fill=\lrcolor] (\n*4/8, \n*0/8) rectangle (\n*5/8, \n*1/8); %
        \draw[fill=\lrcolor] (\n*5/8, \n*0/8) rectangle (\n*6/8, \n*1/8); %

        \draw[fill=\frcolor] (\n*1/8, \n*7/8) rectangle (\n*2/8, \n*8/8); %
        \draw[fill=\frcolor] (\n*2/8, \n*6/8) rectangle (\n*3/8, \n*7/8); %
        \draw[fill=\frcolor] (\n*3/8, \n*5/8) rectangle (\n*4/8, \n*6/8); %
        \draw[fill=\frcolor] (\n*4/8, \n*4/8) rectangle (\n*5/8, \n*5/8); %
        \draw[fill=\frcolor] (\n*5/8, \n*3/8) rectangle (\n*6/8, \n*4/8); %
        \draw[fill=\frcolor] (\n*6/8, \n*2/8) rectangle (\n*7/8, \n*3/8); %
        \draw[fill=\frcolor] (\n*7/8, \n*1/8) rectangle (\n*8/8, \n*2/8); %

        \draw[fill=\frcolor] (\n*0/8, \n*6/8) rectangle (\n*1/8, \n*7/8); %
        \draw[fill=\frcolor] (\n*1/8, \n*5/8) rectangle (\n*2/8, \n*6/8); %
        \draw[fill=\frcolor] (\n*2/8, \n*4/8) rectangle (\n*3/8, \n*5/8); %
        \draw[fill=\frcolor] (\n*3/8, \n*3/8) rectangle (\n*4/8, \n*4/8); %
        \draw[fill=\frcolor] (\n*4/8, \n*2/8) rectangle (\n*5/8, \n*3/8); %
        \draw[fill=\frcolor] (\n*5/8, \n*1/8) rectangle (\n*6/8, \n*2/8); %
        \draw[fill=\frcolor] (\n*6/8, \n*0/8) rectangle (\n*7/8, \n*1/8); %

        \draw[fill=\frcolor] (\n*0/8, \n*7/8) rectangle (\n*1/8, \n*8/8); %
        \draw[fill=\frcolor] (\n*1/8, \n*6/8) rectangle (\n*2/8, \n*7/8); %
        \draw[fill=\frcolor] (\n*2/8, \n*5/8) rectangle (\n*3/8, \n*6/8); %
        \draw[fill=\frcolor] (\n*3/8, \n*4/8) rectangle (\n*4/8, \n*5/8); %
        \draw[fill=\frcolor] (\n*4/8, \n*3/8) rectangle (\n*5/8, \n*4/8); %
        \draw[fill=\frcolor] (\n*5/8, \n*2/8) rectangle (\n*6/8, \n*3/8); %
        \draw[fill=\frcolor] (\n*6/8, \n*1/8) rectangle (\n*7/8, \n*2/8); %
        \draw[fill=\frcolor] (\n*7/8, \n*0/8) rectangle (\n*8/8, \n*1/8); %
    \end{tikzpicture}
\end{equation*}

We will sample the admissible blocks of level 3
by applying $\mA - \mA^{(2)}$ to a set of test matrices
subject to certain conditions on their sparsity structure.
For example, to isolate a sample of $\mA_{8, 10}$
(see \cref{fig:matrix_strong_1d_lvl3}),
we must avoid unwanted contributions from
$\mA_{8, 8}, \mA_{8, 9}, \mA_{8, 11}$,
so we multiply $\mA - \mA^{(2)}$ with a test matrix
whose rows indexed by $I_8, I_9, I_{11}$ are all zeros
and whose rows indexed by $I_{10}$
are filled with random values drawn from the standard normal distribution.
The contents of the other rows are irrelevant
for the purpose of sampling $\mA_{8,10}$
since they will be multiplied with zeros in $\mA - \mA^{(2)}$.
More generally, to sample some admissible block $\mA_{\alpha, \beta}$
of level 3,
we require a test matrix $\mOmega$ that satisfies the following
\emph{sampling constraints}.
\begin{equation}\label{eq:h1constraints}
\begin{aligned}
    \mOmega(I_\beta, :)  &= \mG_\beta & & \\
    \mOmega(I_\gamma, :) &= \mzero &\text{for all }& \gamma \in \alphanear \cup \alphainteraction \setminus \{\beta\}
\end{aligned}
\end{equation}
where $\mG_\beta$ is a random matrix of size $\vert I_\beta \vert$-by-$r$.
If test matrix $\mOmega$ satisfies those sampling constraints,
then the rows of the product $\mA \mOmega - \mA^{(2)} \mOmega$ indexed by $I_\alpha$
will contain a randomized sample of the column space of $\mA_{\alpha, \beta}$.

To sample all of the admissible blocks,
we require a set of test matrices $\{\mOmega_i\}$
such that for every admissible pair $(\alpha, \beta)$,
there is a test matrix within the set
that satisfies the constraints~\eqref{eq:h1constraints}
associated with that pair.
Moreover, we would like that set to be as small as possible
to minimize cost.

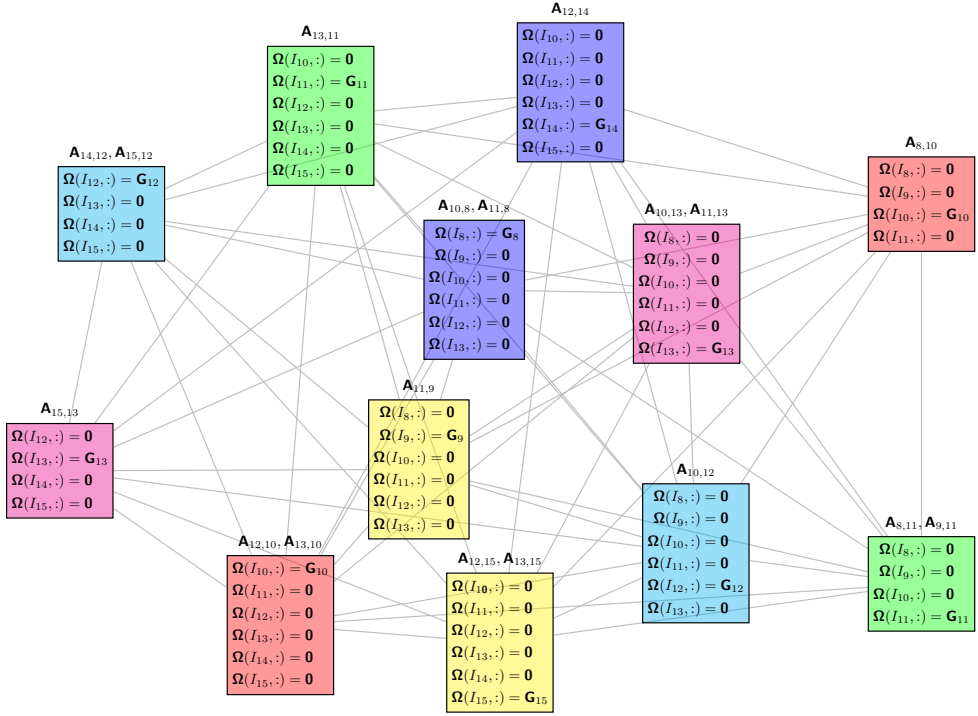
\begin{figure}[tb]
    \centering
    \resizebox{\columnwidth}{!}{
    \begin{tikzpicture}[spring electrical layout, node distance=4.0cm]
    \node (8x10)  [vertexstyle, label={$\mA_{8, 10}$}, fill=red] {
            $\begin{aligned}
                \mOmega(I_{8}, :) &= \mzero \\
                \mOmega(I_{9}, :) &= \mzero \\
                \mOmega(I_{10}, :) &= \mG_{10} \\
                \mOmega(I_{11}, :) &= \mzero \\
            \end{aligned}$
        };
    \node (8x11)  [vertexstyle, label={$\mA_{8, 11}, \mA_{9, 11}$}, fill=green] {
            $\begin{aligned}
                \mOmega(I_{8}, :) &= \mzero \\
                \mOmega(I_{9}, :) &= \mzero \\
                \mOmega(I_{10}, :) &= \mzero \\
                \mOmega(I_{11}, :) &= \mG_{11} \\
            \end{aligned}$
    };
    \node (10x8)  [vertexstyle, label={$\mA_{10, 8}, \mA_{11, 8}$}, fill=blue] {
            $\begin{aligned}
                \mOmega(I_{8}, :) &= \mG_{8} \\
                \mOmega(I_{9}, :) &= \mzero \\
                \mOmega(I_{10}, :) &= \mzero \\
                \mOmega(I_{11}, :) &= \mzero \\
                \mOmega(I_{12}, :) &= \mzero \\
                \mOmega(I_{13}, :) &= \mzero \\
            \end{aligned}$
    };
    \node (10x12) [vertexstyle, label={$\mA_{10, 12}$}, fill=cyan] {
            $\begin{aligned}
                \mOmega(I_{8}, :) &= \mzero \\
                \mOmega(I_{9}, :) &= \mzero \\
                \mOmega(I_{10}, :) &= \mzero \\
                \mOmega(I_{11}, :) &= \mzero \\
                \mOmega(I_{12}, :) &= \mG_{12} \\
                \mOmega(I_{13}, :) &= \mzero \\
            \end{aligned}$
    };
    \node (10x13) [vertexstyle, label={$\mA_{10, 13}, \mA_{11, 13}$}, fill=magenta] {
            $\begin{aligned}
                \mOmega(I_{8}, :) &= \mzero \\
                \mOmega(I_{9}, :) &= \mzero \\
                \mOmega(I_{10}, :) &= \mzero \\
                \mOmega(I_{11}, :) &= \mzero \\
                \mOmega(I_{12}, :) &= \mzero \\
                \mOmega(I_{13}, :) &= \mG_{13} \\
            \end{aligned}$
    };
    \node (11x9)  [vertexstyle, label={$\mA_{11, 9}$}, fill=yellow] {
            $\begin{aligned}
                \mOmega(I_{8}, :) &= \mzero \\
                \mOmega(I_{9}, :) &= \mG_{9} \\
                \mOmega(I_{10}, :) &= \mzero \\
                \mOmega(I_{11}, :) &= \mzero \\
                \mOmega(I_{12}, :) &= \mzero \\
                \mOmega(I_{13}, :) &= \mzero \\
            \end{aligned}$
    };
    \node (12x10) [vertexstyle, label={$\mA_{12, 10}, \mA_{13, 10}$}, fill=red] {
            $\begin{aligned}
                \mOmega(I_{10}, :) &= \mG_{10} \\
                \mOmega(I_{11}, :) &= \mzero \\
                \mOmega(I_{12}, :) &= \mzero \\
                \mOmega(I_{13}, :) &= \mzero \\
                \mOmega(I_{14}, :) &= \mzero \\
                \mOmega(I_{15}, :) &= \mzero \\
            \end{aligned}$
    };
    \node (12x14) [vertexstyle, label={$\mA_{12, 14}$}, fill=blue] {
            $\begin{aligned}
                \mOmega(I_{10}, :) &= \mzero \\
                \mOmega(I_{11}, :) &= \mzero \\
                \mOmega(I_{12}, :) &= \mzero \\
                \mOmega(I_{13}, :) &= \mzero \\
                \mOmega(I_{14}, :) &= \mG_{14} \\
                \mOmega(I_{15}, :) &= \mzero \\
            \end{aligned}$
    };
    \node (12x15) [vertexstyle, label={$\mA_{12, 15}, \mA_{13, 15}$}, fill=yellow] {
            $\begin{aligned}
                \mOmega(I_{1\mzero}, :) &= \mzero \\
                \mOmega(I_{11}, :) &= \mzero \\
                \mOmega(I_{12}, :) &= \mzero \\
                \mOmega(I_{13}, :) &= \mzero \\
                \mOmega(I_{14}, :) &= \mzero \\
                \mOmega(I_{15}, :) &= \mG_{15} \\
            \end{aligned}$
    };
    \node (13x11) [vertexstyle, label={$\mA_{13, 11}$}, fill=green] {
            $\begin{aligned}
                \mOmega(I_{10}, :) &= \mzero \\
                \mOmega(I_{11}, :) &= \mG_{11} \\
                \mOmega(I_{12}, :) &= \mzero \\
                \mOmega(I_{13}, :) &= \mzero \\
                \mOmega(I_{14}, :) &= \mzero \\
                \mOmega(I_{15}, :) &= \mzero \\
            \end{aligned}$
    };
    \node (14x12) [vertexstyle, label={$\mA_{14, 12}, \mA_{15, 12}$}, fill=cyan] {
            $\begin{aligned}
                \mOmega(I_{12}, :) &= \mG_{12} \\
                \mOmega(I_{13}, :) &= \mzero \\
                \mOmega(I_{14}, :) &= \mzero \\
                \mOmega(I_{15}, :) &= \mzero \\
            \end{aligned}$
    };
    \node (15x13) [vertexstyle, label={$\mA_{15, 13}$}, fill=magenta] {
            $\begin{aligned}
                \mOmega(I_{12}, :) &= \mzero \\
                \mOmega(I_{13}, :) &= \mG_{13} \\
                \mOmega(I_{14}, :) &= \mzero \\
                \mOmega(I_{15}, :) &= \mzero \\
            \end{aligned}$
    };

    \path (8x10) edge[edgestyle] (8x11);
    \path (8x10) edge[edgestyle] (10x8);
    \path (8x10) edge[edgestyle] (10x12);
    \path (8x10) edge[edgestyle] (10x13);
    \path (8x10) edge[edgestyle] (11x9);
    \path (8x10) edge[edgestyle] (12x14);
    \path (8x10) edge[edgestyle] (12x15);
    \path (8x10) edge[edgestyle] (13x11);

    \path (8x11) edge[edgestyle] (10x8);
    \path (8x11) edge[edgestyle] (10x12);
    \path (8x11) edge[edgestyle] (10x13);
    \path (8x11) edge[edgestyle] (11x9);
    \path (8x11) edge[edgestyle] (12x10);
    \path (8x11) edge[edgestyle] (12x14);
    \path (8x11) edge[edgestyle] (12x15);

    \path (10x8) edge[edgestyle] (10x12);
    \path (10x8) edge[edgestyle] (10x13);
    \path (10x8) edge[edgestyle] (11x9);
    \path (10x8) edge[edgestyle] (12x10);
    \path (10x8) edge[edgestyle] (13x11);
    \path (10x8) edge[edgestyle] (14x12);
    \path (10x8) edge[edgestyle] (15x13);

    \path (10x12) edge[edgestyle] (10x13);
    \path (10x12) edge[edgestyle] (11x9);
    \path (10x12) edge[edgestyle] (12x10);
    \path (10x12) edge[edgestyle] (12x14);
    \path (10x12) edge[edgestyle] (12x15);
    \path (10x12) edge[edgestyle] (13x11);
    \path (10x12) edge[edgestyle] (15x13);

    \path (10x13) edge[edgestyle] (11x9);
    \path (10x13) edge[edgestyle] (12x10);
    \path (10x13) edge[edgestyle] (12x14);
    \path (10x13) edge[edgestyle] (12x15);
    \path (10x13) edge[edgestyle] (13x11);
    \path (10x13) edge[edgestyle] (14x12);

    \path (11x9) edge[edgestyle] (12x10);
    \path (11x9) edge[edgestyle] (13x11);
    \path (11x9) edge[edgestyle] (14x12);
    \path (11x9) edge[edgestyle] (15x13);

    \path (12x10) edge[edgestyle] (12x14);
    \path (12x10) edge[edgestyle] (12x15);
    \path (12x10) edge[edgestyle] (13x11);
    \path (12x10) edge[edgestyle] (14x12);
    \path (12x10) edge[edgestyle] (15x13);

    \path (12x14) edge[edgestyle] (12x15);
    \path (12x14) edge[edgestyle] (13x11);
    \path (12x14) edge[edgestyle] (14x12);
    \path (12x14) edge[edgestyle] (15x13);

    \path (12x15) edge[edgestyle] (13x11);
    \path (12x15) edge[edgestyle] (14x12);
    \path (12x15) edge[edgestyle] (15x13);

    \path (13x11) edge[edgestyle] (14x12);
    \path (13x11) edge[edgestyle] (15x13);

    \path (14x12) edge[edgestyle] (15x13);

    \end{tikzpicture}
    }
    \caption{The constraint incompatibility graph corresponding to the 18 admissible blocks belonging to level 3
    of the matrix shown in \fref{fig:matrix_strong_1d_lvl3}.
    Each vertex corresponds to a distinct set of sampling constraints \eqref{eq:h1constraints}.
    Edges connect pairs of vertices that are incompatible.
    The number of vertices is less than the number of admissible blocks
    since some admissible blocks share the same set of sampling constraints.}
    \label{fig:graph_strong_1d_lvl3}
\end{figure}

In order to minimize the number of test matrices,
we aim to form a small number of groups of compatible sampling constraints.
We say that two sets of sampling constraints are \emph{compatible}
if it is possible to form a test matrix $\mOmega$
that satisfies both sets of constraints.
We then define a
\emph{constraint incompatibility graph},
which represents compatibility relationships between pairs of constraint sets.
The graph corresponding to the 18 admissible blocks belonging to level 3
is depicted in \fref{fig:graph_strong_1d_lvl3}.

\begin{definition}[Constraint incompatibility graph]
\label{def:incompatibility_graph}
The constraint incompatibility graph for level $l$ of the tree
is the graph
in which each vertex corresponds
to a distinct constraint set~\eqref{eq:h1constraints},
and pairs of vertices are connected by an edge
if their corresponding constraint sets are incompatible.
\end{definition}

We then compute a vertex coloring of the constraint incompatibility graph
using the DSatur algorithm (\cref{alg:dsatur}).
For a valid coloring of the graph,
each subset of vertices sharing the same color
represents a mutually compatible collection of sampling constraints.
Then for each color, we can define one test matrix
that satisfies all of the sampling constraints
associated with the vertices of that color.
The coloring shown in \fref{fig:graph_strong_1d_lvl3}
yields test matrices with the following structures.
\begin{equation}
\label{eq:neartestmatrices}
    \begin{bmatrix}
        \mOmega_1 & \mOmega_2 & \mOmega_3 & \mOmega_4 & \mOmega_5 & \mOmega_6
    \end{bmatrix}
    =
    \begin{bmatrix}
        \mG_{8}  & \mzero   & \mzero   & \mzero   & \mzero   & \mzero   \\
        \mzero   & \mG_{9}  & \mzero   & \mzero   & \mzero   & \mzero   \\
        \mzero   & \mzero   & \mG_{10} & \mzero   & \mzero   & \mzero   \\
        \mzero   & \mzero   & \mzero   & \mG_{11} & \mzero   & \mzero   \\
        \mzero   & \mzero   & \mzero   & \mzero   & \mG_{12} & \mzero   \\
        \mzero   & \mzero   & \mzero   & \mzero   & \mzero   & \mG_{13} \\
        \mG_{14} & \mzero   & \mzero   & \mzero   & \mzero   & \mzero   \\
        \mzero   & \mG_{15} & \mzero   & \mzero   & \mzero   & \mzero   \\
    \end{bmatrix}
\end{equation}

As in \cref{sec:h1_level2},
we evaluate the samples $\mY_i = \mA \mOmega_i - \mA^{(2)} \mOmega_i$
and orthonormalize the relevant blocks
to obtain orthonormal bases $\ulong_{\alpha, \beta}$
of the column spaces of the admissible blocks.
A similar process yields
orthormal bases $\vlong_{\alpha, \beta}$
of the row spaces.
Finally, we solve for $\mB_{\alpha, \beta}$
again using \eqref{eq:solveforb}.

\subsubsection{Extracting inadmissible blocks of the leaf level}
\label{sec:leaf}

\begin{figure}[tb]
    \centering
    \resizebox{\columnwidth}{!}{
    \begin{tikzpicture}[scale=2.5]
        \node (8x8)  [vertexstyle, label={$\mA_{8, 8}$}, fill=red] at (0, 7) {
            $\begin{aligned}
                \mOmega(I_{8}, :) &= \mI \\
                \mOmega(I_{9}, :) &= \mzero \\
            \end{aligned}$
        };
        \node (8x9)  [vertexstyle, label={$\mA_{8, 9}$}, fill=green] at (1, 7) {
            $\begin{aligned}
                \mOmega(I_{8}, :) &= \mzero \\
                \mOmega(I_{9}, :) &= \mI \\
            \end{aligned}$
        };

        \node (9x8)  [vertexstyle, label={$\mA_{9, 8}$}, fill=red] at (0, 6) {
            $\begin{aligned}
                \mOmega(I_{8}, :) &= \mI \\
                \mOmega(I_{9}, :) &= \mzero \\
                \mOmega(I_{10}, :) &= \mzero \\
            \end{aligned}$
        };
        \node (9x9)  [vertexstyle, label={$\mA_{9, 9}$}, fill=green] at (1, 6) {
            $\begin{aligned}
                \mOmega(I_{8}, :) &= \mzero \\
                \mOmega(I_{9}, :) &= \mI \\
                \mOmega(I_{10}, :) &= \mzero \\
            \end{aligned}$
        };
        \node (9x10)  [vertexstyle, label={$\mA_{9, 10}$}, fill=blue] at (2, 6) {
            $\begin{aligned}
                \mOmega(I_{8}, :) &= \mzero \\
                \mOmega(I_{9}, :) &= \mzero \\
                \mOmega(I_{10}, :) &= \mI \\
            \end{aligned}$
        };

        \node (10x9)  [vertexstyle, label={$\mA_{10, 9}$}, fill=green] at (1, 5) {
            $\begin{aligned}
                \mOmega(I_{9}, :) &= \mI \\
                \mOmega(I_{10}, :) &= \mzero \\
                \mOmega(I_{11}, :) &= \mzero \\
            \end{aligned}$
        };
        \node (10x10)  [vertexstyle, label={$\mA_{10, 10}$}, fill=blue]  at (2, 5) {
            $\begin{aligned}
                \mOmega(I_{9}, :) &= \mzero \\
                \mOmega(I_{10}, :) &= \mI \\
                \mOmega(I_{11}, :) &= \mzero \\
            \end{aligned}$
        };
        \node (10x11)  [vertexstyle, label={$\mA_{10, 11}$}, fill=red] at (3, 5) {
            $\begin{aligned}
                \mOmega(I_{9}, :) &= \mzero \\
                \mOmega(I_{10}, :) &= \mzero \\
                \mOmega(I_{11}, :) &= \mI \\
            \end{aligned}$
        };

        \node (11x10)  [vertexstyle, label={$\mA_{11, 10}$}, fill=blue] at (2, 4) {
            $\begin{aligned}
                \mOmega(I_{10}, :) &= \mI \\
                \mOmega(I_{11}, :) &= \mzero \\
                \mOmega(I_{12}, :) &= \mzero \\
            \end{aligned}$
        };
        \node (11x11)  [vertexstyle, label={$\mA_{11, 11}$}, fill=red] at (3, 4) {
            $\begin{aligned}
                \mOmega(I_{10}, :) &= \mzero \\
                \mOmega(I_{11}, :) &= \mI \\
                \mOmega(I_{12}, :) &= \mzero \\
            \end{aligned}$
        };
        \node (11x12)  [vertexstyle, label={$\mA_{11, 12}$}, fill=green] at (4, 4) {
            $\begin{aligned}
                \mOmega(I_{10}, :) &= \mzero \\
                \mOmega(I_{11}, :) &= \mzero \\
                \mOmega(I_{12}, :) &= \mI \\
            \end{aligned}$
        };

        \node (12x11)  [vertexstyle, label={$\mA_{12, 11}$}, fill=red] at (3, 3) {
            $\begin{aligned}
                \mOmega(I_{11}, :) &= \mI \\
                \mOmega(I_{12}, :) &= \mzero \\
                \mOmega(I_{13}, :) &= \mzero \\
            \end{aligned}$
        };
        \node (12x12)  [vertexstyle, label={$\mA_{12, 12}$}, fill=green] at (4, 3) {
            $\begin{aligned}
                \mOmega(I_{11}, :) &= \mzero \\
                \mOmega(I_{12}, :) &= \mI \\
                \mOmega(I_{13}, :) &= \mzero \\
            \end{aligned}$
        };
        \node (12x13)  [vertexstyle, label={$\mA_{12, 13}$}, fill=blue] at (5, 3) {
            $\begin{aligned}
                \mOmega(I_{11}, :) &= \mzero \\
                \mOmega(I_{12}, :) &= \mzero \\
                \mOmega(I_{13}, :) &= \mI \\
            \end{aligned}$
        };

        \node (13x12)  [vertexstyle, label={$\mA_{13, 12}$}, fill=green] at (4, 2) {
            $\begin{aligned}
                \mOmega(I_{12}, :) &= \mI \\
                \mOmega(I_{13}, :) &= \mzero \\
                \mOmega(I_{14}, :) &= \mzero \\
            \end{aligned}$
        };
        \node (13x13)  [vertexstyle, label={$\mA_{13, 13}$}, fill=blue] at (5, 2) {
            $\begin{aligned}
                \mOmega(I_{12}, :) &= \mzero \\
                \mOmega(I_{13}, :) &= \mI \\
                \mOmega(I_{14}, :) &= \mzero \\
            \end{aligned}$
        };
        \node (13x14)  [vertexstyle, label={$\mA_{13, 14}$}, fill=red] at (6, 2) {
            $\begin{aligned}
                \mOmega(I_{12}, :) &= \mzero \\
                \mOmega(I_{13}, :) &= \mzero \\
                \mOmega(I_{14}, :) &= \mI \\
            \end{aligned}$
        };

        \node (14x13)  [vertexstyle, label={$\mA_{14, 13}$}, fill=blue] at (5, 1) {
            $\begin{aligned}
                \mOmega(I_{13}, :) &= \mI \\
                \mOmega(I_{14}, :) &= \mzero \\
                \mOmega(I_{15}, :) &= \mzero \\
            \end{aligned}$
        };
        \node (14x14)  [vertexstyle, label={$\mA_{14, 14}$}, fill=red] at (6, 1) {
            $\begin{aligned}
                \mOmega(I_{13}, :) &= \mzero \\
                \mOmega(I_{14}, :) &= \mI \\
                \mOmega(I_{15}, :) &= \mzero \\
            \end{aligned}$
        };
        \node (14x15)  [vertexstyle, label={$\mA_{14, 15}$}, fill=green] at (7, 1) {
            $\begin{aligned}
                \mOmega(I_{13}, :) &= \mzero \\
                \mOmega(I_{14}, :) &= \mzero \\
                \mOmega(I_{15}, :) &= \mI \\
            \end{aligned}$
        };

        \node (15x14)  [vertexstyle, label={$\mA_{15, 14}$}, fill=red] at (6, 0) {
            $\begin{aligned}
                \mOmega(I_{14}, :) &= \mI \\
                \mOmega(I_{15}, :) &= \mzero \\
            \end{aligned}$
        };
        \node (15x15)  [vertexstyle, label={$\mA_{15, 15}$}, fill=green] at (7, 0) {
            $\begin{aligned}
                \mOmega(I_{14}, :) &= \mzero \\
                \mOmega(I_{15}, :) &= \mI \\
            \end{aligned}$
        };

    \path (8x8) edge[edgestyle] (8x9);
    \path (8x8) edge[edgestyle] (9x9);
    \path (8x8) edge[edgestyle] (10x9);
    \path (8x8) edge[edgestyle] (9x10);

    \path (9x8) edge[edgestyle] (8x9);
    \path (9x8) edge[edgestyle] (9x9);
    \path (9x8) edge[edgestyle] (10x9);
    \path (9x8) edge[edgestyle] (9x10);
    \path (9x8) edge[edgestyle] (10x10);
    \path (9x8) edge[edgestyle] (11x10);

    \path (8x9) edge[edgestyle] (9x10);
    \path (8x9) edge[edgestyle] (10x10);
    \path (8x9) edge[edgestyle] (10x11);

    \path (9x9) edge[edgestyle] (9x10);
    \path (9x9) edge[edgestyle] (10x10);
    \path (9x9) edge[edgestyle] (11x10);
    \path (9x9) edge[edgestyle] (10x11);

    \path (10x9) edge[edgestyle] (9x10);
    \path (10x9) edge[edgestyle] (10x10);
    \path (10x9) edge[edgestyle] (11x10);
    \path (10x9) edge[edgestyle] (10x11);
    \path (10x9) edge[edgestyle] (11x11);
    \path (10x9) edge[edgestyle] (12x11);

    \path (9x10) edge[edgestyle] (10x11);
    \path (9x10) edge[edgestyle] (11x11);
    \path (9x10) edge[edgestyle] (11x12);

    \path (10x10) edge[edgestyle] (10x11);
    \path (10x10) edge[edgestyle] (11x11);
    \path (10x10) edge[edgestyle] (12x11);
    \path (10x10) edge[edgestyle] (11x12);

    \path (11x10) edge[edgestyle] (10x11);
    \path (11x10) edge[edgestyle] (11x11);
    \path (11x10) edge[edgestyle] (12x11);
    \path (11x10) edge[edgestyle] (11x12);
    \path (11x10) edge[edgestyle] (12x12);
    \path (11x10) edge[edgestyle] (13x12);

    \path (10x11) edge[edgestyle] (11x12);
    \path (10x11) edge[edgestyle] (12x12);
    \path (10x11) edge[edgestyle] (12x13);

    \path (11x11) edge[edgestyle] (11x12);
    \path (11x11) edge[edgestyle] (12x12);
    \path (11x11) edge[edgestyle] (13x12);
    \path (11x11) edge[edgestyle] (12x13);

    \path (12x11) edge[edgestyle] (11x12);
    \path (12x11) edge[edgestyle] (12x12);
    \path (12x11) edge[edgestyle] (13x12);
    \path (12x11) edge[edgestyle] (12x13);
    \path (12x11) edge[edgestyle] (13x13);
    \path (12x11) edge[edgestyle] (14x13);

    \path (11x12) edge[edgestyle] (12x13);
    \path (11x12) edge[edgestyle] (13x13);
    \path (11x12) edge[edgestyle] (13x14);

    \path (12x12) edge[edgestyle] (12x13);
    \path (12x12) edge[edgestyle] (13x13);
    \path (12x12) edge[edgestyle] (14x13);
    \path (12x12) edge[edgestyle] (13x14);

    \path (13x12) edge[edgestyle] (12x13);
    \path (13x12) edge[edgestyle] (13x13);
    \path (13x12) edge[edgestyle] (14x13);
    \path (13x12) edge[edgestyle] (13x14);
    \path (13x12) edge[edgestyle] (14x14);
    \path (13x12) edge[edgestyle] (15x14);

    \path (12x13) edge[edgestyle] (13x14);
    \path (12x13) edge[edgestyle] (14x14);
    \path (12x13) edge[edgestyle] (14x15);

    \path (13x13) edge[edgestyle] (13x14);
    \path (13x13) edge[edgestyle] (14x14);
    \path (13x13) edge[edgestyle] (15x14);
    \path (13x13) edge[edgestyle] (14x15);

    \path (14x13) edge[edgestyle] (13x14);
    \path (14x13) edge[edgestyle] (14x14);
    \path (14x13) edge[edgestyle] (15x14);
    \path (14x13) edge[edgestyle] (14x15);
    \path (14x13) edge[edgestyle] (15x15);

    \path (13x14) edge[edgestyle] (14x15);
    \path (13x14) edge[edgestyle] (15x15);

    \path (14x14) edge[edgestyle] (14x15);
    \path (14x14) edge[edgestyle] (15x15);

    \path (15x14) edge[edgestyle] (14x15);
    \path (15x14) edge[edgestyle] (15x15);

    \end{tikzpicture}
    }
    \caption{Incompatibility graph for inadmissible blocks belonging to level $L$.}
    \label{fig:graph_strong_1d_leaf}
\end{figure}
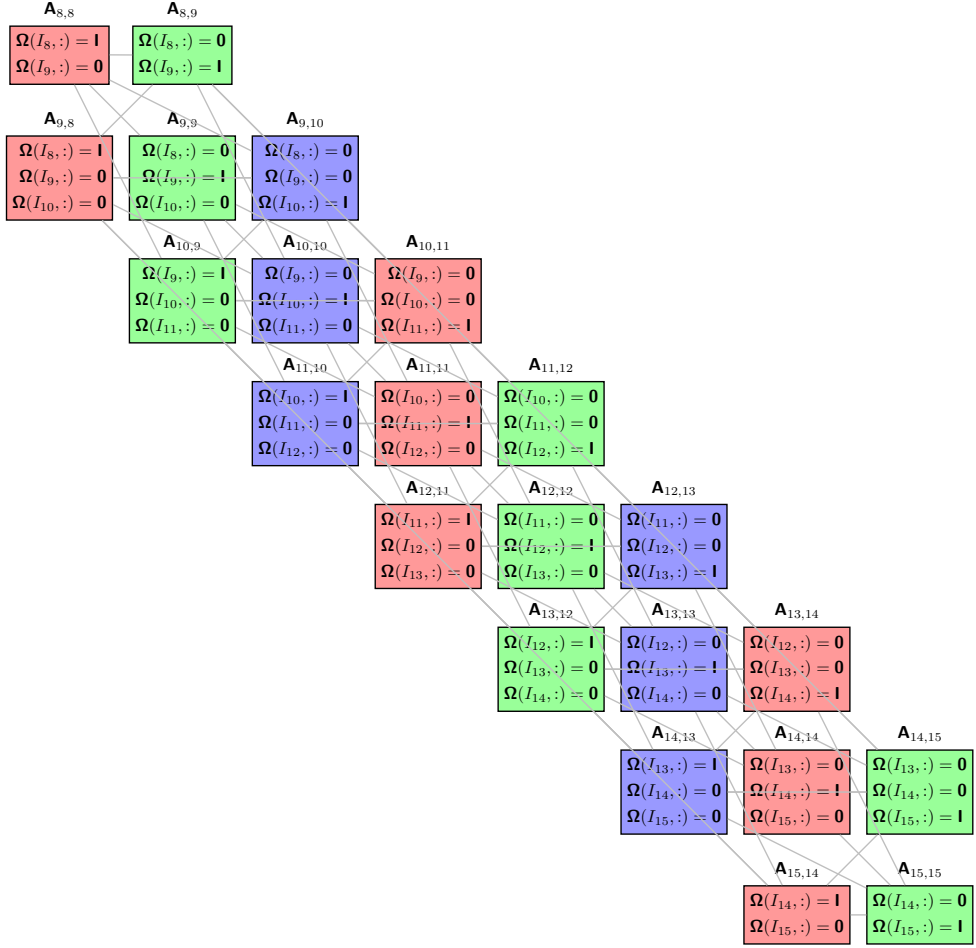

Once we have computed low-rank approximations of the admissible blocks for every level,
we finally extract the inadmissible blocks of the leaf level.
Since the inadmissible blocks are not necessarily low-rank,
they cannot be recovered from a small number of randomized samples.
Instead, we use test matrices that will multiply the inadmissible blocks
with appropriately sized identity matrices.
Also, at this point we already have low-rank approximations
of the admissible blocks belonging to level $L$,
so the test matrices
only need to avoid contributions from other inadmissible blocks,
resulting in fewer constraints than \eqref{eq:h1constraints}.
To extract inadmissible block $\mA_{\alpha, \beta}$ of level $L$,
we require a test matrix $\mOmega$ that satisfies the following
sampling constraints.
\begin{equation*}
\begin{aligned}
    \mOmega(I_\beta, :)  &= \mI & & \\
    \mOmega(I_\gamma, :) &= \mzero &\text{for all }& \gamma \in \alphanear \setminus \{\beta\}
\end{aligned}
\end{equation*}
If test matrix $\mOmega$ satisfies these constraints,
then the rows of the product $\mA \mOmega - \mA^{(2)} \mOmega$ indexed by $I_\alpha$
will contain $\mA_{\alpha, \beta}$.
The graph corresponding to the 22 inadmissible blocks belonging to level 3
is depicted in \fref{fig:graph_strong_1d_leaf},
and its coloring produces test matrices with the following structures.

\begin{equation*}
    \begin{bmatrix}
        \mOmega_1 & \mOmega_2 & \mOmega_3
    \end{bmatrix}
    =
    \begin{bmatrix}
        \mI_{8}  & \mzero   & \mzero   \\
        \mzero   & \mI_{9}  & \mzero   \\
        \mzero   & \mzero   & \mI_{10} \\
        \mI_{11} & \mzero   & \mzero   \\
        \mzero   & \mI_{12} & \mzero   \\
        \mzero   & \mzero   & \mI_{13} \\
        \mI_{14} & \mzero   & \mzero   \\
        \mzero   & \mI_{15} & \mzero   \\
    \end{bmatrix}
\end{equation*}

The entire process of compressing an $\h^1$ matrix
is summarized in \cref{alg:h1_pseudocode}.

\begin{algorithm}
\caption{Randomized compression of an $\h^1$ matrix}
\label{alg:h1_pseudocode}
\begin{algorithmic}
    \For{level $l \in [2, ..., L]$}
        \State{\underline{Compute randomized samples of $\mA - \mA^{(l)}$}}
        \State{Construct structured random test matrices $\{\mOmega_i\}$
            of size $N \times (k + p)$ as in \cref{sec:h1_level3} or \cref{sec:h1general}}
        \ForAll{$\mOmega_i$}
            \State{Multiply $\mY_i = \mA \mOmega_i - \mA^{(l)} \mOmega_i$}
        \EndFor
        \State
        
        \State{\underline{Compute orthonormal basis matrices $\ulong_{\alpha, \beta}$}}
        \ForAll{interacting pairs $\alpha, \beta$ in level $l$}
            \State{Identify $\mY_i$ that contains a sample of $\mA_{\alpha, \beta}$
                in $\mY_i(I_\alpha, :)$}
            \State{Orthonormalize $\ulong_{\alpha, \beta}
                = \qr(\mY_i(I_\alpha, :), k)$}
        \EndFor
        \State
        
        \State{\underline{Compute randomized samples of
            $\mA^* - {\mA^{(l)}}^*$}}
        \State{Construct structured random test matrices $\{\mPsi_i\}$
            of size $N \times (k + p)$ as in \cref{sec:h1_level3} or \cref{sec:h1general}
            (if the rank structure of $\mA$ is symmetric,
            then the test matrices $\{\mOmega_i\}$ can be reused)}
        \ForAll{$\mPsi_i$}
            \State{Multiply $\mZ_i = \mA^* \mPsi_i - {\mA^{(l)}}^* \mPsi_i$}
        \EndFor
        \State{}
        
        \State{\underline{Compute orthonormal basis matrices $\vlong_{\alpha, \beta}$}}
        \ForAll{interacting pairs $\alpha, \beta$ in level $l$}
            \State{Identify $\mZ_i$ that contains a sample of $\mA_{\alpha, \beta}^*$
                in $\mZ_i(I_\beta, :)$}
            \State{Orthonormalize $\vlong_{\alpha, \beta}
                = \qr(\mZ_i(I_\beta, :), k)$}
        \EndFor
        \State
        
        \State{\underline{Solve for $\mB$}}
        \ForAll{interacting pairs $\alpha, \beta$ in level $l$}
            \State{$\mB_{\alpha, \beta}
                = \left( \mG_\alpha \ulong_\alpha \right)^\dagger
                \mG_\alpha \mA_{\alpha, \beta} \mG_\beta
                \left( \vlong_\beta \mG_\beta \right)^\dagger$}
        \EndFor
    \EndFor
    \State
    
    \State{\underline{Extract the inadmissible blocks of level $L$}}
    \State{Construct structured random test matrices $\{\mOmega_i\}$ of size
        $N \times m$
        as in \cref{sec:leaf} or \cref{sec:h1general}}
    \ForAll{$\mOmega_i$}
        \State{Multiply $\mY_i = \mA \mOmega_i - \mA^{(L)} \mOmega_i$}
    \EndFor
    \ForAll{neighbor pairs $\alpha, \beta$ in level $L$}
        \State{Identify $\mY_i$ that contains $\mA_{\alpha, \beta}$}
        \State{Extract the block $\mA_{\alpha, \beta}$
            from $\mY_i(I_\alpha, :)$}
    \EndFor
\end{algorithmic}
\end{algorithm}

\subsubsection{General patterns in the test matrices for \texorpdfstring{$\h^1$}{H1} compression}
\label{sec:h1general}

In this section, we describe sets of test matrices
that are sufficient for compressing
any $\h^1$ matrix based on points in arbitrary dimensions.
Though sufficient, these sets of test matrices
are not necessarily the smallest possible,
and the proposed method of constructing problem-specific test matrices
via graph coloring often yields a smaller set of test matrices,
resulting in fewer matrix-vector products.
The test matrices described in this section
establish an upper bound on the number of matrix-vector products
that improves on the number given in \cite{lin2011fast},
and they also imply a bound on the chromatic number of the graph,
which appears in the estimate of asymptotic complexity.

\paragraph{A general set of test matrices for sampling admissible blocks}
The structures of the test matrices derived in \cref{sec:h1_level3,sec:leaf}
exhibit patterns that generalize to finer and higher-dimensional grids.
The sampling constraints~\cref{eq:h1constraints}
for admissible block $\mA_{\alpha, \beta}$
specify that the rows of the test matrix corresponding to $\beta$
must be filled with random values,
and the rows of the test matrix
corresponding to the other boxes in $\alphanear \cup \alphainteraction$
must be filled with zeros.
For a 1-dimensional problem,
the indices corresponding to $\alphanear \cup \alphainteraction$
form 6 contiguous blocks of the test matrix.
Accordingly, each of the test matrices defined in \cref{eq:h1testmatrices}
has every sixth block filled with random values
and the other blocks filled with zeros,
a pattern that generalizes to an arbitrarily fine grid in one dimension.

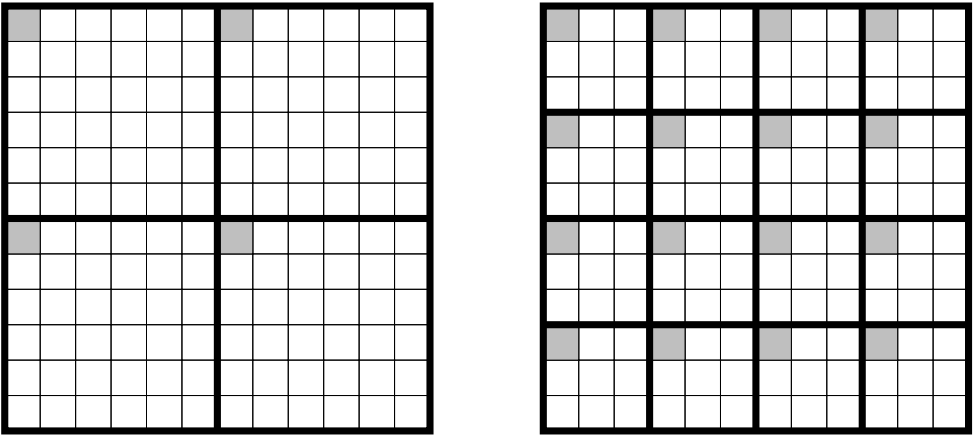
\begin{figure}[tb]
    \centering

    \begin{subfigure}{0.45\textwidth}
    \centering
    \resizebox{\textwidth}{!}{
    \begin{tikzpicture}
        [box/.style={rectangle,draw=black,minimum size=1cm}]

    \foreach \x in {0,1,...,11}{
        \foreach \y in {0,1,...,11}
            \node[box] at (\x,\y){};
    }

    \node[box,fill=lightgray] at (0, 11){};
    \node[box,fill=lightgray] at (6, 11){};
    \node[box,fill=lightgray] at (0, 5){};
    \node[box,fill=lightgray] at (6, 5){};

    \node[box,minimum size=12cm,line width=2.0mm] at (5.5, 5.5){};
    \draw[line width=2.0mm] (5.5, -0.5) -- (5.5, 11.5);
    \draw[line width=2.0mm] (-0.5, 5.5) -- (11.5, 5.5);

    \end{tikzpicture}
    }
    \end{subfigure}
    \hfill
    \begin{subfigure}{0.45\textwidth}
    \centering
    \resizebox{\textwidth}{!}{
    \begin{tikzpicture}
        [box/.style={rectangle,draw=black,minimum size=1cm}]

    \foreach \x in {0,1,...,11}{
        \foreach \y in {0,1,...,11}
            \node[box] at (\x,\y){};
    }

    \foreach \x in {0, 3, 6, 9} {
        \foreach \y in {2, 5, 8, 11}
            \node[box,fill=lightgray] at (\x, \y){};
    }

    \node[box,minimum size=12cm,line width=2.0mm] at (5.5, 5.5){};
    \draw[line width=2.0mm] (5.5, -0.5) -- (5.5, 11.5);
    \draw[line width=2.0mm] (-0.5, 5.5) -- (11.5, 5.5);
    \draw[line width=2.0mm] (2.5, -0.5) -- (2.5, 11.5);
    \draw[line width=2.0mm] (-0.5, 2.5) -- (11.5, 2.5);
    \draw[line width=2.0mm] (8.5, -0.5) -- (8.5, 11.5);
    \draw[line width=2.0mm] (-0.5, 8.5) -- (11.5, 8.5);

    \end{tikzpicture}
    }
    \end{subfigure}
    \caption{
        Patterns representing a test matrix for sampling admissible blocks (left)
        and a test matrix for sampling inadmissible blocks (right)
        for a problem in 2 dimensions.
        Blocks of the test matrices corresponding to gray boxes are filled with random values,
        and those corresponding to white boxes are filled with zeros.
        The other test matrices are obtained by shifting these pattern
        horizontally and vertically.
    }
    \label{fig:h1activations}
\end{figure}

To describe a general set of test matrices for a problem in 2 dimensions,
we partition the domain
into tiles, each of which covers $6 \times 6$ boxes.
We define 36 test matrices,
each activating a set of
boxes that share the same position
within their respective tiles.
The set of boxes activated by one of the test matrices
is shown in \cref{fig:h1activations}.
The sampling constraints in \cref{eq:h1constraints}
apply to $\alphanear \cup \alphainteraction$,
which form a square of $6 \times 6$ boxes,
and they specify that exactly one of those boxes must be activated.
Therefore, those constraints will be satisfied by one of the 36 test matrices.
By a similar argument, this pattern generalizes to higher-dimensional problems,
requiring at most $6^d$ test matrices
to sample one level of admissible blocks for a problem in $d$ dimensions.
Furthermore, since these test matrices satisfy \cref{eq:h1constraints},
they must correspond to a valid coloring
of the graph described in \cref{sec:h1_level3},
bounding the chromatic number of the graph by
\begin{equation*}
    \colorsnonunif \leq 6^d.
\end{equation*}

\paragraph{A general set of test matrices for extracting inadmissible blocks}
A similar argument proves that
extracting inadmissible blocks of the leaf level
can be accomplished using $3^d$ test matrices.
The set of boxes activated by one of the test matrices
for a problem in two dimensions
is shown in \cref{fig:h1activations}.
Therefore, the chromatic number
of the graph described in \cref{sec:leaf}
is bounded by
\begin{equation*}
    \colorsnear \leq 3^d.
\end{equation*}

\subsubsection{Asymptotic complexity}
\label{sec:h1_complexity}

Let
$L \sim \log{N}$ denote the depth of the tree,
$\tmult$ denote the time to apply $\mA$ or $\mA^*$ to a vector,
and
$\tflop$ denote the time to execute one floating point operation.
There are $2^{dl}$ boxes belonging to level $l$ of the tree,
and the interaction list of each box consists of up to $6^d - 3^d$ other boxes,
so there are approximately $(6^d - 3^d) 2^{dl}$ admissible blocks
associated with level $l$.
The cost of applying the low-rank approximation of an admissible block
associated with level $l$ to a vector is $\sim k N / 2^{dl}$ operations.
Therefore, the cost of applying the level-$l$ truncated matrix
$\mA^{(l)}$ (or its transpose) to a vector is
\begin{equation*}
    \tAl
    \sim
    \tflop \times \sum_{j=0}^l (6^d - 3^d) 2^{dj} k \frac{N}{2^{dj}}
    \sim
    \tflop \times (6^d - 3^d) l k N.
\end{equation*}

The cost of compressing the admissible blocks associated with level $l$ of the tree is
\begin{equation*}
    \tl
    \sim
    \left( \tmult + \tAl \right) \times 2 \colorsnonunif k
    +
    \tflop \times (6^d - 3^d) 2^{dl} k^2 \frac{N}{2^{dl}},
\end{equation*}
since we require $\sim \colorsnonunif k$ applications of $\mA - \mA^{(l)}$ and its transpose,
where $\colorsnonunif$ denotes the chromatic number
of the graph described in \cref{sec:h1_level3},
and we require an additional
$\sim k^2 N / 2^{dl}$ operations for each admissible block
to compute low-rank approximations.

Finally, summing $\tl$ over each level
gives the total time to construct an $\h^1$ representation as follows.
The cost of extracting inadmissible blocks
associated with the leaf level of the tree is omitted as it
only contributes a lower order term to the overall cost.
We also omit the costs associated with the graph coloring problem
since it is only worthwhile for problems that exhibit low-dimensional structure,
in which case the costs would be much lower than the worst-case analysis would suggest,
and we observe in practice that the cost of graph coloring
represents a small portion of the total cost of compression.
For problems without low-dimensional structure,
one may skip the graph coloring step
and instead use the general set of test matrices described in \cref{sec:h1general}.

\begin{equation*}
    \label{eq:h1complexity}
    \tcompress \sim \tmult \times 2 \colorsnonunif k \log{N} \\
    + \tflop \times (6^d - 3^d) 2 \colorsnonunif k^2 N \left( \log{N} \right)^2
\end{equation*}

\subsection{Uniform \texorpdfstring{$\h^1$}{H1} matrix compression}
\label{sec:unif_h1}

A uniform $\h^1$ approximation requires for each box $\tau$
a column basis matrix $\ulong_\tau$
and a row basis matrix $\vlong_\tau$
such that
\begin{equation*}
    \mA_{\alpha, \beta} \approx \ulong_\alpha \ulong_\alpha^* \mA_{\alpha, \beta} \vlong_\beta \vlong_\beta^*
\end{equation*}
for all admissible pairs $(\alpha, \beta)$.
In other words, $\ulong_\alpha$ must span the column space of
$\mA \left( I_\alpha, \cup_{\beta \in \alphainteraction} I_\beta \right)$,
the submatrix of interactions between $\alpha$ and the boxes in its interaction list.
Similarly, $\vlong_\beta$ must span the row space of
$\mA \left( \cup_{\alpha \in \betainteraction} I_\alpha, I_\beta \right)$.
The algorithms for compressing $\h^1$ and uniform $\h^1$ matrices
have a two important differences.
For a uniform $\h^1$ matrix,
we use \cref{alg:compress_two_stage},
rather than \cref{alg:compress_double_random}
to compute low-rank approximations.
Also, instead of sampling the column space
of each admissible block $\mA_{\alpha, \beta}$ separately,
we will sample the interactions between each box
and all of the boxes in its interaction list together.

We return to the task of compressing the admissible blocks
associated with level 3 of the tree as described in \cref{sec:h1_level3}.
As before, we will sample $\mA - \mA^{(2)}$
with a set of test matrices.
To sample those interactions for some box $\alpha$,
we require a test matrix $\mOmega$ that satisfies the following
sampling constraints.

\begin{equation}\label{eq:unifh1constraints}
\begin{aligned}
    \mOmega(I_\beta, :)  &= \mG_\beta &\text{for all }& \beta \in \alphainteraction \\
    \mOmega(I_\gamma, :) &= \mzero    &\text{for all }& \gamma \in \alphanear
\end{aligned}
\end{equation}

If test matrix $\mOmega$ satisfies the above sampling constraints,
then the rows of the product $\mA \mOmega - \mA^{(2)} \mOmega$ indexed by $I_\alpha$
will contain a randomized sample of the column space
of $\mA \left( I_\alpha, \cup_{\beta \in \alphainteraction} I_\beta \right)$,
the submatrix of interactions between $\alpha$ and the boxes in its interaction list.

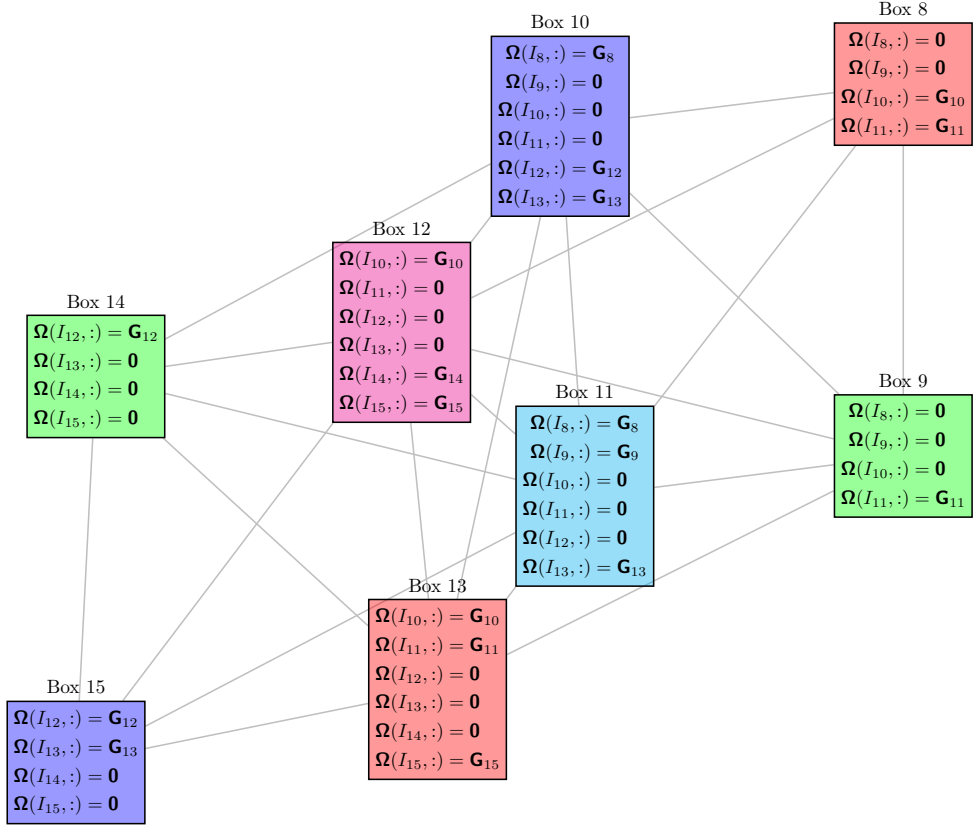
\begin{figure}[tb]
    \centering
    \resizebox{\columnwidth}{!}{
    \begin{tikzpicture}[spring electrical layout, electric charge=60]
    \node (8)  [vertexstyle, label={Box 8}, fill=red] {
            $\begin{aligned}
                \mOmega(I_{8}, :) &= \mzero \\
                \mOmega(I_{9}, :) &= \mzero \\
                \mOmega(I_{10}, :) &= \mG_{10} \\
                \mOmega(I_{11}, :) &= \mG_{11} \\
            \end{aligned}$
        };
    \node (9)  [vertexstyle, label={Box 9}, fill=green] {
            $\begin{aligned}
                \mOmega(I_{8}, :) &= \mzero \\
                \mOmega(I_{9}, :) &= \mzero \\
                \mOmega(I_{10}, :) &= \mzero \\
                \mOmega(I_{11}, :) &= \mG_{11} \\
            \end{aligned}$
        };
    \node (10)  [vertexstyle, label={Box 10}, fill=blue] {
            $\begin{aligned}
                \mOmega(I_{8}, :) &= \mG_{8} \\
                \mOmega(I_{9}, :) &= \mzero \\
                \mOmega(I_{10}, :) &= \mzero \\
                \mOmega(I_{11}, :) &= \mzero \\
                \mOmega(I_{12}, :) &= \mG_{12} \\
                \mOmega(I_{13}, :) &= \mG_{13} \\
            \end{aligned}$
        };
    \node (11)  [vertexstyle, label={Box 11}, fill=cyan] {
            $\begin{aligned}
                \mOmega(I_{8}, :) &= \mG_{8} \\
                \mOmega(I_{9}, :) &= \mG_{9} \\
                \mOmega(I_{10}, :) &= \mzero \\
                \mOmega(I_{11}, :) &= \mzero \\
                \mOmega(I_{12}, :) &= \mzero \\
                \mOmega(I_{13}, :) &= \mG_{13} \\
            \end{aligned}$
        };
    \node (12)  [vertexstyle, label={Box 12}, fill=magenta] {
            $\begin{aligned}
                \mOmega(I_{10}, :) &= \mG_{10} \\
                \mOmega(I_{11}, :) &= \mzero \\
                \mOmega(I_{12}, :) &= \mzero \\
                \mOmega(I_{13}, :) &= \mzero \\
                \mOmega(I_{14}, :) &= \mG_{14} \\
                \mOmega(I_{15}, :) &= \mG_{15} \\
            \end{aligned}$
        };
    \node (13)  [vertexstyle, label={Box 13}, fill=red] {
            $\begin{aligned}
                \mOmega(I_{10}, :) &= \mG_{10} \\
                \mOmega(I_{11}, :) &= \mG_{11} \\
                \mOmega(I_{12}, :) &= \mzero \\
                \mOmega(I_{13}, :) &= \mzero \\
                \mOmega(I_{14}, :) &= \mzero \\
                \mOmega(I_{15}, :) &= \mG_{15} \\
            \end{aligned}$
        };
    \node (14)  [vertexstyle, label={Box 14}, fill=green] {
            $\begin{aligned}
                \mOmega(I_{12}, :) &= \mG_{12} \\
                \mOmega(I_{13}, :) &= \mzero \\
                \mOmega(I_{14}, :) &= \mzero \\
                \mOmega(I_{15}, :) &= \mzero \\
            \end{aligned}$
        };
    \node (15)  [vertexstyle, label={Box 15}, fill=blue] {
            $\begin{aligned}
                \mOmega(I_{12}, :) &= \mG_{12} \\
                \mOmega(I_{13}, :) &= \mG_{13} \\
                \mOmega(I_{14}, :) &= \mzero \\
                \mOmega(I_{15}, :) &= \mzero \\
            \end{aligned}$
        };

    \path (8) edge[edgestyle] (9);
    \path (8) edge[edgestyle] (10);
    \path (8) edge[edgestyle] (11);
    \path (8) edge[edgestyle] (12);

    \path (9) edge[edgestyle] (10);
    \path (9) edge[edgestyle] (11);
    \path (9) edge[edgestyle] (12);
    \path (9) edge[edgestyle] (13);

    \path (10) edge[edgestyle] (11);
    \path (10) edge[edgestyle] (12);
    \path (10) edge[edgestyle] (13);
    \path (10) edge[edgestyle] (14);

    \path (11) edge[edgestyle] (12);
    \path (11) edge[edgestyle] (13);
    \path (11) edge[edgestyle] (14);
    \path (11) edge[edgestyle] (15);

    \path (12) edge[edgestyle] (13);
    \path (12) edge[edgestyle] (14);
    \path (12) edge[edgestyle] (15);

    \path (13) edge[edgestyle] (14);
    \path (13) edge[edgestyle] (15);

    \path (14) edge[edgestyle] (15);

    \end{tikzpicture}
    }
    \caption{The constraint incompatibility graph corresponding to the 8 boxes
    on level 3 of the matrix shown in \fref{fig:matrix_strong_1d_lvl3}
    along with the constraints \eqref{eq:unifh1constraints} associated with sampling a uniform $\h^1$ matrix.}
    \label{fig:graph_strong_1d_lvl3_unif}
\end{figure}

We have one set of sampling constraints for each box $\alpha$,
and we use them to form the constraint incompatibility graph
(\cref{def:incompatibility_graph})
shown in \cref{fig:graph_strong_1d_lvl3_unif}.
As in \cref{sec:h1_level3},
we compute a vertex coloring of the graph.
The coloring of the graph in \cref{fig:graph_strong_1d_lvl3_unif}
specifies test matrices with the following structures.
\begin{equation*}
    \begin{bmatrix}
        \mOmega_1 & \mOmega_2 & \mOmega_3 & \mOmega_4 & \mOmega_5
    \end{bmatrix}
    =
    \begin{bmatrix}
        \mzero   & \mzero   & \mG_{8}  & \mG_{8}  & *      \\
        \mzero   & \mzero   & \mzero   & \mG_{9}  & *      \\
        \mG_{10} & \mzero   & \mzero   & \mzero   & \mG_{10} \\
        \mG_{11} & \mG_{11} & \mzero   & \mzero   & \mzero \\
        \mzero   & \mG_{12} & \mG_{12} & \mzero   & \mzero \\
        \mzero   & \mzero   & \mG_{13} & \mG_{13} & \mzero \\
        \mzero   & \mzero   & \mzero   & *        & \mG_{14} \\
        \mG_{15} & \mzero   & \mzero   & *        & \mG_{15} \\
    \end{bmatrix}
\end{equation*}

We evaluate the samples $\mY_i = \mA \mOmega_i - \mA^{(2)} \mOmega_i$
for each test matrix $\mOmega_i$,
and orthonormalize the relevant blocks of the sample matrices $\mY_i$ to obtain
orthonormal column basis matrix
$\ulong_\alpha$ for each box $\alpha$.
That is, if $\mOmega_i$ is the test matrix
satisfying the sampling constraints of box $\alpha$,
then
\begin{equation*}
    \mY_i(I_\alpha, :) = \sum_{\beta \in \alphainteraction} \mA_{\alpha, \beta} \mG_\beta,
\end{equation*}
so we compute
\begin{equation*}
    \ulong_\alpha = \qr \left( \mY_i(I_\alpha, :) \right).
\end{equation*}

For the second stage of compression using \cref{alg:compress_two_stage},
we must compute $\mA_{\alpha, \beta}^* \ulong_\alpha$
for each admissible block $\mA_{\alpha, \beta}$.
To do so, we use the procedure for sampling an $\h^1$ matrix
with non-uniform basis matrices (\cref{sec:h1_level3}),
but rather than filling blocks of the test matrices with random values,
we fill them with appropriately chosen uniform basis matrices.
With that modification, we have the following sampling constraints
for each admissible pair $(\alpha, \beta)$.
\begin{equation*}
\begin{aligned}
    \mPsi(I_\beta, :)  &= \ulong_\beta & & \\
    \mPsi(I_\gamma, :) &= \mzero &\text{for all }& \gamma \in \alphanear \cup \alphainteraction \setminus \{\beta\}
\end{aligned}
\end{equation*}

We construct a suitable set of test matrices via graph coloring
and evaluate the products
$\mZ_i = \left( \mA - \mA^{(2)} \right)^* \mPsi_i$
for each test matrix $\mPsi_i$.
Then for each admissible pair $(\alpha, \beta)$,
there is some $\mZ_i$ such that
$\mZ_i(I_\beta, :) = \mA_{\alpha, \beta}^* \ulong_\alpha$.
With those results,
we compute the uniform basis matrices $\vlong_\alpha$
and the matrices $\mB_{\alpha, \beta}$ as follows.
\begin{gather*}
    \vlong_\beta
        = \qr \left( \sum_{\alpha \in \betainteraction} \mZ_i(I_\beta, :)  \right)
        = \qr \left( \sum_{\alpha \in \betainteraction} \mA_{\alpha, \beta}^* \ulong_\alpha \right) \\
    \mB_{\alpha, \beta}
        = \mZ_i(I_\beta, :)^* \vlong_\beta
        = (\mA_{\alpha, \beta}^* \ulong_\alpha)^* \vlong_\beta
\end{gather*}

Note that it would have been possible
to construct $\vlong_\beta$
by taking combined samples of the interactions
between $\beta$ and all of the boxes in its interaction list,
similarly to how we obtain $\ulong_\beta$.
However, then we would not have
$\mA_{\alpha, \beta}^* \ulong_\alpha$,
which is later required when computing $\mB_{\alpha, \beta}$.
The advantage of sampling $\mA^*$ the way we do
is that we obtain what is needed
to compute both $\vlong_\beta$ and $\mB_{\alpha, \beta}$.

\begin{remark}
An alternative approach for uniform $\h^1$ compression~\cite{lin2011fast}
is obtained by adding a uniformization step
to the $\h^1$ compression algorithm.
With that approach, one first computes
low-rank approximations
$\mA_{\alpha, \beta}
=
\ulong_{\alpha, \beta} \mB_{\alpha, \beta} \vlong_{\alpha, \beta}^*$
with non-uniform basis matrices.
Those low-rank approximations are recompressed
to obtain uniform basis matrices,
followed by a change of basis.
\begin{equation*}
    \mA_{\alpha, \beta}
    \approx
    \ulong_\alpha \left(
    \ulong_\alpha^* \ulong_{\alpha, \beta} \mB_{\alpha, \beta} \vlong_{\alpha, \beta}^* \vlong_\beta
    \right) \vlong_\beta^*
\end{equation*}
One disadvantage of that approach
is that $\ulong_\alpha$ and $\vlong_\alpha$ are based on samples of the factorizations
$\ulong_{\alpha, \beta} \mB_{\alpha, \beta} \vlong_{\alpha, \beta}^*$,
which may only be approximate,
rather than direct samples of the blocks $\mA_{\alpha, \beta}$ themselves,
potentially resulting in lower accuracy.
Another difference is that the procedure for sampling a uniform $\h^1$ matrix
presented in this section requires fewer matrix-vector products
than the procedure for sampling an $\h^1$ matrix,
as described in \cref{sec:h1general,sec:unifh1chromatic}.
\end{remark}

\subsubsection{General patterns in the test matrices for uniform \texorpdfstring{$\h^1$}{H1} compression}
\label{sec:unifh1chromatic}

\begin{figure}[tb]
    \centering

    \begin{subfigure}{0.45\textwidth}
    \centering
    \resizebox{\textwidth}{!}{
    \begin{tikzpicture}
        [box/.style={rectangle,draw=black,minimum size=1cm}]

    \foreach \x in {0,1,...,14}{
        \foreach \y in {0,1,...,14}
            \node[box,fill=lightgray] at (\x,\y){};
    }

    \foreach \x in {0,1,2}{
        \foreach \y in {2,3,4}{
            \foreach \dx in {0, 5, 10}{
                \foreach \dy in {0, 5, 10}{
                    \node[box,fill=white] at (\x + \dx, \y + \dy){};
                }
            }
        }
    }

    \node[box,minimum size=15cm,line width=2.0mm] at (7, 7){};

    \end{tikzpicture}
    }
    \end{subfigure}
    \hfill
    \caption{
        A pattern representing one test matrix
        for the first stage of uniform $\h^1$ sampling for a problem in 2 dimensions.
        Blocks of the test matrix corresponding to gray boxes are filled with random values,
        and those corresponding to white boxes are filled with zeros.
        The other test matrices are obtained by shifting this pattern
        horizontally and vertically.
    }
    \label{fig:unifh1activations}
\end{figure}
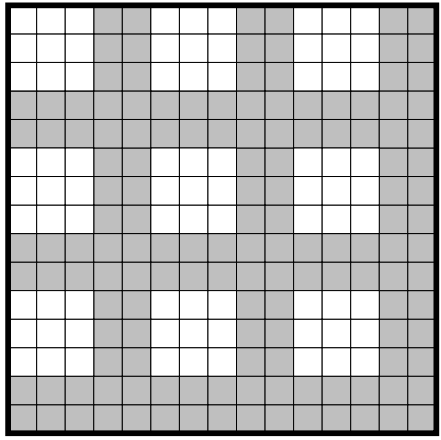

As in \cref{sec:h1general},
we describe a general set of test matrices
that is applicable to finer and higher-dimensional grids.
To sample the interactions of box $\alpha$
for a problem in 2 dimensions,
the blocks of the test matrix
corresponding to the boxes in $\alphanear$ must be filled with zeros,
and the blocks corresponding to the boxes in $\alphainteraction$
must be filled with random values \cref{eq:unifh1constraints}.
The pattern of activated boxes for one test matrix is shown in \cref{fig:unifh1activations}.
The complete set of 25 test matrices
is obtained by shifting the pattern
horizontally and vertically.
The pattern generalizes to higher-dimensional problems,
requiring at most $5^d$ test matrices
to carry out the first stage of uniform $\h^1$ sampling
for a problem in $d$ dimensions.
Therefore, we establish the upper bound
\begin{equation*}
    \colorsunif \leq 5^d.
\end{equation*}

\subsection{\texorpdfstring{$\h^2$}{H2} matrix compression}
\label{sec:h2}

To obtain an algorithm for computing an $\h^2$ representation,
we apply three modifications to the algorithm for computing a uniform $\h^1$ representation.
First, we take singular value decompositions of samples,
rather than orthonormalizing via QR factorizations.
Second, we enrich the basis matrices of each box
so that they span the relevant parts of the basis matrices of its parent.
For box $\alpha$ with parent $\tau$,
we define $I_{\tau, \alpha}$ to be the positions of the indices of $I_\alpha$ within $I_\tau$
so that $I_\tau(I_{\tau, \alpha}) = I_\alpha$.
When computing the basis matrix $\ulong_\alpha$,
we augment $\mY_i(I_\alpha, :)$,
the sample of $\mA \left( I_\alpha, \cup_{\beta \in \alphainteraction} I_\beta \right)$,
with $\ulong_\tau(I_{\tau, \alpha}, :) \mS_\tau^\text{in} \mG_k$,
a sample of $\mA \left( I_\alpha, \left( \cup_{\beta \in \taunear} I_\beta \right)^c \right)$,
where $\mG_k$ is a $k$-by-$k$ Gaussian random matrix.
\begin{equation*}
\begin{aligned}
    \mY_i(I_\alpha, ;) &= \sum_{\beta \in \alphainteraction} \mA_{\alpha, \beta} \mG_\beta \\
    \lbrack \ulong_\alpha, \mS_\alpha^\text{in}, \sim \rbrack &= \svd \left( \mY_i(I_\alpha, :) + \ulong_\tau(I_{\tau, \alpha}, :) \mS_\tau^\text{in} \mG_k \right) \\
\end{aligned}
\end{equation*}
We use a similarly augmented sample when computing $\vlong_\alpha$.
\begin{equation*}
\begin{aligned}
    \lbrack \vlong_\alpha, \mS_\alpha^\text{out}, \sim \rbrack &= \svd \left( \left( \sum_{\beta \in \alphainteraction} \mA_{\beta, \alpha}^* \ulong_\beta \right) + \vlong_\tau(I_{\tau, \alpha}, :) \mS_\tau^\text{out} \mG_k \right) \\
\end{aligned}
\end{equation*}
For some boxes, there may be no contribution from the parent with which to augment the sample.
For example, if $\alpha$ is a box belonging to level 2,
then $\left( \cup_{\beta \in \taunear} I_\beta \right)^c$ is the empty set.
In such cases, we compute the SVD without augmenting.
Third, after we have computed long basis matrices of the children of $\tau$,
we compute the short basis matrix $\ushort_\tau$ by solving~\eqref{eq:nested_basis},
and then we can discard $\ulong_\tau$.
Similarly, we compute $\vshort_\tau$ and discard $\vlong_\tau$.
The processes of compressing uniform $\h^1$ and $\h^2$ matrices
are summarized in \cref{alg:h2_pseudocode}.
The algorithm for applying a level-truncated approximation
of an $\h^2$ matrix to a vector
is given in \cref{alg:apply_h2_level_truncated}.

\begin{algorithm}
\caption{Randomized compression of a uniform $\h^1$ or $\h^2$ matrix}
\label{alg:h2_pseudocode}
\begin{algorithmic}
    \For{level $l \in [2, ..., L]$}
        \State{\underline{Compute randomized samples of $\mA - \mA^{(l)}$}}
        \State{Construct structured random test matrices $\{\mOmega_i\}$
            of size $N \times (k + p)$ as in \cref{sec:unif_h1} or \cref{sec:unifh1chromatic}}
        \ForAll{$\mOmega_i$}
            \State{Multiply $\mY_i = \mA \mOmega_i - \mA^{(l)} \mOmega_i$}
        \EndFor
        \State
        
        \State{\underline{Compute uniform orthonormal basis matrices $\ulong_\alpha$}}
        \ForAll{boxes $\alpha$ in level $l$}
            \State{Identify $\mY_i$ that contains a sample of
                $\mA(I_\alpha, \cup_{\beta \in \alphainteraction} I_\beta)$
                in $\mY_i(I_\alpha, :)$}
            \If{$\mA$ has uniform $\h^1$ structure}
                \State{Orthonormalize the sample:
                    $\ulong_\alpha = \qr(\mY_i(I_\alpha, :), k)$}
            \ElsIf{$\mA$ has $\h^2$ structure}
                \State{With $\tau$ as the parent of $\alpha$, compute an SVD of the augmented sample:}
                \State{
                    $\lbrack \ulong_\alpha, \mS_\alpha^\text{in}, \sim \rbrack
                    = \svd(\mY_i(I_\alpha, :) + \ulong_\tau(I_{\tau, \alpha}, :) \mS_\tau^\text{in} \mG_k, k)$}
            \EndIf
        \EndFor
        \State
        
        \State{\underline{Compute randomized samples of
            $\mA^* - {\mA^{(l)}}^*$}}
        \State{Construct structured random test matrices $\{\mPsi_i\}$
            of size $N \times (k + p)$ as in \cref{sec:h1_level3} or \cref{sec:h1general}
            with the modification in \cref{sec:unif_h1}}
        \ForAll{$\mPsi_i$}
            \State{Multiply $\mZ_i = \mA^* \mPsi_i - {\mA^{(l)}}^* \mPsi_i$}
        \EndFor
        \State{}
        
        \State{\underline{Compute uniform orthonormal basis matrices $\vlong_\beta$}}
        \ForAll{boxes $\beta$ in level $l$}
            \State{Identify $\mZ_i$ that contains
                $\mA_{\alpha, \beta}^* \ulong_\alpha$
                in $\mZ_i(I_\beta, :)$ for $\alpha \in \betainteraction$}
            \If{$\mA$ has uniform $\h^1$ structure}
                \State{Sum samples and orthonormalize:
                    $\vlong_\beta = \qr(\sum_{\alpha \in \betainteraction} \mZ_i(I_\beta, :), k)$
                }
            \ElsIf{$\mA$ has $\h^2$ structure}
                \State{With $\tau$ as the parent of $\beta$, sum samples, augment, and compute an SVD:}
                \State{
                    $\lbrack \vlong_\beta, \mS_\beta^\text{out}, \sim \rbrack
                    = \svd((\sum_{\alpha \in \betainteraction} \mZ_i(I_\beta, :)) + \vlong_\tau(I_{\tau, \beta}, :) \mS_\tau^\text{out} \mG_k, k)$
                }
            \EndIf
        \EndFor
        \State
        
        \State{\underline{Compute $\mB_{\alpha, \beta}$}}
        \ForAll{interacting pairs $\alpha, \beta$ in level $l$}
            \State{Identify $\mZ_i$ that contains
                $\mA_{\alpha, \beta}^* \ulong_\alpha$
                in $\mZ_i(I_\beta, :)$ for $\alpha \in \betainteraction$}
            \State{$\mB_{\alpha, \beta}
                = \mZ_i^*(I_\beta, :) \vlong_\beta$}
        \EndFor
        \State
        
        \State{\underline{Nest the basis matrices}}
        \If{$\mA$ has $\h^2$ structure}
            \ForAll{boxes $\alpha$ in level $l$}
                \State{Compute $\ushort_\alpha, \vshort_\alpha$
                    using \cref{eq:nested_basis}}
                \State{Discard $\ulong_\alpha, \vlong_\alpha$}
            \EndFor
        \EndIf
    \EndFor
    \State
    
    \State{Extract the inadmissible blocks of level $L$
        in the same way as in \mbox{\cref{alg:h1_pseudocode}}}
\end{algorithmic}
\end{algorithm}

\begin{algorithm}
\caption{Applying an $\h^2$ level-truncated approximation $\mA^{(l)}$ to vector $\qq$}
\label{alg:apply_h2_level_truncated}
\begin{algorithmic}
    \State{\underline{Build outgoing expansions in level $l$.}}
    \ForAll{boxes $\tau$ in level $l$}
        \State{$\hat{\qq}_\tau = \vlong_\tau^* \qq(I_\tau)$}
    \EndFor
    \State
    
    \State{\underline{Build outgoing expansions
        for levels coarser than $l$ (upward pass).}}
    \ForAll{levels $i \in [l - 1, ... 2]$}
        \ForAll{boxes $\tau$ in level $i$}
            \ForAll{children $\alpha$ of $\tau$}
                \State{$\hat{\qq}_\tau (I_\alpha, :)
                    = \vshort_\tau^*(I_\alpha, :) \hat{\qq}_\alpha$}
            \EndFor
        \EndFor
    \EndFor
    \State
    
    \State{\underline{Build incoming expansions for boxes in level 2.}}
    \ForAll{boxes $\alpha$ in level 2}
        \State{$\hat{\uu}_\alpha
            = \sum_{\beta \in \alphainteraction} \mB_{\alpha, \beta} \hat{\qq}_\beta$}
    \EndFor
    \State
    
    \State{\underline{Build incoming expansions
        for levels finer than 2 (downward pass).}}
    \ForAll{levels $i \in [3, ..., l - 1]$}
        \ForAll{boxes $\alpha$ in level $i$}
            \State{Let $\tau$ be the parent of $\alpha$}
            \State{$\hat{\uu}_\alpha
            = \sum_{\beta \in \alphainteraction} \mB_{\alpha, \beta} \hat{\qq}_\beta + \ushort_\tau(I_\alpha, :) \hat{\uu}_\tau$}
        \EndFor
    \EndFor
    \State
    
    \State{\underline{Build incoming expansions for level $l$.}}
    \ForAll{boxes $\alpha$ in level $l$}
        \State{$\uu(I_\alpha)
            = \ulong_\tau \hat{\uu}_\tau
            + \sum_{\beta \in \alphanear} \mA(I_\alpha, I_\beta) \qq(I_\beta)$}
    \EndFor
\end{algorithmic}
\end{algorithm}

\subsubsection{Asymptotic complexity}

The asymptotic complexity of the $\h^2$ algorithm
is similar to that of the $\h^1$ algorithm.
Since we are now using uniform basis matrices,
the cost of applying $T_{A^{(l)}}$ is lower.
\begin{equation*}
\begin{aligned}
    \tAl
    &\sim
    \tflop \times
        \left(
        2^{dl} k \frac{N}{2^{dl}}
        + \sum_{j=0}^{l} (6^d - 3^d) k^2 2^{dj}
        \right)
    \\
    &\sim \tflop \times \left( k N + (6^d - 3^d) k^2 2^{dl} \right)
\end{aligned}
\end{equation*}

The number of matrix-vector products to carry out sampling
for uniform basis matrices is also lower.
\begin{equation*}
    \tl
    \sim
    \left( \tmult + \tAl \right) \times (\colorsunif + \colorsnonunif) k
    +
    \tflop \times (6^d - 3^d) 2^{dl} k^2 \frac{N}{2^{dl}},
\end{equation*}

Summing $\tl$ over each level,
we find that the number of floating point operations
is lower than that of \cref{eq:h1complexity} by a factor of $\bigO(\log{N})$.
As in \cref{sec:h1_complexity}, we omit the costs associated with graph coloring.
\begin{multline*}
    \tcompress \sim \tmult \times \left( \colorsunif + \colorsnonunif \right) k \log{N} \\
    + \tflop \times k^2 N
        \left(
            \left( \colorsunif + \colorsnonunif + 6^d - 3^d \right) \log{N} \right. \\
            \left. + (6^d - 3^d) (\colorsunif + \colorsnonunif)
        \right)
\end{multline*}

\section{Numerical experiments} \label{sec:experiments}

In this section, we present a selection of numerical results.
The experiment in \cref{sec:low_dim}
demonstrates the ability of the graph coloring approach
to exploit low-dimensional structure.
In \cref{sec:bie,sec:operator_mult,sec:fmm,sec:nesteddiss},
we report the following quantities
for a number of test problems and rank structure formats:
(1) the time to compress the operator,
(2) the time to apply the compressed representation to a vector,
(3) the relative accuracy of the compressed representation,
and (4) the storage requirements of the compressed representation
measured as the number of floating point values per degree of freedom.
For compression time, we report both the total time taken for compression
as well as the ``net time,''
which does not include the time spent by the black-box multiplication routine.
The algorithms for compressing matrices and applying compressed representations
are written in Python,
and the black-box multiplication routines are written in MATLAB.
The experiments were carried out on a workstation
with two Intel Xeon Gold 6254 processors with 18 cores each
and 754GB of memory.

The matrices are compressed
with the $\h^1$, uniform $\h^1$, and  $\h^2$ formats.
For the uniform $\h^1$ format,
we take two approaches:
the ``$\h^1$\,+\,unif.'' approach uses the $\h^1$ sampling procedure
followed by a uniformization step (as in \cite{lin2011fast}),
and the ``unif. $\h^1$'' approach uses the technique
described in \cref{sec:unif_h1},
which uses a procedure that samples entire interaction lists
and thus avoids the uniformization step.

We measure the accuracy of the compressed matrices
using the relative error
\[
E = \frac{\|\mtx{A}_{\rm computed} - \mtx{A}\|}{\|\mtx{A}\|}
\]
computed via 20 iterations of the power method.
We also report the maximum leaf node size $m$
and the number $r$ of random vectors per test matrix,
which are inputs to the compression algorithm.

\subsection{Exploiting low-dimensional structure}
\label{sec:low_dim}

\begin{figure}[tb]
    \includegraphics[width=\textwidth]{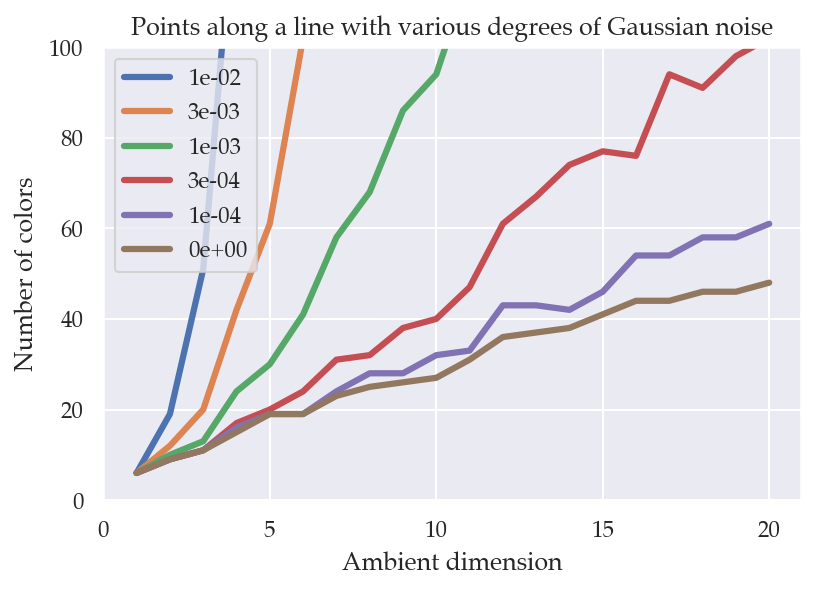}
    \caption{
        Number of colors for the incompatibility graphs
        that arise from sampling one level of admissible blocks of an $\h^1$ matrix
        based on a uniform grid along a randomly oriented line through
        $\lbrack 0, 1 \rbrack^d$ with added perturbation
        over a range of dimensions $d$.
    }
    \label{fig:line_with_noise}
\end{figure}

A major advantage of the graph coloring approach
is that it tailors the test matrices to the problem at hand,
exploiting low-dimensional structure to minimize the number of matrix-vector products.
In \cref{sec:h1,sec:unif_h1,sec:h2}, we establish upper bounds on the chromatic numbers of the graphs
in terms of the dimension $d$ of the computational domain.
However, if the geometry of the points exhibits lower-dimensional structure
(e.g., a discretization of a surface in three-dimensional space,
a machine learning dataset consisting of observations
belonging to a high-dimensional feature space),
the chromatic numbers may be much lower.

To demonstrate this effect, we color the graphs
that arise from sampling one level of admissible blocks of an $\h^1$ matrix
based on a uniform grid along a randomly oriented line through $\lbrack 0, 1 \rbrack^d$
over a range of dimensions $d$.
We perturb the data by adding Gaussian noise to the coordinates of the points
in order to simulate geometries that are not perfectly one-dimensional.
With zero noise, the points lie exactly on a line.
As more and more noise is added,
the effective dimensionality of the data increases
from one to the full ambient dimension.
The results reported in \cref{fig:line_with_noise}
demonstrate that the number of colors grows modestly
in the presence of low-dimensional structure
and very rapidly in the absence of low-dimensional structure.

\subsection{Boundary integral equation}
\label{sec:bie}
\begin{figure}
\begin{center}
\includegraphics[height=30mm]{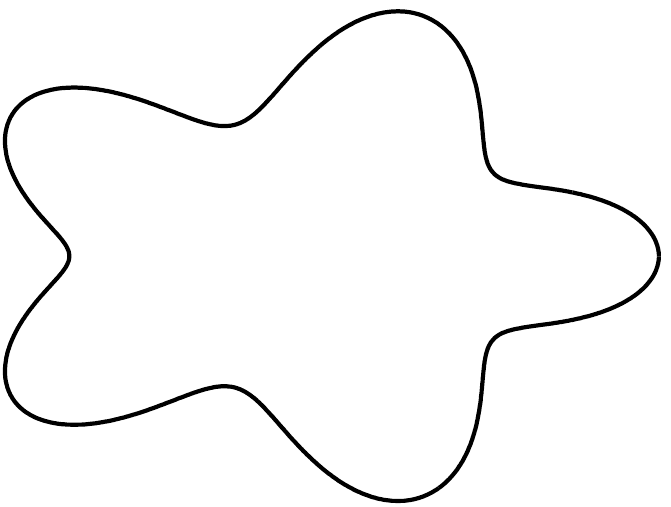}
\end{center}
\caption{Contour $\Gamma$ on which the BIE (\ref{eq:BIE}) is defined.}
\label{fig:BIE}
\end{figure}

We consider a matrix arising from the discretization of the 
Boundary Integral Equation (BIE)
\begin{equation}
\label{eq:BIE}
\frac{1}{2}q(\pvct{x}) +
\int_{\Gamma}
\frac{(\pvct{x} - \pvct{y})\cdot\pvct{n}(\pvct{y})}
     {4\pi|\pvct{x} - \pvct{y}|^{2}}\,q(\pvct{y})\,ds(\pvct{y}) = f(\pvct{x}),
\qquad\pvct{x} \in \Gamma,
\end{equation}
where $\Gamma$ is the simple closed contour in the plane shown in Figure \ref{fig:BIE}, 
and where
$\pvct{n}(\pvct{y})$ is the outwards pointing unit normal of $\Gamma$ at
$\pvct{y}$. The BIE (\ref{eq:BIE}) is a standard integral equation formulation
of the Laplace equation with boundary condition $f$ on the domain interior
to $\Gamma$. The BIE (\ref{eq:BIE}) is discretized using the
Nystr\"om method on $N$ equispaced points on $\Gamma$, with the Trapezoidal
rule as the quadrature (since the kernel in (\ref{eq:BIE}) is smooth, 
the Trapezoidal rule has exponential convergence). 

The fast matrix-vector multiplication is in this case furnished by the
recursive skeletonization (RS) procedure of \cite{2005_martinsson_fastdirect}.
To avoid spurious effects due to the rank structure inherent in RS, we 
compute the matrix-vector products
at close to double precision accuracy, and with an entirely
uncorrelated tree structure.

\begin{figure}[tb]
    \begin{subfigure}{0.48\textwidth}
        \includegraphics[width=\textwidth]{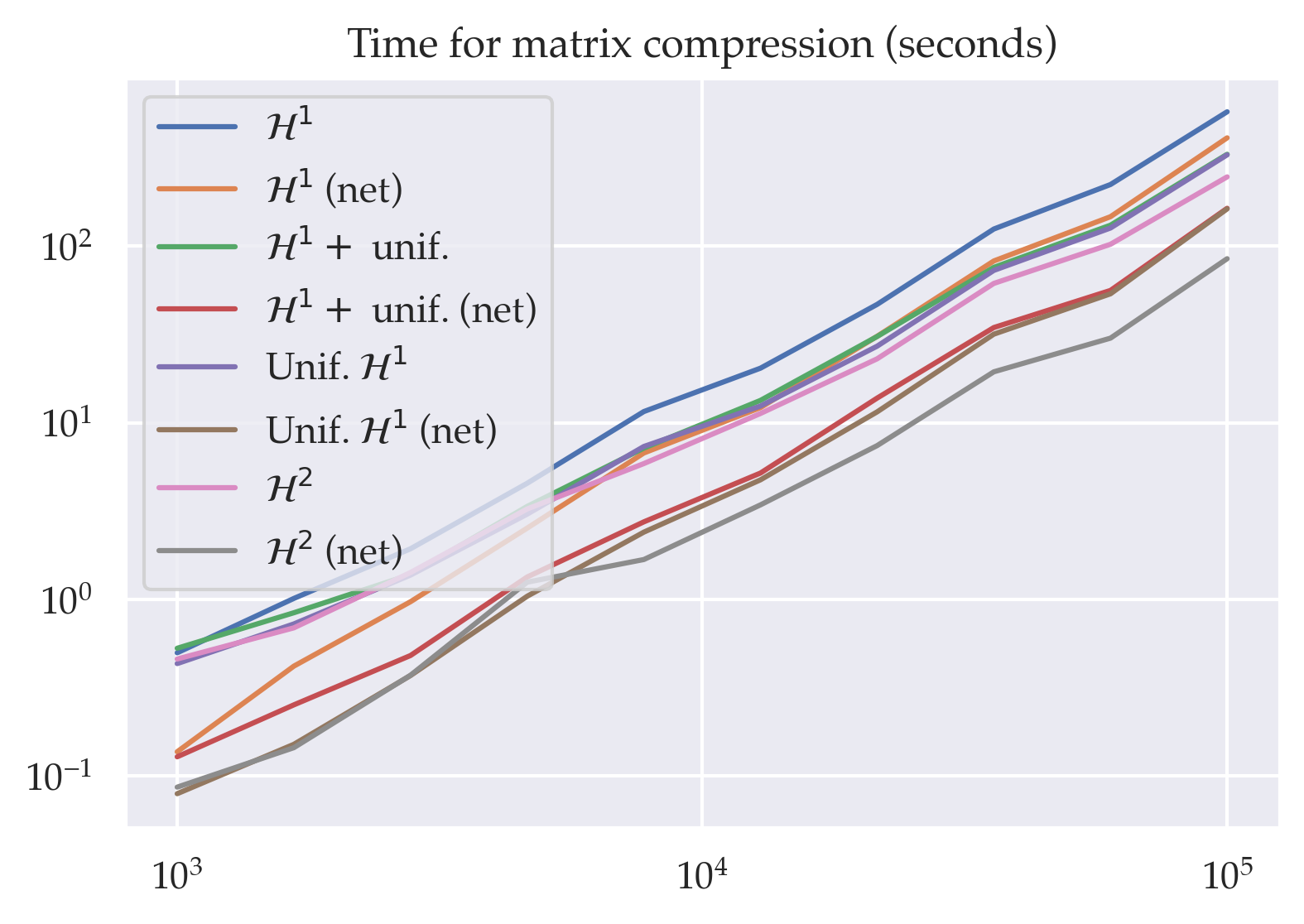}
    \end{subfigure}
    \hspace*{\fill}
    \begin{subfigure}{0.48\textwidth}
        \includegraphics[width=\textwidth]{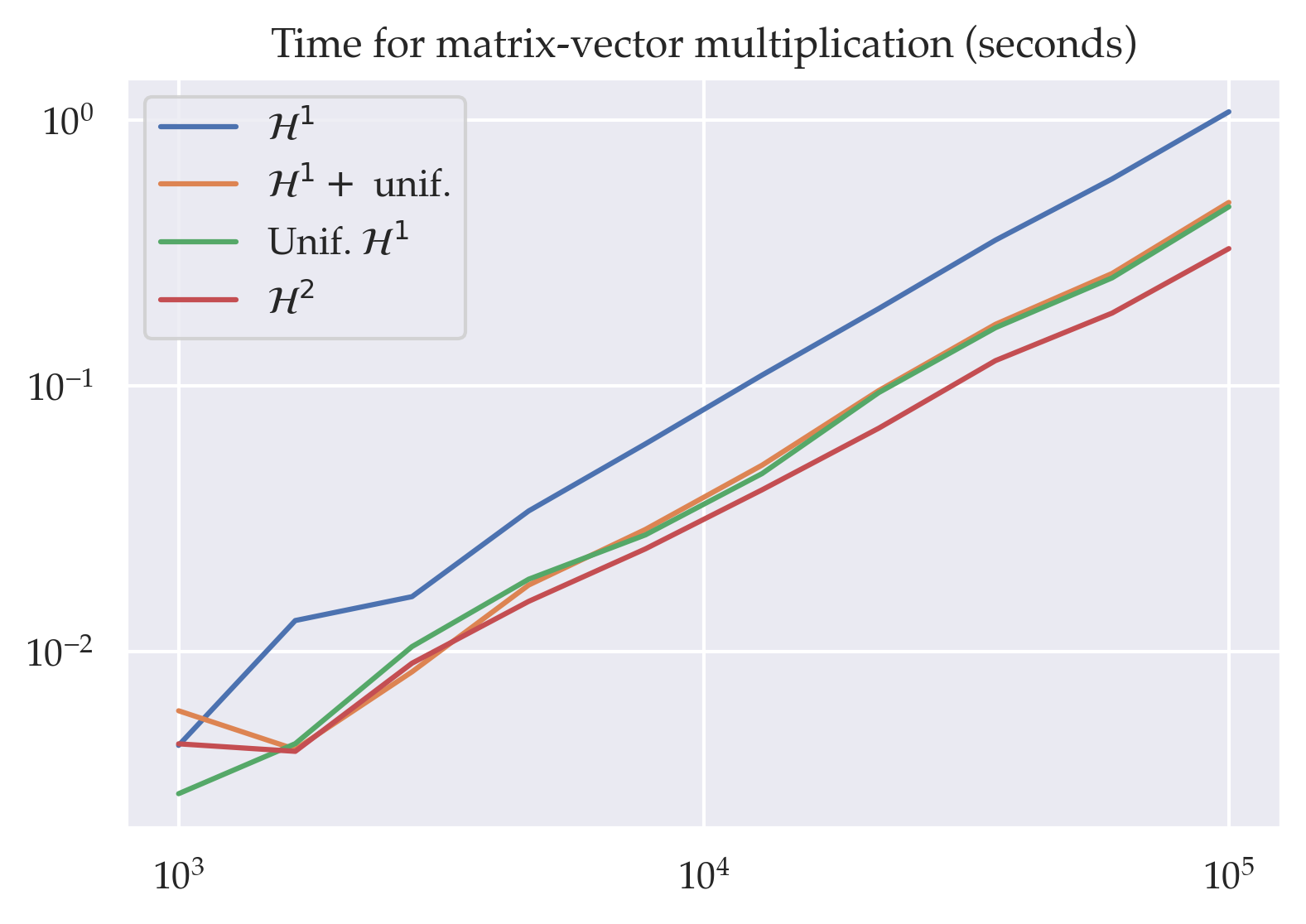}
    \end{subfigure}

    \medskip
    \begin{subfigure}{0.48\textwidth}
        \includegraphics[width=\textwidth]{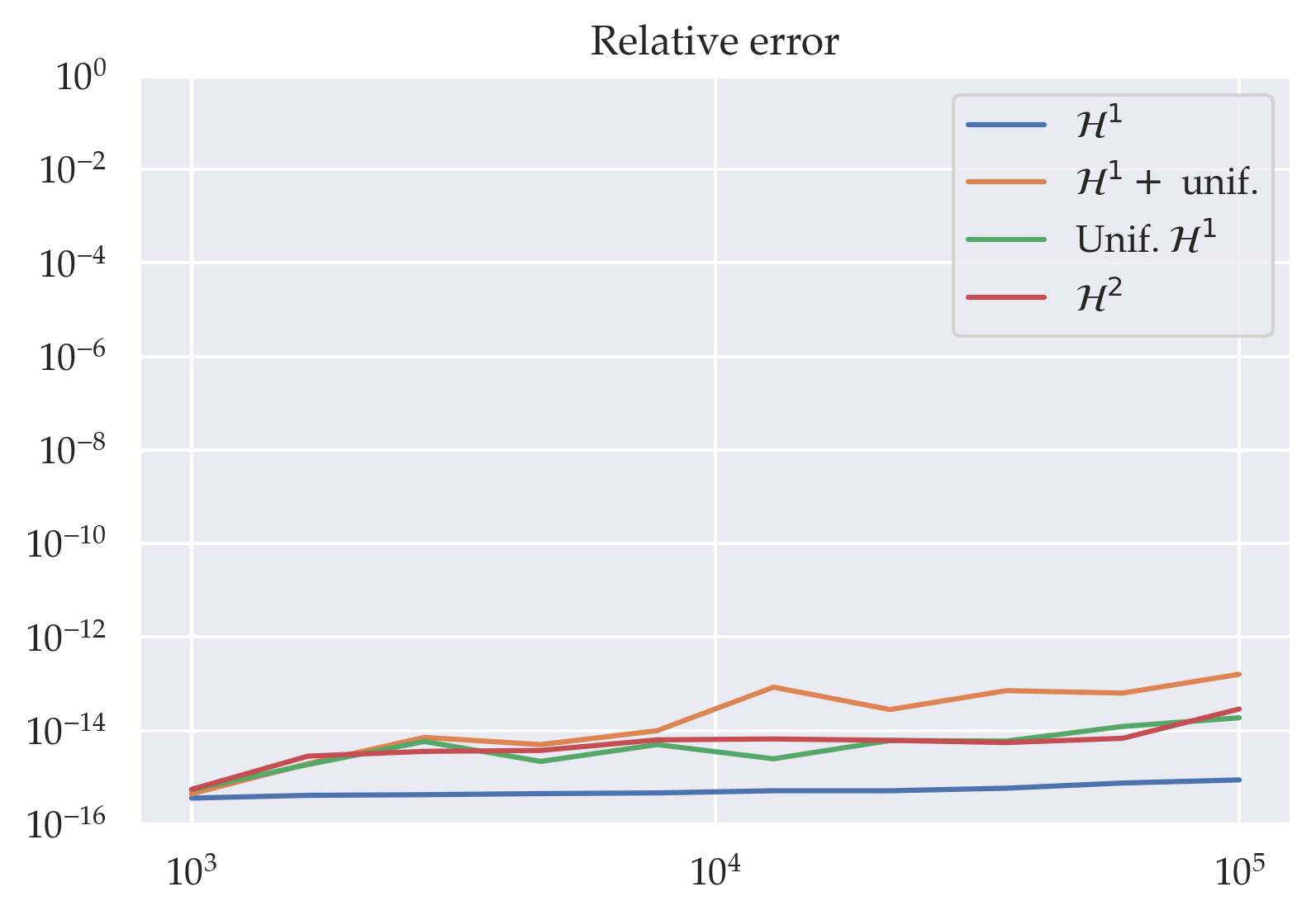}
    \end{subfigure}
    \hspace*{\fill}
    \begin{subfigure}{0.48\textwidth}
        \includegraphics[width=\textwidth]{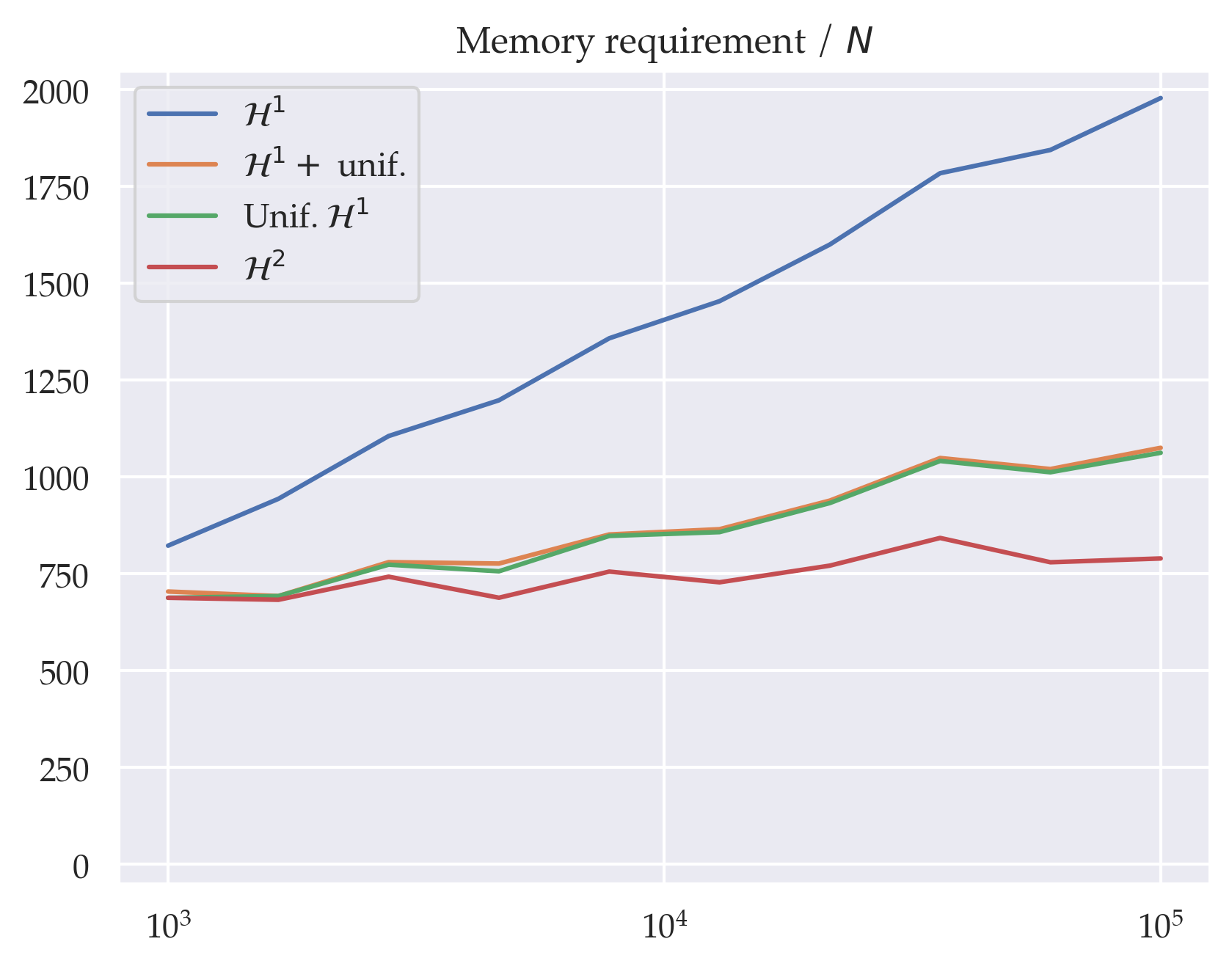}
    \end{subfigure}
    \caption{
        Results from applying peeling algorithms
        to a double layer potential on a simple contour in the plane.
        Here $r = 20$ and $m = 200$,
        and the horizontal axes represent the problem size $N$.
    }
    \label{fig:dl_results}
\end{figure}

Results are given in \cref{fig:dl_results}.

\begin{remark}
This problem under consideration here is included to illustrate the asymptotic performance of the proposed method. However, it is artificial in two ways: (1) There is no actual need to use more than a couple of hundred points to resolve (\ref{eq:BIE})
numerically to double precision accuracy. 
(2) In this case, we have alternate means of computing a rank structured representation of the matrix.
Subsequent examples illustrate more realistic use cases.
\end{remark}

\subsection{Operator multiplication}
\label{sec:operator_mult}

We next investigate how the proposed technique performs on a matrix matrix
multiplication problem. Specifically, we determine the Neumann-to-Dirichlet
operator $T$ for the contour shown in Figure \ref{fig:BIE} using the well
known formula
\[
T = S\left(\frac{1}{2}I + D^*\right)^{-1},
\]
where $S$ is the single layer operator 
$
[Sq](\pvct{x}) = \int_{\Gamma}-\frac{1}{2\pi}\,\log|\pvct{x} - \pvct{y}|\,q(\pvct{y})\,ds(\pvct{y}),
$
and where $D^{*}$ is the adjoint of the double-layer operator
$
[D^{*}q](\pvct{x}) = \int_{\Gamma}\frac{\pvct{n}(\pvct{x}) \cdot (\pvct{x} - \pvct{y})}{2\pi |\pvct{x} - \pvct{y}|^{2}}
\,q(\pvct{y})\,ds(\pvct{y}).
$
The operators $S$ and $D$ are again discretized using a Nystr\"om method
on equispaced points (with sixth order Kapur-Rokhlin \cite{1997_kapur_rokhlin}
corrections to handle the singularity in $S$), resulting in matrices
$\mtx{S}$ and $\mtx{D}$. The $\mtx{S}\bigl(\tfrac{1}{2}\mtx{I} + \mtx{D}^{*}\bigr)^{-1}$
is again applied using the recursive skeletonization procedure of 
\cite{2005_martinsson_fastdirect}.

\begin{figure}[tb]
    \begin{subfigure}{0.48\textwidth}
        \includegraphics[width=\textwidth]{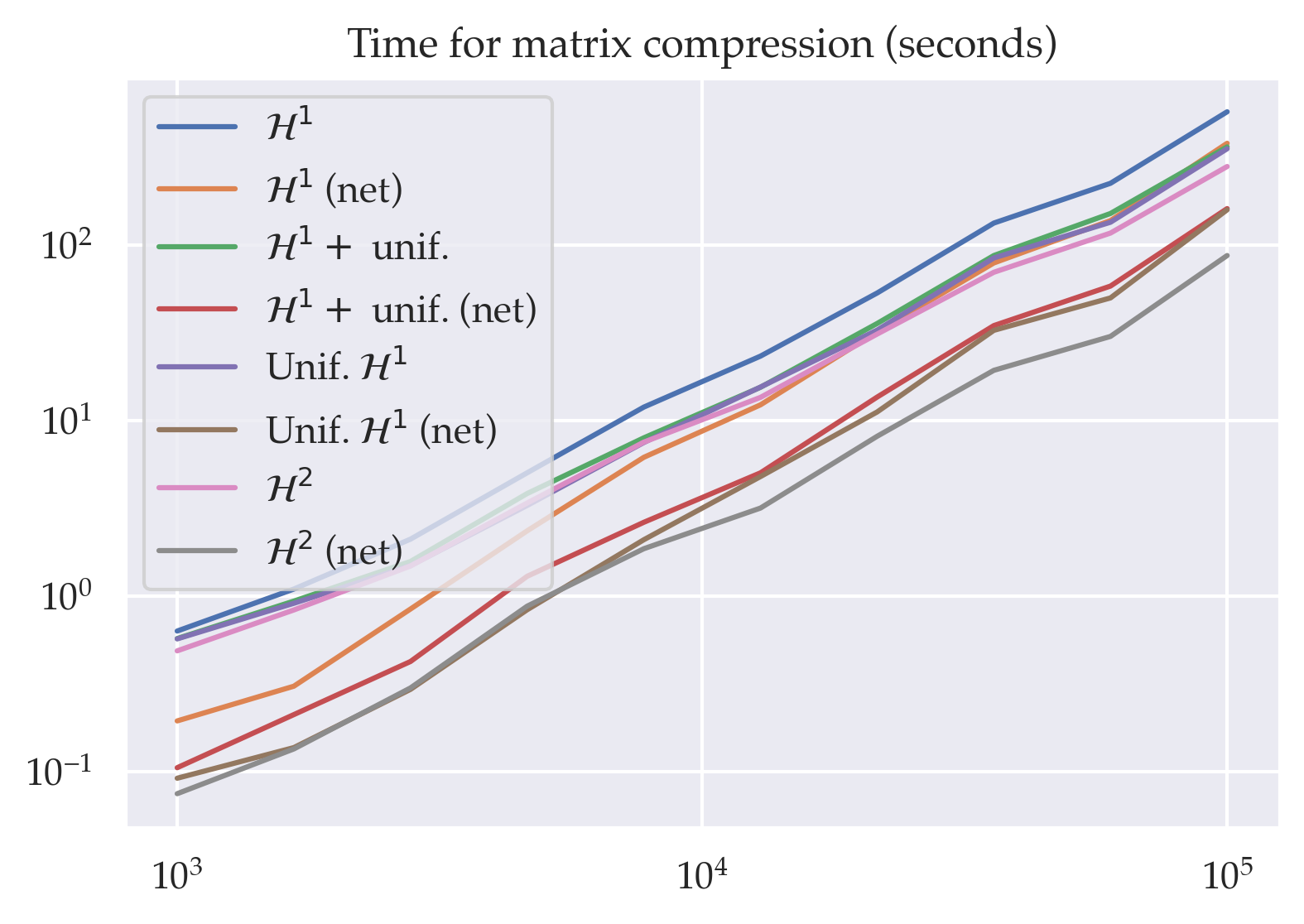}
    \end{subfigure}
    \hspace*{\fill}
    \begin{subfigure}{0.48\textwidth}
        \includegraphics[width=\textwidth]{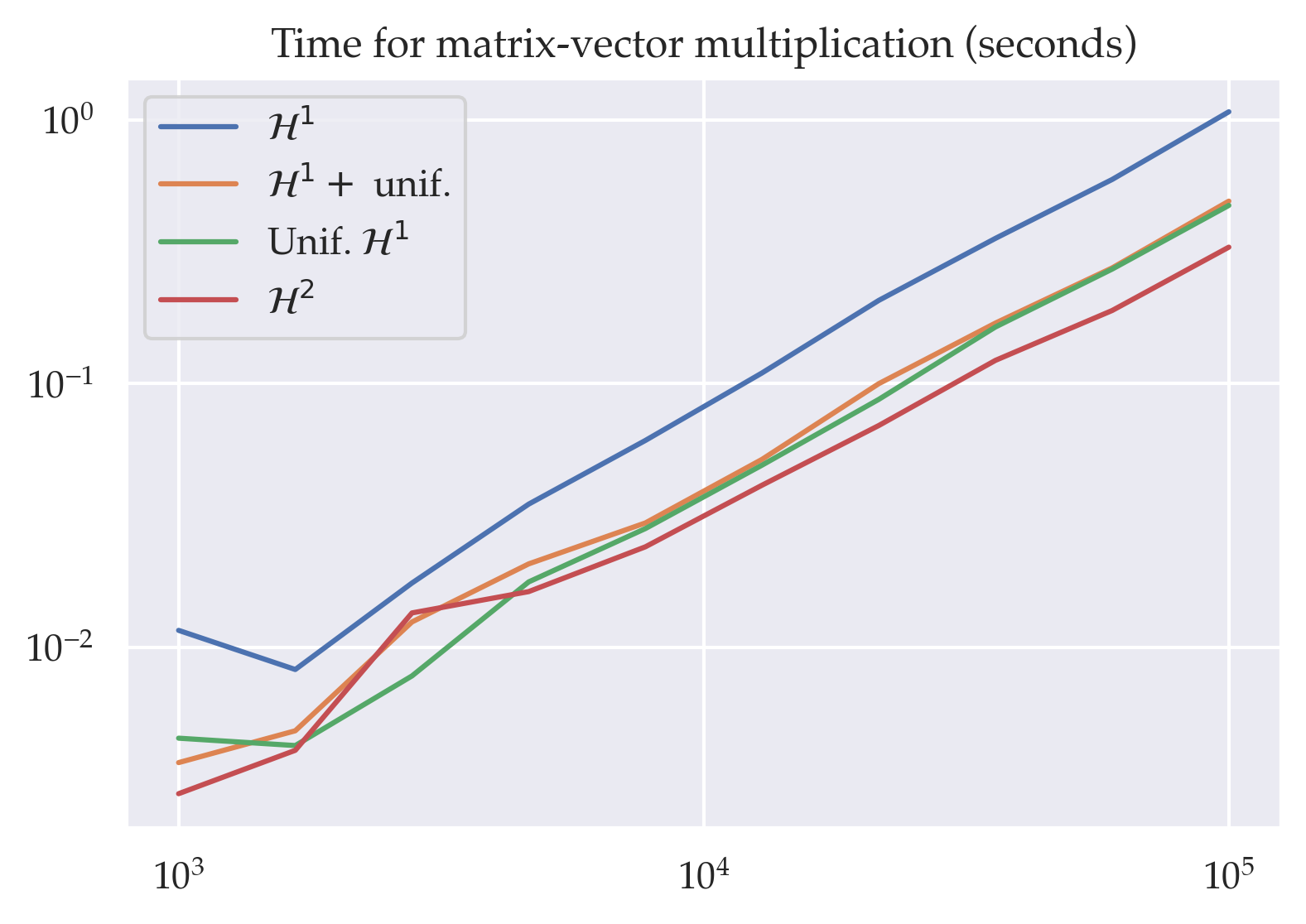}
    \end{subfigure}

    \medskip
    \begin{subfigure}{0.48\textwidth}
        \includegraphics[width=\textwidth]{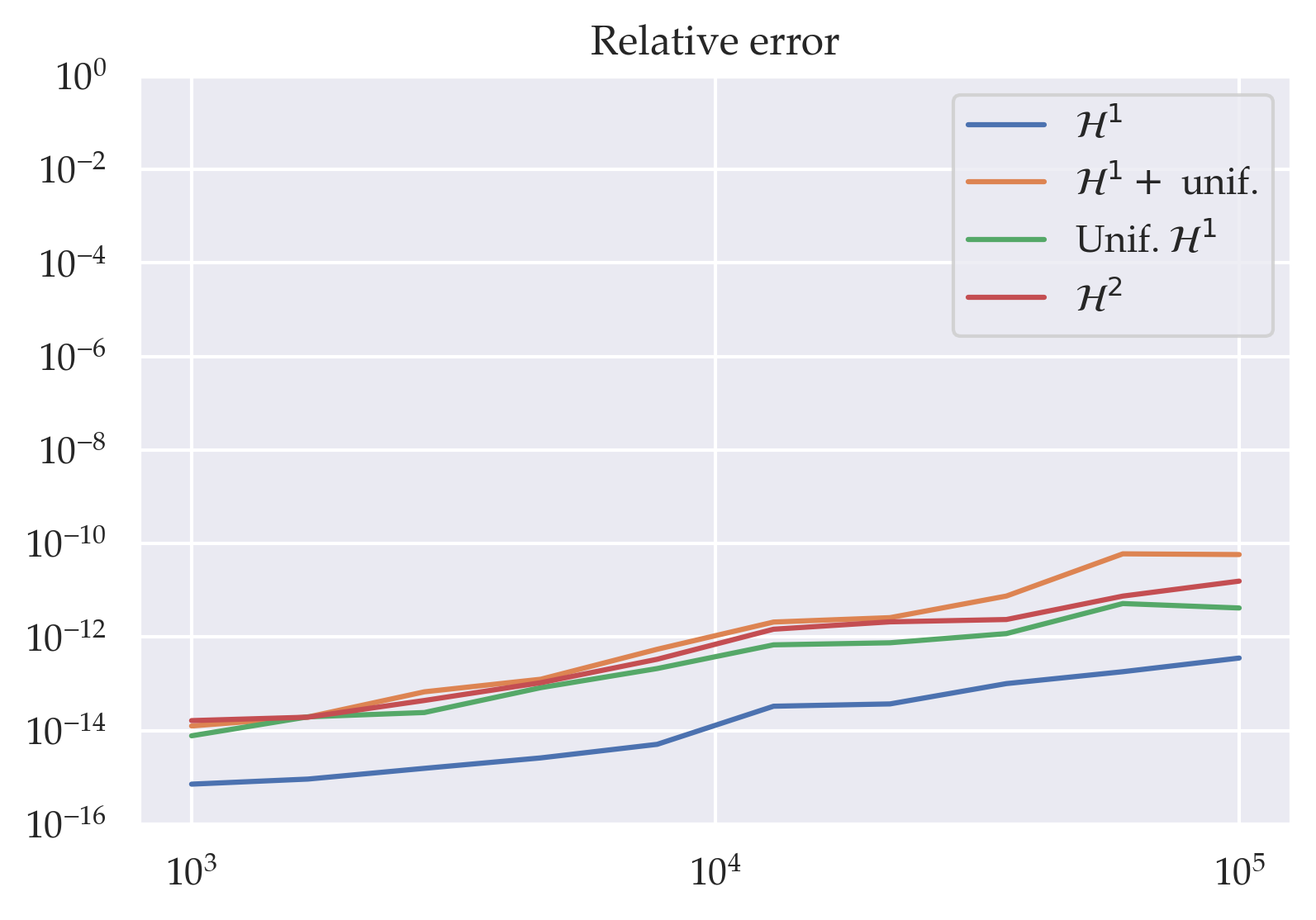}
    \end{subfigure}
    \hspace*{\fill}
    \begin{subfigure}{0.48\textwidth}
        \includegraphics[width=\textwidth]{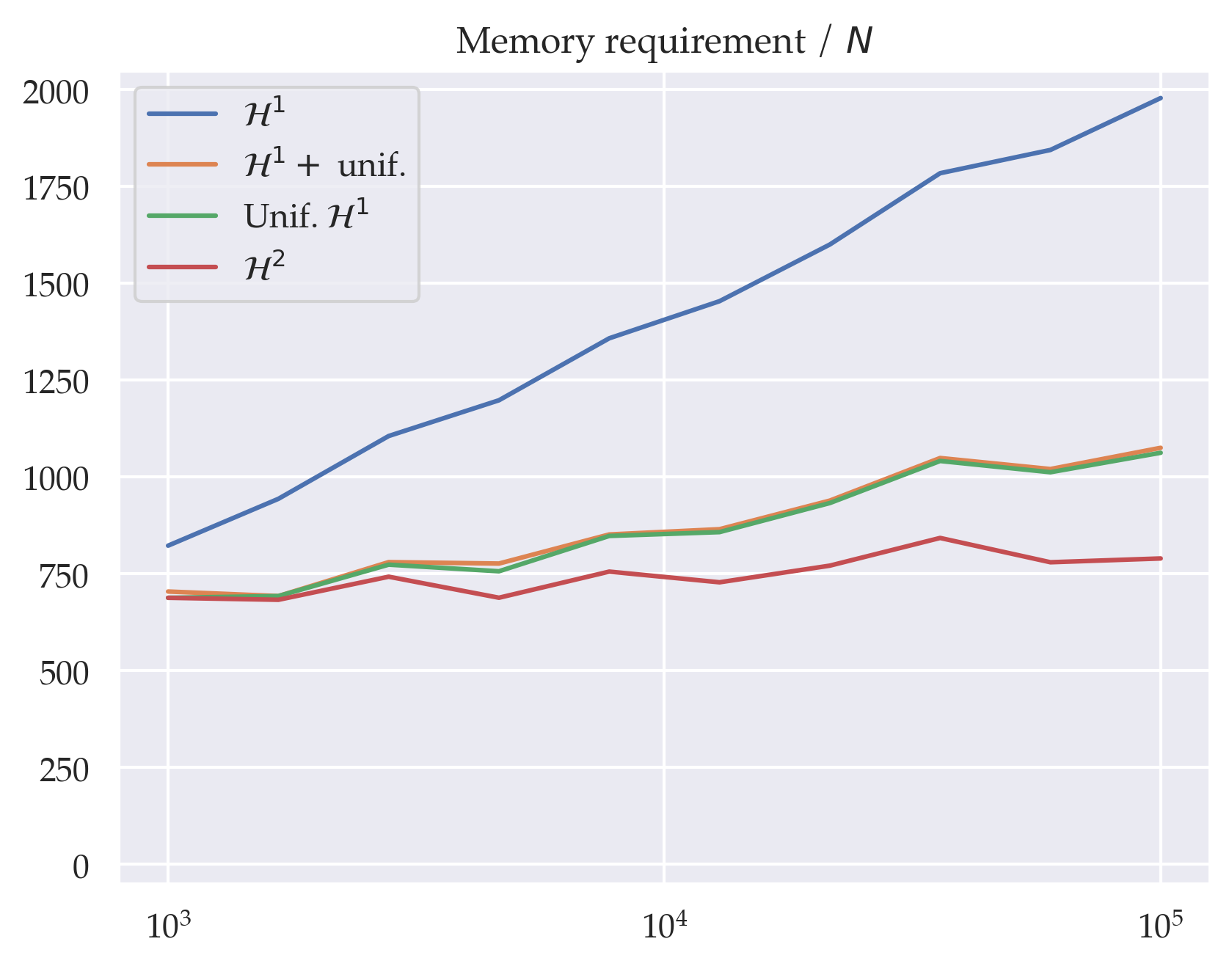}
    \end{subfigure}
    \caption{
        Results from applying peeling algorithms
        to the Neumann-to-Dirichlet operator.
        Here $r = 20$ and $m = 200$,
        and the horizontal axes represent the problem size $N$.
    }
    \label{fig:dfn_results}
\end{figure}

Results are given in \cref{fig:dfn_results}.

\subsection{Fast multipole method}
\label{sec:fmm}

We consider a kernel matrix
representing $N$-body Laplace interactions in three dimensions,
where the interaction kernel is defined by
\[
\kfunc(x, y) = \sum_{i \neq j} \frac{c_j}{\Vert x_i - x_j\Vert}
\]
for sets of $N$ points $\{x_i\}$ and charges $\{c_i\}$.
To simulate a problem with low-dimensional structure,
we distribute the points uniformly at random
on the surface of the unit sphere.
We use an implementation of the fast multipole method
included in the Flatiron Institute Fast Multipole Libraries~\cite{fmm3d}
to efficiently apply the matrix to vectors.
The operator is computed to a target accuracy of $\varepsilon = 10^{-4}$.

\begin{figure}[tb]
    \begin{subfigure}{0.48\textwidth}
        \includegraphics[width=\textwidth]{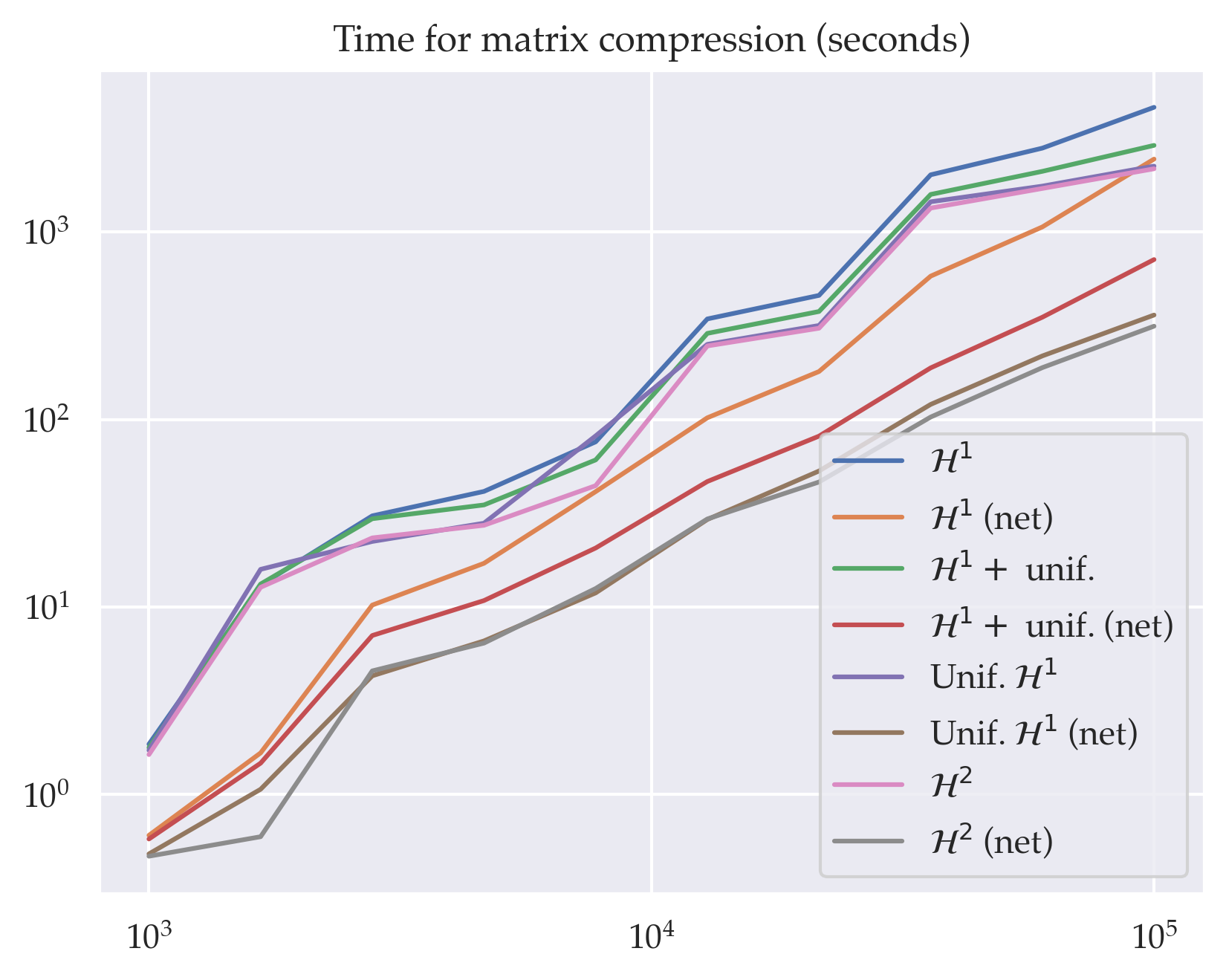}
    \end{subfigure}
    \hspace*{\fill}
    \begin{subfigure}{0.48\textwidth}
        \includegraphics[width=\textwidth]{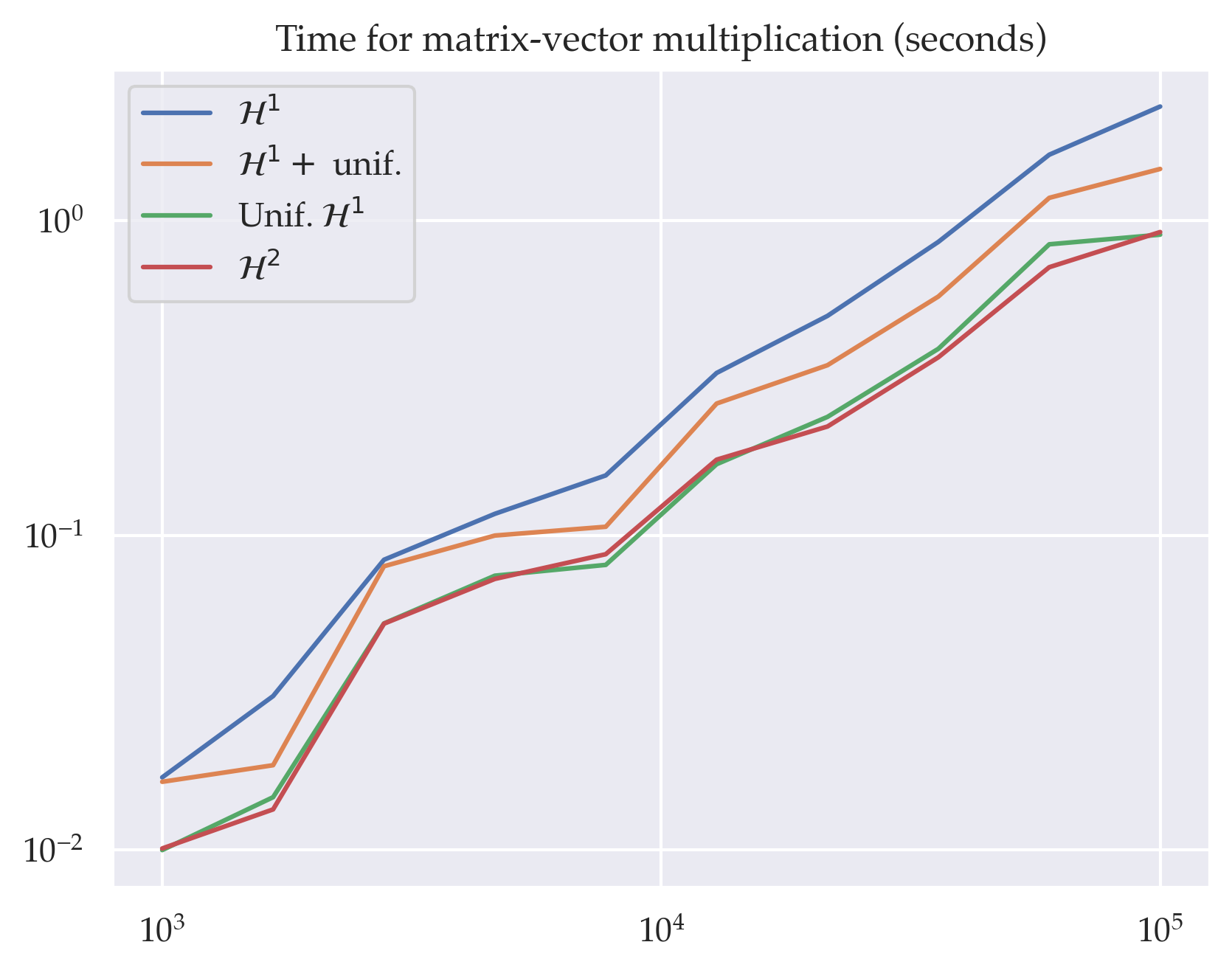}
    \end{subfigure}

    \medskip
    \begin{subfigure}{0.48\textwidth}
        \includegraphics[width=\textwidth]{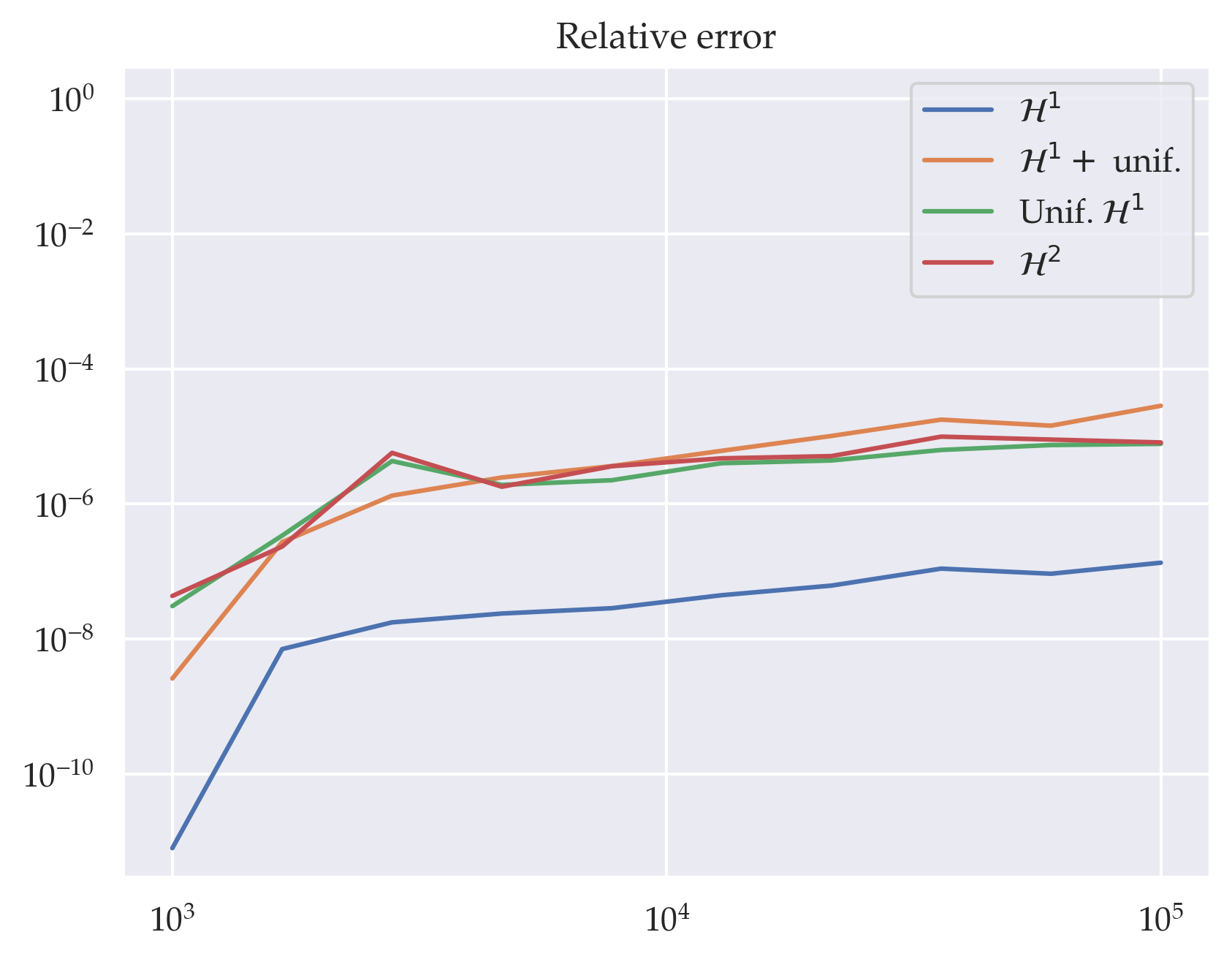}
    \end{subfigure}
    \hspace*{\fill}
    \begin{subfigure}{0.48\textwidth}
        \includegraphics[width=\textwidth]{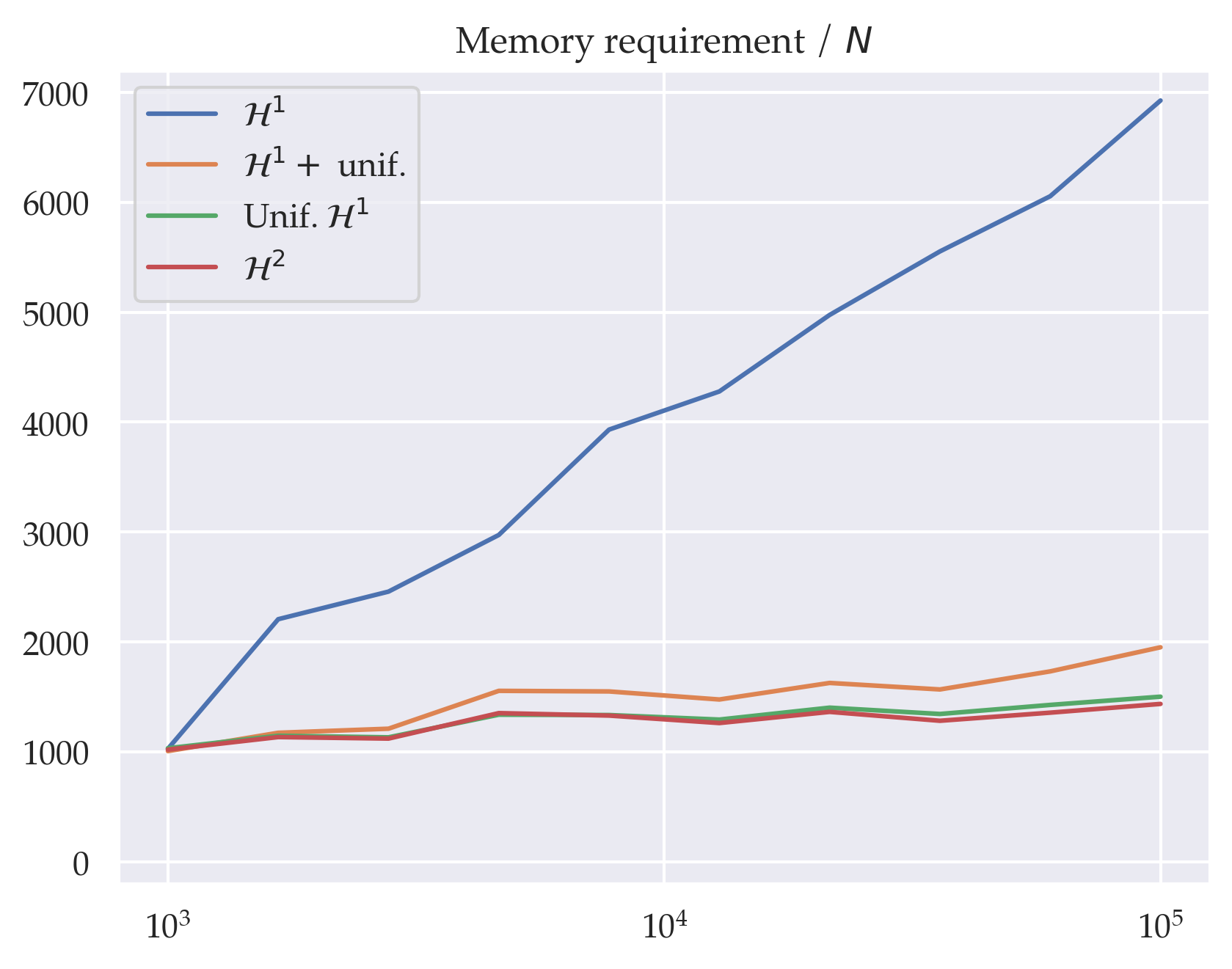}
    \end{subfigure}
    \caption{
        Results from applying peeling algorithms
        to the 3D fast multipole method operator.
        Here $r = 20$ and $m = 50$,
        and the horizontal axes represent the problem size $N$.
    }
    \label{fig:fmm3d_results}
\end{figure}

Results are given in \cref{fig:fmm3d_results}.

\subsection{Frontal matrices in nested dissection}
\label{sec:nesteddiss}
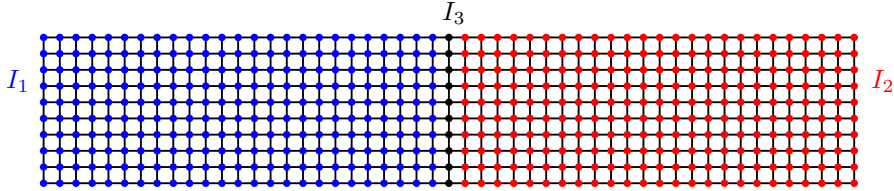
\begin{figure}
\begin{center}
\setlength{\unitlength}{1mm}
\begin{picture}(125,26)
\put(005,00){\includegraphics[width=110mm]{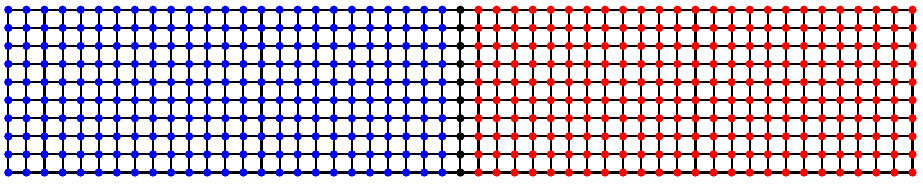}}
\put(001,14){$\color{blue}I_{1}$}
\put(116,14){$\color{red}I_{2}$}
\put(059,23){$\color{black}I_{3}$}
\end{picture}
\end{center}
\caption{An example of the grid in the sparse LU example
described in Section \ref{sec:nesteddiss}. 
There are $N \times n$ points in the grid, shown for $N = 10, n = 51$.
}
\label{fig:nesteddiss_geom}
\end{figure}

Our next example is a simple model problem that illustrates the behavior
of the proposed method in the context of sparse direct solvers. The idea
here is to use rank structure to compress the increasingly large Schur
complements that arise in the LU factorization of a sparse matrix arising
from the finite element or finite difference discretization of an
elliptic PDE, cf.~\cite[Ch.~21]{2019_martinsson_fast_direct_solvers}.
As a model problem, we consider
a $51N \times 51N$ matrix
$\mtx{C}$ that
encodes the stiffness matrix for the standard five-point stencil
finite difference approximation to the Poisson equation on a rectangle
using a grid with $N \times 51$ nodes. We partition the grid
into three sets $\{1,2,3\}$, as shown in \cref{fig:nesteddiss_geom},
and then tessellate $\mtx{C}$ accordingly,
\[
\mtx{C} =
\begin{bmatrix}
\mtx{C}_{11} & \mtx{0}      & \mtx{C}_{13} \\
\mtx{0}      & \mtx{C}_{22} & \mtx{C}_{23} \\
\mtx{C}_{31} & \mtx{C}_{32} & \mtx{C}_{33}
\end{bmatrix},
\]
where
$\mtx{C}_{11}$
and
$\mtx{C}_{22}$
are $25N \times 25N$,
and
$\mtx{C}_{33}$
is $N \times N$.
The matrix we seek to compress is the $N \times N$ Schur complement
\[
\mtx{A} = \mtx{C}_{33}
- \mtx{C}_{31}\mtx{C}_{11}^{-1}\mtx{C}_{13}
- \mtx{C}_{32}\mtx{C}_{22}^{-1}\mtx{C}_{23}.
\]
In our example, we apply $\mtx{A}$ to a vector by calling standard
sparse direct solvers for the left and the right subdomains, respectively.
We use sparse matrix routines that are provided by MATLAB
and that rely on UMFPACK~\cite{davis2004algorithm} for the sparse solves.

\begin{figure}[tb]
    \begin{subfigure}{0.48\textwidth}
        \includegraphics[width=\textwidth]{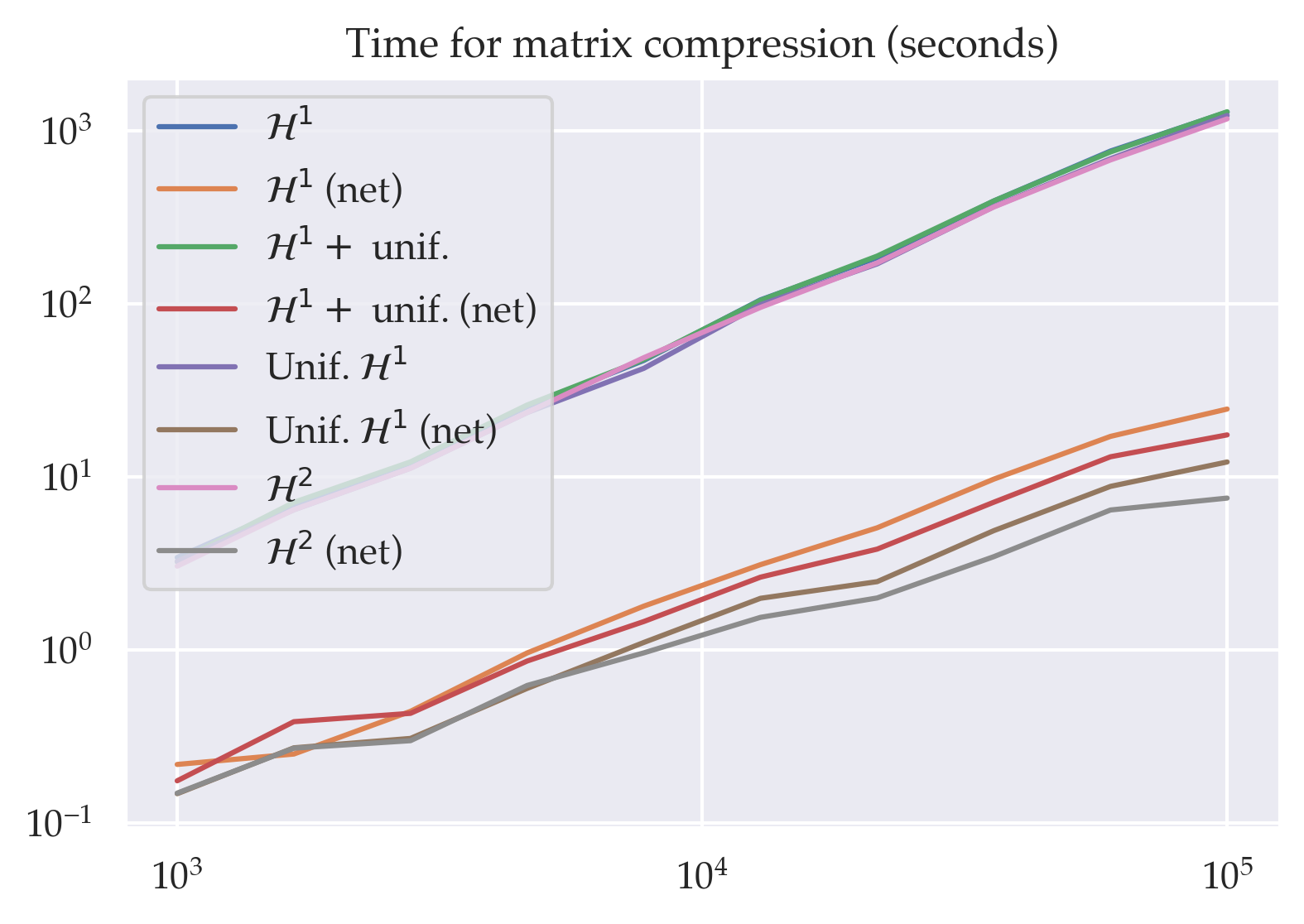}
    \end{subfigure}
    \hspace*{\fill}
    \begin{subfigure}{0.48\textwidth}
        \includegraphics[width=\textwidth]{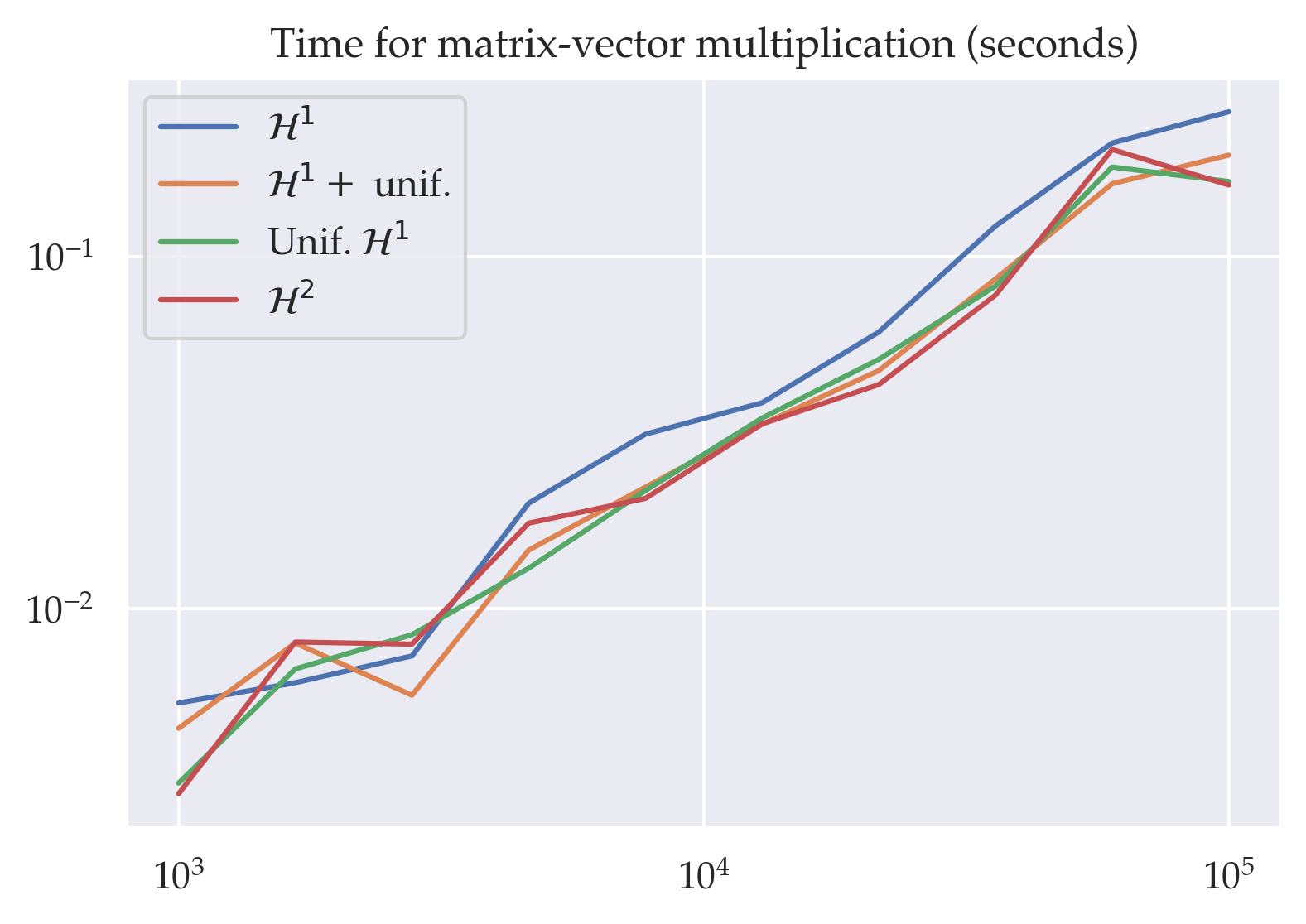}
    \end{subfigure}

    \medskip
    \begin{subfigure}{0.48\textwidth}
        \includegraphics[width=\textwidth]{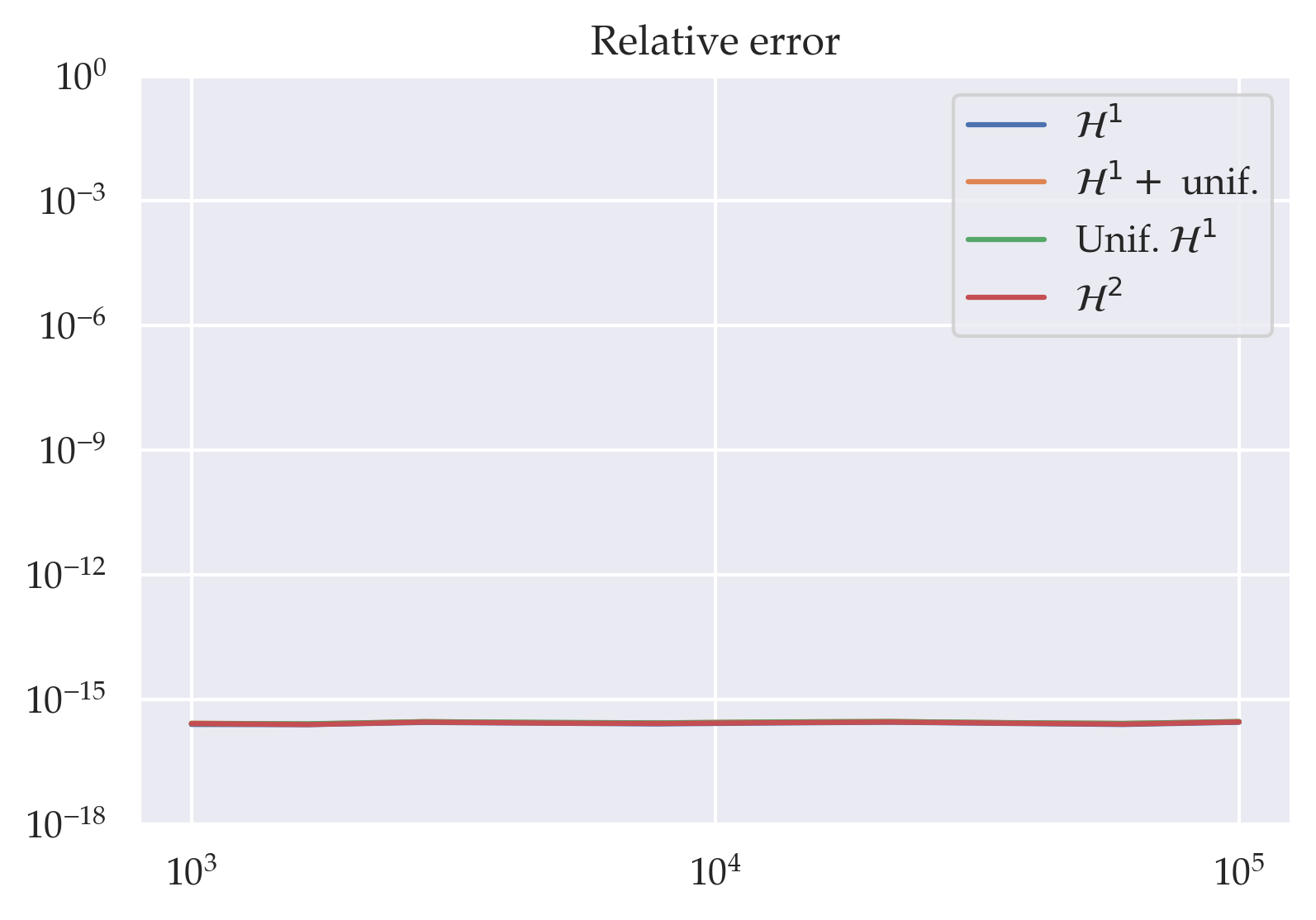}
    \end{subfigure}
    \hspace*{\fill}
    \begin{subfigure}{0.48\textwidth}
        \includegraphics[width=\textwidth]{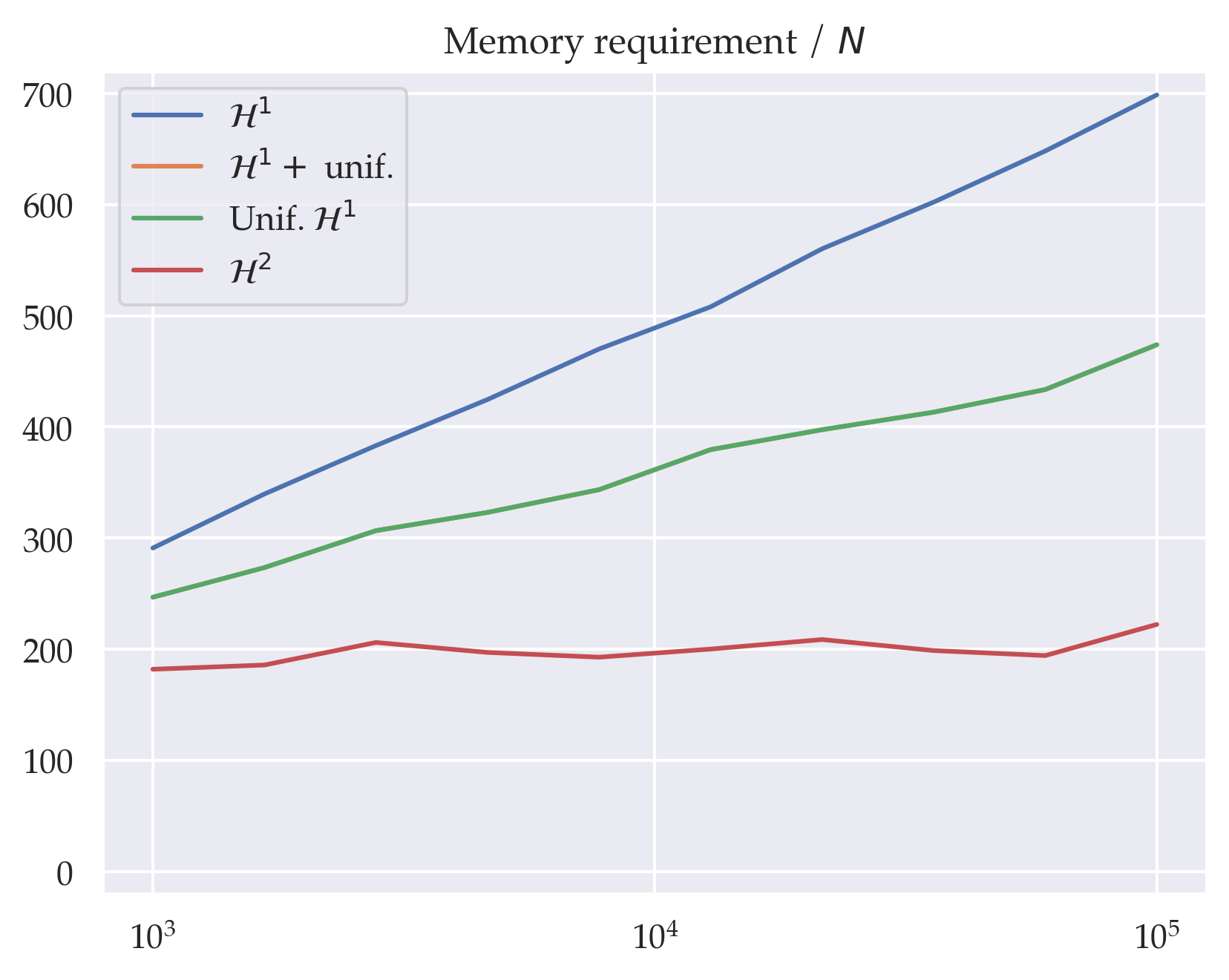}
    \end{subfigure}
    \caption{
        Results from applying peeling algorithms
        to frontal matrices in the nested dissection algorithm.
        Here $r = 10$ and $m = 50$,
        and the horizontal axes represent the problem size $N$.
    }
    \label{fig:nd_results}
\end{figure}

Results are given in \cref{fig:nd_results}.

\subsection{Summary of observations}
\begin{itemize}
\item The numerical results support our claims of quasilinear
    computational complexity for all steps of the algorithms presented.
\item The unif. $\h^1$ compression scheme consistently outperforms
    the $\h^1$\,+\,unif. scheme.
    The advantage appears not only in shorter compression times,
    but also in higher accuracy and lower storage requirements.
\item The approximations achieve high accuracy in every case,
    with the exception of the 3D FMM, for which the accuracy
    matches or exceeds the accuracy of the operator.
\item The count of floating point numbers per degree of freedom
    remains constant or grows very slowly for the $\h^2$ format,
    reflecting linear asymptotic storage cost,
    and it grows logarithmically for the other formats,
    reflecting quasilinear asymptotic storage costs.
\item In every case, the majority of the compression time
    is spent in the black-box matrix-vector multiplication routines,
    highlighting the importance of minimizing the number of matrix-vector products.
\end{itemize}

\section{Conclusions} \label{sec:conclusions}
This paper presents algorithms for randomized compression
of rank-structured matrices.
The algorithms only access the matrix
via black-box matrix-vector multiplication routines.
We formulate a graph coloring problem
to design sets of test matrices
that are tailored to the given matrix
and to minimize the number of matrix-vector multiplications required.
Numerical experiments indicate
that the algorithms are accurate 
and much more efficient than prior works,
particularly when the underlying geometry
exhibits low-dimensional structure.

\subsection*{Acknowledgments}
The work reported was supported by the Office of Naval Research (N00014-18-1-2354),
by the National Science Foundation (DMS-1952735 and DMS-2012606),
and by the Department of Energy ASCR (DE-SC0022251).

\bibliographystyle{siamplain}
\bibliography{references}
\end{document}